\documentclass[10pt]{amsart}
\usepackage{graphicx,color}

\usepackage{url} 
\newtheorem{theorem}{Theorem}[section]
\newtheorem{lemma}[theorem]{Lemma}

\theoremstyle{definition}
\newtheorem{definition}[theorem]{Definition}

\theoremstyle{remark}
\newtheorem{remark}[theorem]{Remark}

\numberwithin{equation}{section}

\subjclass[2000]{Primary~81U40, Secondary~47A40}
\keywords{Schr{\"o}dinger operator, lattice, scattering theory}

\title[Schr{\"o}dinger Operators on Perturbed Lattices]
{Spectral Properties of Schr{\"o}dinger Operators on Perturbed Lattices}
\author{Kazunori ANDO}
\address{Department of Mathematics, Inha University, Incheon, 402-751, Republic of  Korea}
\author{Hiroshi ISOZAKI}
\address{Institute of Mathematics,
University of Tsukuba, 
Tsukuba, 305-8571, Japan}
\author{Hisashi MORIOKA}
\address{Center for Promotion of Educational Innovation, Shibaura Institute of Technology, 307 Fukasaku, Minuma, Saitama-shi, Saitama, 337-8570, Japan}
\address{
ando@inha.ac.kr, isozakih@math.tsukuba.ac.jp, hisashimorioka@gmail.com}

\date{April, 3, 2015}

\begin{document}
\maketitle

\begin{abstract}
We study the spectral properties of Schr{\"o}dinger operators on perturbed lattices. We shall prove the non-existence or the discreteness of embedded eigenvalues, the limiting absorption principle for the resolvent, construct a spectral representation, and define the S-matrix. Our theory covers the square, triangular,  diamond, Kagome lattices, as well as the ladder, the graphite and the subdivision of square lattice.
\end{abstract}

\section{Introduction}

In this and the forthcoming articles, we shall investigate
 the spectral properties of Schr{\"o}dinger  operators on perturbed periodic lattices of dimension $d\ge2$. 
The physical background is the scattering phenomenon. Sending waves from the infinity of a perturbed periodic structure, we observe the behavior of scattered waves at infinity.
The mapping from the free wave in the remote past to the scattered wave in the remote future is the scattering matrix (S-matrix). Our goal is the inverse scattering, 
 in particular, we  aim at the reconstruction of the perturbed periodic structure from the S-matrix of a fixed energy. In the present paper, we devote ourselves to the forward problem, i.e. that
 of the continuous spectrum of the Schr{\"o}dinger operator describing the scattering process,  more precisely, the limiting absorption principle for the resolvent,  construction of spectral representations and S-matrix are the main ingredients. The inverse problem will be discussed in the 2nd part \cite{AndIsoMor}.
 
We start from the Laplacian $\widehat H_0$ on a lattice in ${\bf R}^d$, which is a matrix whose entries are shift operators. Passing to the Fourier series, it is transformed to an operator of multiplication by an hermitian matrix $H_0(x)$ acting on a vector bundle on the flat torus ${\bf T}^d = {\bf R}^d/(2\pi{\bf Z})^d$. Then  its spectral properties
boil down to the characteristic polynomial $p(x,\lambda) = \det (H_0(x) - \lambda)$.

A lot of remarkable methods have been  found in the long history of scattering theory, e.g. \cite{Ka59}, \cite{Eidus}, \cite{Ag70}, \cite{KatoKuroda}, \cite{Kuroda}, \cite{IkeSa},  \cite{Enss}, \cite{Mourre}, \cite{DeGe}, \cite{Yaf}, \cite{MazMe}, \cite{Me}, ranging over a variety of areas of mathematical analysis. It is worthwhile to recall here the role of Fourier transform in the study of the continuous spectrum of a self-adjoint differential operator 
$H = P(D_x) + V(x,D_x)$ on ${\bf R}^d$. The  first important step is the 
Rellich type theorem, or Rellich-Vekua theorem, which gives the minimal growth rate for solutions to the equation
$(P(D_x) + V(x,D_x) - \lambda)u = 0$ as $|x|\to\infty$ (see \cite{Re43}, \cite{Vek43}). When $V(x,D_x)=0$, by passing to the Fourier transform, this is reduced to the algebraic properties of the polynomial $P(\xi)$ and the Paley-Wiener Theorem (see \cite{Lit66}, \cite{Lit70}, \cite{Hor73}, \cite{Mur76}). 
 To see 
 the continuous spectrum of $H$, Agmon \cite{Ag75} derived a theorem of division by $P(\xi)$ in some weighted $L^2({\bf R}^d)$-space and proved the existence of the limit $\lim_{\epsilon\to0}(H - \lambda \mp i\epsilon)^{-1}$. This is the limiting absorption principle, the key step to clarify 
 the detailed spectral strucure. Agmon-H{\"o}rmander \cite{AgHo76}, \cite{HoVol1}, \cite{HoVol2} supplemented this approach by introducing Besov spaces $\mathcal B, \mathcal B^{\ast}$, which are optimal for the existence of the limit $(H - \lambda \mp i0)^{-1}$, and introduced the radiation condition in the form of pseudo-differential operators ($\Psi$DO) to guarantee the uniqueness of the solution.
 Our  strategy is to extend this Agmon-H{\"o}rmander's approach to discrete problems.

Inverse potential scattering for discrete Schr{\"o}dinger operators 
have already been considered in \cite{IsKo} on the square lattice (see also \cite{Es})
and  in \cite{Ando13} on the hexagonal lattice, where knowledge 
of the S-matrix for all energies is used  for the reconstruction of the potential by the method of complex Born approximation.
In \cite{IsMo2}, inverse potential scattering on the square lattice 
from the S-matrix of one fixed energy is studied. The main idea of \cite{IsMo2} is to
reduce the issue to an inverse boundary value problem on a bounded domain,
and the reconstruction is done through the Dirichlet-Neumann map. 
To relate the scattering matrix with the D-N map for the bounded domain, 
 the discrete analogue of
the Rellich type uniqueness theorem plays an  important role (\cite{IsMo1}). 
The radiation  condition is closely related to the asymptotic behavior of the Green function as space variables tend to infinity. It requires the strict convexity of the Fermi surface $M_{\lambda} = \{x \in {\bf T}^d\, ; \,p(x,\lambda)=0\}$ for the Laplacian on the unperturbed periodic lattice, hence restricts the energy region for the reconstruction procedure to be valid. In the present work, we introduce the radiation condition in the form of wavefront set as in \cite{AgHo76}, and instead of the behavior of the resolvent near infinity of the lattice space, we consider the singularity expansion of the resolvent on the torus. This change of view point makes it possible to remove the above mentioned assumption of strict convexity 
and the restriction for the energy level. We also mention that  \cite{IsMo1} relies on  \cite{Sha} which gives  basic ideas from the theory of functions of several complex variables and algebraic geometry in the proof of Rellich
 type theorem.

The plan of the paper is as follows. Our final result in this paper is 
Theorem \ref{S7FinalTheorem} in \S 7 
on the characterization of the solution space of the Helmholtz equation.

\smallskip
\hskip10mm\S 2.  Basic properties of graph \\
\hskip14mm\S 3. Examples of periodic lattices \\
\hskip14mm\S 4. Distributions on the torus \\
\hskip14mm\S 5. Rellich type theorem on the torus \\
\hskip14mm\S 6. Spectral properties of the free Hamiltonian \\
\hskip14mm\S 7. Hamiltonians on the perturbed lattice

\smallskip

In \S 2, we introduce the Laplacian on the lattice. Passing to the Fourier series, it is transferred to a matrix on the torus, whose characteristic polynomial is crucial for the spectral properties,  
and we pick up two typical cases. \S 3 is devoted to the exposition of various examples of lattices, and basic properties of their Laplacians. The main analytical tool to study the spectrum is a division theorem for distributions on the torus, which we discuss in \S 4. The main aim of \S 5 is to prove the Rellich
 type theorem for the discrete Laplacian on the lattice by using Hilbert Nullstellensatz. In \S 6, we study the spectral properties for the unperturbed lattice. Based on these preparations, we develop in \S 7 
the spectral and scattering theory for the Laplacian on the perturbed lattice. The main results are the resolvent estimates, spectral representations, unitarity of the S-matrix and 
the structure of the solution space for the Helmholtz equation.

The spectral properties for discrete Schr{\"o}dinger operators,
 more generally perturbed Laplacians on non-compact graphs, have been discussed a long time with an abundance of references,  and are nowadays becoming  more active issues. We cite here  rather recent articles, \cite{KoShiSu98}, \cite{S99}, \cite{HiShi04}, \cite{HiNo09}, \cite{HSSS12},  \cite{Su13}, 
\cite{BuIvKu13}, \cite{ColFra13}, \cite{KoSa13}, \cite{Na14}, which are directly related to this paper. 

\medskip
Let us give some remarks on the notation. 
${\bf T}^d_{\bf C}$ denotes the complex torus 
\begin{equation}
{\bf T}^d_{\bf C} = {\bf C}^d/(2\pi{\bf Z})^d.
\label{S1TdC}
\end{equation}
 For $f \in \mathcal S'({\bf R}^d)$, 
$\widetilde f(\xi)$ denotes its Fourier transform 
\begin{equation}
\widetilde f(\xi) = (2\pi)^{-d/2}\int_{{\bf R}^d}e^{-ix\cdot\xi}
f(x)dx.
\label{S1Fouriertransf}
\end{equation}
On the other hand, for $f(x) \in \mathcal S'({\bf T}^d)$, $\widehat f(n)$ denotes its Fourier coefficients
\begin{equation}
\widehat f(n) = (2\pi)^{-d/2}\int_{{\bf T}^d}e^{-ix\cdot n}f(x)dx.
\label{S1Fouriercoeffi}
\end{equation}
We also use $\widehat f = \big(\widehat f(n)\big)_{n\in{\bf Z}^d}$ to denote a function on ${\bf Z}^d$, and by $\mathcal U$ the mapping
\begin{equation}
\mathcal U : {\mathcal S}'({\bf Z}^d) \ni \big(\widehat f(n)\big)_{n\in{\bf Z}^d} \to 
f(x) = (2\pi)^{-d/2}\sum_{n\in{\bf Z}^d}\widehat f(n)e^{in\cdot x} \in 
{\mathcal S}'({\bf T}^d).
\label{S1Fourierseries}
\end{equation}
For Banach spaces $X$ and $Y$, ${\bf B}(X;Y)$ denotes the set of all bounded operators from $X$ to $Y$. For a self-adjoint operator $A$, $\sigma(A), \sigma_p(A), \sigma_d(A), \sigma_e(A)$  denote its spectrum, point spectrum, discrete spectrum 
 and essential spectrum, respectively. $\mathcal H_{ac}(A)$ is the absolutely continuous subspace for $A$, and $\mathcal H_{p}(A)$ is the closure of the linear hull of eigenvectors of $A$.
For an interval $I \subset {\bf R}$ and a Hilbert space ${\bf h}$, $L^2(I,{\bf h},d\lambda)$ denotes the set of all ${\bf  h}$-valued $L^2$-functions on $I$ with respect to the measure $d\lambda$. $S^m_{1,0}$ denotes the standard H{\"o}rmander class of symbols for pseudo-differential operators ($\Psi$DO), i.e. $|\partial_x^{\alpha}\partial_{\xi}^{\beta}p(x,\xi)| \leq C_{\alpha\beta}(1 + |\xi|)^{m-\beta}$ (see e.g. \cite{HoVol3}).


\section{Basic properties of graph} \label{Basic}


\subsection{Vertices and edges}
We consider an infinite, connected  graph $\{\mathcal V, \mathcal E\}$, where $\mathcal V$ is a vertex set and $\mathcal E$ an edge set.
We assume that the graph is simple, i.e. there are neither self-loop, which is an edge connecting a vertex to itself, nor multiple edges, which are two or more edges connecting the same vertices.
 For an edge $e = (v,w)\in \mathcal E$, we denote
$$
o(e) = v, \quad t(e) = w, \quad \overline{e} = (w,v), \quad
v\sim w.
$$
In the following, we assume that $e \in \mathcal E \Longrightarrow \overline{e} \in \mathcal E$.
We put
\begin{equation}
{\mathcal N}_v = \{w \in \mathcal V \, ; \, v \sim w\},
\label{C1S1Adjacentpoints}
\end{equation}
which will be called the set of points {\it adjacent} to $v$.
The {\it degree} of $v \in \mathcal V$ is then defined by
\begin{equation}
{\rm deg}\, (v) = {\sharp}\, {\mathcal N}_v = {\sharp}\,\{e \in \mathcal E\, ; \, o(e) = v\},
\nonumber
\end{equation}
which is assumed to be finite for all $v \in \mathcal V$.
A function $f : {\mathcal V} \to {\bf C}$ is denoted as $f = \left(f(v)\right)_{v\in\mathcal V}$. 
Let $\ell^2(\mathcal V)$ be the set of ${\bf C}$-valued functions $f$ on $\mathcal V$ satisfying
\begin{equation}
\|f\|_{\rm deg}^2 := {\mathop\sum_{v\in\mathcal V}}|f(v)|^2 \,{\rm deg}\,(v) < \infty.
\nonumber
\end{equation}
Equipped with the inner product
\begin{equation}
(f,g)_{\rm deg} = {\mathop\sum_{v\in\mathcal V}}f(v)\overline{g(v)}\,{\rm deg}\,(v),
\nonumber
\end{equation}
$\ell^2(\mathcal V)$ is a Hilbert space. 


\subsection{Laplacian on the periodic graph}\label{subsectionperiodicgraph}
A periodic graph in ${\bf R}^d$ is a triple $\Gamma_0 = \{\mathcal L_0, \mathcal V_0, \mathcal E_0\}$, where  $\mathcal L_0$ is a lattice of rank $d \geq 2$ in ${\bf R}^d$ with basis ${\bf v}_j, j= 1,\cdots,d$, i.e. 
\begin{equation}
\mathcal L_0 = \big\{{\bf v}(n)\, ; \, n \in 
{\bf Z}^d\big\}, \quad
{\bf v}(n) = \sum_{j=1}^dn_j{\bf v}_j, \quad n =(n_1,\cdots,n_d) \in {\bf Z}^d,
\nonumber
\end{equation}
and the vertex set is defined by
\begin{equation}
\mathcal V_0 = {\mathop\cup_{j=1}^s}\big(p_j + \mathcal L_0\big),
\nonumber
\end{equation}
and where $p_j$, $j = 1,\cdots,s$, are the points in ${\bf R}^d$ satisfying
\begin{equation}
p_i - p_j \not\in \mathcal L_0, \quad {\rm if}\quad i\neq j.
\label{C1pi-pj}
\end{equation}
By (\ref{C1pi-pj}), there exists a bijection $\mathcal V_0 \ni a \to (j(a),n(a)) 
\in \{1,\cdots,s\}\times{\bf Z}^d$ such that
\begin{equation}
a = p_{j(a)} + {\bf v}(n(a)).
\nonumber
\end{equation}
The group ${\bf Z}^d$ acts on ${\mathcal V}_0$ as follows :
\begin{equation}
{\bf Z}^d\times{\mathcal V}_0 \ni (m,a) \to m\cdot a := p_{j(a)}+{\bf v}(m+n(a)) \in {\mathcal V}_0.
\nonumber
\end{equation}
The edge set  $\mathcal E_0 \subset \mathcal V_0 \times \mathcal V_0$ is assumed to satisfy
\begin{equation}
\mathcal E_0 \ni (a,b) \Longrightarrow  
(m\cdot a, m\cdot b) \in \mathcal E_0,
\quad \forall m \in {\bf Z}^d.
\nonumber
\end{equation}
Then ${\rm deg}\,(p_j + {\bf v}(n))$ depends only on $j$, and is denoted by ${\rm deg}_0(j)$ :
\begin{equation}
{\rm deg}_0(j) ={\rm  deg}\,(p_j + {\bf v}(n)).
\label{S2Definedeg(j)}
\end{equation}
Any function $\widehat f$ on $\mathcal V_0$ is written as
$\widehat f(n) = (\widehat f_1(n),\cdots,\widehat f_s(n)), \  n \in {\bf Z}^d$, where $\widehat f_j(n)$ is identified with a function on $p_j + \mathcal L_0$.  Hence $\ell^2(\mathcal V_0)$ is the Hilbert space equipped with the inner product
$$
(\widehat f,\widehat g)_{\ell^2(\mathcal V_0)} = 
\sum_{j=1}^s(\widehat f_j,\widehat g_j)_{{\rm deg}_0(j)}.
$$
We then define a unitary operator $\mathcal U_{\mathcal L_0} : \ell^2(\mathcal V_0) \to L^2({\bf T}^d)^s$
\begin{equation}
\big(\mathcal U_{\mathcal L_0}\widehat f\big)_j = (2\pi)^{-d/2}
\sqrt{{\rm deg}_0(j)}\sum_{n\in {\bf Z}^d}\widehat f_j(n)e^{in\cdot x},
\label{S1UdDefine}
\end{equation}
where $L^2({\bf T}^d)^s$ is equipped with the inner product
\begin{equation}
(f,g)_{L^2({\bf T}^d)^s} = \sum_{j=1}^s\int_{{\bf T}^d}f_j(x)\overline{g_j(x)}dx.
\label{S2L2bfTdinnerproduct}
\end{equation}

Recall that the shift operator $\widehat S_j$ acts on a sequence $\big(a(n)\big)_{n\in{\bf Z}^d}$ as follows :
$$
\big(\widehat S_ja\big)(n) = a(n + {\bf e}_j),
$$
where ${\bf e}_1 = (1,0,\cdots,0), \cdots, {\bf e}_d= (0,\cdots,0,1)$. Then we have 
\begin{equation}
\mathcal U \widehat S_j = e^{-ix_j}\mathcal U.
\label{S2Ud0SjUd0}
\end{equation}
The Laplacian $\widehat \Delta_{\Gamma_0}$ on the graph $\Gamma_0 = 
\{\mathcal L_0, \mathcal V_0, \mathcal E_0\}$ is defined by the following formula
\begin{equation}
\begin{split}
  (\widehat\Delta_{\Gamma_0} \widehat f)(n) & = (\widehat g_1(n),\cdots,\widehat g_s(n)),\\
\widehat g_i(n) & = \frac{1}{{\rm deg}_0(i)}\sum_{b\sim p_i + {\bf v}(n)}\widehat f_{j(b)}(n(b)),
\end{split}
\label{S2DefineLaplacian}
\end{equation}
where $b = p_{j(b)} + {\bf v}(n(b))$. 
Passing to the Fourier series, we rewrite it into the following form :
\begin{equation}
\mathcal U_{\mathcal L_0}(- \widehat\Delta_{\Gamma_0}) (\mathcal U_{\mathcal L_0})^{-1} f = H_0(x)f(x), \quad 
f \in L^2({\bf T}^d)^s,
\nonumber
\end{equation}
where $H_0(x)$ is an $s\times s$ Hermitian matrix whose entries are trigonometric functions. Let $D$ be the $s\times s$ diagonal matrix whose $(j,j)$ entry is $\sqrt{{\rm deg}_0(j)}$. Then $\mathcal U_{\mathcal L_0} = D\mathcal U$, hence
\begin{equation}
H_0(x) = DH_0^0(x)D^{-1}, \quad H_0^0(x) = \mathcal U(- \widehat\Delta_{\Gamma_0}){\mathcal U}^{-1},
\label{S2H0(x)rewritten}
\end{equation}
and $H_0^0(x)$ is computed by (\ref{S2Ud0SjUd0}).


\subsection{Preliminary facts}
In this subsection, we consider geometric and algebraic properties of the following functions $a_d(z)$ and $b_d(z)$ : 
\begin{equation}
a_d(z) = \sum_{j=1}^d\cos z_j,
\label{C1S4Definead(z)}
\end{equation}
\begin{equation}
b_d(z) = \sum_{j=1}^d\cos z_j + \sum_{1\leq j<k \leq d}\cos(z_j - z_k),
\label{C1S4Definebd(z)}
\end{equation}
since all the characteristic polynomials of the examples of lattices to be presented in the next section are reduced to them. 

For an analytic function $f$ on ${\bf T}^d_{\bf C}$, we put
\begin{equation}
\begin{split}
&S_{a}^{\bf C}(f) = \left\{z \in {\bf T}^d_{\bf C}\, ; \, f(z)=a\right\}, \\
& S_{a,reg}^{\bf C}(f) = \left\{z \in S_{a}^{\bf C}(f) ; \, \nabla_z f(z)\neq 0\right\}, \\
& S_{a,sng}^{\bf C}(f) = \{z \in S_{a}^{\bf C}(f)\, ; \, \nabla_zf(z)=0\}, \\
& SV(f) = \left\{f(z)\, ; \, z \in \bigcup_{a \in \bf C}S_{a,sng}^{\bf C}(f)\right\},\\
& f({\bf T}^d) = \left\{f(x)\, ; \, x \in {\bf T}^d\right\}.
\end{split}
\nonumber
\end{equation}


\begin{lemma}\label{LemmaSquareLattice}
(1) $\bigcup_{a \in \bf C}S_{a,sng}^{\bf C}(a_d) = \left(\pi{\bf Z}\right)^d\cap {\bf T}^d_{\bf C}$. \\
\noindent
(2) $SV(a_d)  = \{-d, - d + 2, \cdots, d-2, d\}$. \\
\noindent
(3) $a_d({\bf T}^d) = [-d,d]$. \\
\noindent
(4)  For $-d < a < d$, each connected component of $S_{a,reg}^{\bf C}(a_d)$ intersects with ${\bf T}^d$, and the intersection is a $(d-1)$-dimensional real analytic submanifold of ${\bf T}^d$.
\end{lemma}

Proof. Since $\partial a_d/\partial z_j= - \sin z_j$, the assertions (1), (2) and (3) are easy to prove.
Let ${\bf C}^{\ast} = {\bf C}\setminus\{(-\infty,-1]\cup[1,\infty)\}$, and recall that 
$\cos\zeta$ maps $\{0<{\rm Re}\,\zeta < \pi\}$ conformally to ${\bf C}^{\ast}$. 
To prove (4), we take $z^{(0)} = (z^{(0)}_1,\cdots,z^{(0)}_d) \in S_{a,reg}^{\bf C}(a_d)$ arbitrarily.
Then, we can construct continuous curves $c_j^{\ast}(t), (0\leq t\leq1),j = 1,\cdots,d-1,$ such that 
$c_j^{\ast}(0)=\cos z_j^{(0)}$, $c_j^{\ast}(t) \in {\bf C}^* $ for $0 < t < 1$, and $c_j^{\ast}(1) \in (0,1)$, moreover 
$$
c_d^{\ast}(t) :=  a- \sum_{j=1}^{d-1}c_j^{\ast}(t) \in {\bf C}^{\ast}, \quad 0 < t < 1.
$$
Putting $c_j(t) = \arccos c_j^{\ast}(t)$ and $c(t)= (c_1(t),\cdots,c_d(t))$ for $0 \leq t \leq 1$, we have $c(t) \in S_{a,reg}^{{\bf C}} (a_d)$ and $ c(1) \in S_{a,reg}^{{\bf C}} (a_d )\cap {\bf T}^d $.
This proves that $c(t) $ is the desired curve. \qed


\begin{lemma}\label{LemmaDiamond}
(1) For $d =$ even, 
\begin{equation}
SV(b_d) =  
 \Big\{\frac{(\ell +1)^2}{2} - \frac{d+1}{2}\, ; \, \ell = -d, -d+2, \cdots, d-2, d\Big\} \cup \big\{- \frac{d+1}{2}\big\},
\nonumber
\end{equation}
and for $d=$ odd,
\begin{equation}
SV(b_d)= \Big\{\frac{(\ell +1)^2}{2} - \frac{d+1}{2}\, ; \, \ell = -d, -d+2, \cdots, d-2, d\Big\}.
\nonumber
\end{equation}
(2) For  $a \neq -(d+1)/2$, $S^{\bf C}_{a,sng}(b_d) \subset 
\left(\pi{\bf Z}\right)^d\cap {\bf T}^d_{\bf C}$. \\
\noindent
(3) For $a = -(d+1)/2$, $S_{a}^{\bf C}(b_d)$ is a union of analytic submanifolds of complex dimension $d-1$, $d-2$ and a discrete set. If the  discrete set appears, it is in ${\bf T}^d$.
In particular, $S_{-(d+1)/2}^{\bf C}(b_d)\cap{\bf T}^d$ is a union of real analytic submanifolds of real dimension $\leq d-2$.   \\
\noindent
(4) \ 
$b_d({\bf T}^d) = \left[-(d+1)/2,d(d+1)/2\right].$ \\
\noindent
(5)  Assume that $ - (d+1)/2 < a \leq d(d+1)/2$, and let
\begin{equation}
\begin{split}
&\widetilde S_{a,j}(b_d) = \{z \in S_{a}^{\bf C}(b_d)\, ; \, 1 + \sum_{k\neq j}e^{-iz_k} \neq 0\},\quad 1 \leq j \leq d,\\
& \widetilde S_{a,0}(b_d) = \{z \in S_{a}^{\bf C}(b_d) \, ; \, 1 + \sum_{k\neq j}e^{-iz_k}=0,\, \ \forall j\}.
\end{split}
\nonumber
\end{equation}

(i) For any $z^{(0)} \in \widetilde S_{a,0}(b_d)$, there exist $j \neq 0$, $z^{(j)} \in \widetilde S_{a,j}(b_d)$  and an $S_{a,reg}^{\bf C}(b_d)$-valued continuous curve $c(t), 0 \leq t \leq 1,$ such that $c(0) = z^{(0)}$ and $c(1)= z^{(j)}$. 

(ii)  For any $j \neq 0$, $\widetilde S_{a,j}(b_d)$ is arcwise connected and
$\widetilde S_{a,j}(b_d)\cap{\bf T}^d \neq \emptyset$.\\
\noindent
(6) For $ - (d+1)/2 < a < d(d+1)/2$, each connected component of
$S_{a,reg}^{\bf C}(b_d)$ intersects with ${\bf T}^d$, and the intersection is a $(d-1)$-dimensional real analytic submanifold of ${\bf T}^d$.
\end{lemma}

Proof. Letting $f_d(z) = 1 + e^{iz_1} + \cdots + e^{iz_d}$, we have the following factorization
\begin{equation}
b_d(z) + \frac{d+1}{2} = \frac{1}{2}\,f_d(z)f_d(-z).
\label{C1S4fdzFactorize}
\end{equation}
Let us compute $S_{a,sng}^{\bf C}(b_d)$. The above equation implies
$$
\nabla_zb_d(z)=0 \Longleftrightarrow e^{iz_j}f_d(-z) = e^{-iz_j}f_d(z), \quad 
1 \leq j \leq d.
$$
{\it Case 1} : $f_d(z) \neq 0$. In this case, $e^{2iz_1} = \cdots = e^{2iz_d}$. Letting this value to be $w^2$, we have $e^{iz_j} = \pm w$. Using $e^{2iz_j}f_d(-z) = f_d(z)$, we then have
$$
w^2\Big(1 + \frac{\ell}{w}\Big) = 1 + \ell w, 
$$
where $\ell$ is an integer satisfying $- d \leq \ell \leq d$. Therefore, $w =\pm 1 =  \pm e^{iz_j}$, hence  $z_j = 0$ or $\pi$. Since $1 + \sum_{j=1}^d e^{iz_j}\neq 0$, we have the  restriction that $\sharp\{j\, ; \, z_j=\pi\} \neq (d+1)/2$ when $d$ is odd.  Therefore, when $d$ is even, 
\begin{equation}
\begin{split}
&\cos z_1 + \cdots + \cos z_d = \ell, \quad \ell = d, d-2, \cdots, -d, \\
& \sin z_1 + \cdots + \sin z_d=0.
\end{split}
\label{C1S4nablazfzequation}
\end{equation}
Taking the square and adding them, we have $\sum_{i<j}\cos(z_i-z_j) = (\ell^2-d)/2$. Hence, 
$$
\sum_{i=1}^d\cos z_i + \sum_{i<j}\cos(z_i-z_j) = \frac{(\ell+1)^2}{2} - \frac{ d+ 1}{2}.
$$
When $d$ is odd,  (\ref{C1S4nablazfzequation}) also holds, however, $\ell \neq d - 2\dfrac{(d+1)}{2} = -1$.

\noindent
{\it Case 2} : $f_d(z)=0$. In this case, in view of (\ref{C1S4fdzFactorize}),
$b_d(z) = -(d+1)/2$.

\medskip

The assertions (1) and (2) now follow from the above observation.

When $a = - (d+1)/2$, by  (\ref{C1S4fdzFactorize}), $ S_{a}^{\bf C}(b_d)$ is a union of two analytic manifolds $\{f_d(z) = 0\}$ and $\{f_d(-z)=0\}$. If they intersect, we can assume without loss of generality that 
$f_{d-1}(-z')\neq 0$, $z' = (z_1,\cdots,z_{d-1})$, at the intersection point. In fact, if $f_{d-1}(-z') = 0$ for all $d-1$ variables $z'$, adding them, we have $\sum_{j=1}^de^{-iz_j} = d/(1-d)$, which is a contradiction.

Then we have $e^{-iz_d} = - f_{d-1}(-z')$, hence $f_{d-1}(z') - 1/f_{d-1}(-z')=0$, which implies $b_{d-1}(z') = - d/2 + \frac{1}{2}f_{d-1}(z')f_{d-1}(-z') = - (d-1)/2$. By (2), $z'$ is on an analytic submanifold of dimension $d-2$ or a discrete set. Therefore $(z',z_d)$ form an analytic submanifold of dimension $d-2$ or a discrete set.
This proves the assertion for $S_{a}^{\bf C}(b_d)$ of (3). The assertion for $S_{a}^{\bf C}(b_d)\cap{\bf T}^d$ is obtained by applying the implicit function theorem for $1 + \sum_{j=1}^d\cos x_j = -1$, $\sum_{j=1}^d\sin x_j=0$. 

By (1), the minimum of $b_d(x)$ on ${\bf T}^d$ is $ - (d+1)/2$. The maximum is easily seen to be $d(d+1)/2$. This proves (4). 

Let us prove (i) of (5). Take $z^{(0)} = (z_1^{(0)},\cdots, z_d^{(0)}) \in \widetilde S_{a,0}(b_d)$. Adding $1 + \sum_{k\neq j}e^{-iz^{(0)}_k}=0$, we obtain $d + (d-1)\sum_{j=1}^de^{-iz^{(0)}_j}=0$, hence
$$
e^{-iz^{(0)}_1} = \cdots = e^{-iz^{(0)}_d} = \frac{1}{1-d}.
$$
This implies $\cos z^{(0)}_j = (d^2 - 2d +2)/(2(1-d))$, $\cos(z^{(0)}_i-z^{(0)}_j) = 1$. Therefore $b_d(z^{(0)}) = 
d/(2(1-d))$. Since the elements of $SV(b_d)$ are half-integers, we have $b_d(z^{(0)}) \not\in SV(b_d)$, which implies $\nabla_z b_d(z^{(0)}) \neq 0$. Near $z^{(0)}$, $S_{a}^{\bf C}(b_d)$ is represented as, say, $z_d = g(z_1,\cdots,z_{d-1})$, where $g$ is analytic. Then one can find $z^{(d)} \in \widetilde S_{a,d}(b_d)$ and a continuous curve in $S_{a,reg}^{\bf C}(b_d)$ with end points $z^{(0)}$ and $z^{(d)}$.

To prove (ii) for the case $j=d$, we let $w = e^{iz_d}$, and rewrite the equation $b_d(z)=a$ as
\begin{equation}
\begin{split}
\big(1 + \sum_{j=1}^{d-1}e^{-iz_j}\big)w^2 + 2(A-a)w + \big(1 + \sum_{j=1}^{d-1}e^{iz_j}\big)=0,
\end{split}
\label{C1S4fd(z)quadratic}
\end{equation}
where $A = \sum_{j=1}^{d-1}\cos z_j + \sum_{1\leq j<k\leq d-1}\cos(z_j-z_k) = b_{d-1}(z_1,\cdots,z_{d-1})$. 
The discriminant $D$ of (\ref{C1S4fd(z)quadratic}) is given by
$$
D/4 = A^2 - 2(a+1)A + a^2 - d.
$$
Then $D=0$ when $A= a+1\pm\sqrt{2a+d+1}$. 

A simple computation yields
\begin{equation}
\begin{split}
& - \frac{d-1}{2} < a + 1 + \sqrt{2a+d+1} \leq  - \frac{d}{2} + 2 \quad {\rm for} \quad 
- \frac{d+1}{2} < a \leq - \frac{d}{2}, \\
& - \frac{d}{2} < a + 1 - \sqrt{2a+d+1} < \frac{d(d-1)}{2} \quad {\rm for} \quad 
- \frac{d}{2} < a <\frac{d(d+1)}{2}.
\end{split}
\nonumber
\end{equation}
By (4), $b_{d-1}(x')\  (x' \in {\bf T}^{d-1})$ varies over the interval $[-d/2,d(d-1)/2]$. Hence if $-(d+1)/2 < a \leq -d/2$, there exists $x' = (x_1,\cdots,x_{d-1}) \in {\bf T}^{d-1}$ such that 
$$
b_{d-1}(x') = a + 1 + \sqrt{2a+d+1},
$$
and if $-d/2 < a < (d+1)d/2$, there exists $x' = (x_1,\cdots,x_{d-1}) \in {\bf T}^{d-1}$ such that 
$$
b_{d-1}(x') = a + 1 -\sqrt{2a+d+1}.
$$

For such $x' = (x_1,\cdots,x_{d-1})$, $1 + e^{-ix_1} +\cdots + e^{-ix_{d-1}} \neq 0$. Otherwise, by (\ref{C1S4fdzFactorize}) with $d$ replaced by $d-1$, we have $b_{d-1}(x') = - d/2$. Therefore, we have either $a+1 + \sqrt{2a+d+1} = - d/2$ (for the case $-(d+1)/2 < a \leq - d/2$), or 
$a + 1 - \sqrt{2a+d+1} = - d/2$ (for the case $-d/2 < a < d(d+1)/2$).
Then we have $a + 1 +d/2 = \pm\sqrt{2a+d+1}$. Hence $a + d/2=0$, which leads to a contradiction.

Then, the equation (\ref{C1S4fd(z)quadratic}) has a double root $w \in {\bf C}$ such that $|w|=1$. Therefore, $w = e^{ix_d}$ for some $x_d \in {\bf T}^1$, hence $x = (x_1,\cdots,x_d) \in 
S_{a}^{\bf C}(b_d)\cap{\bf T}^d$.

Now, take $\zeta = (\zeta_1,\cdots,\zeta_d) \in \widetilde S_{a,d}(b_d)$ so that $1 + e^{-i\zeta_1} + \cdots + e^{-i\zeta_{d-1}} \neq 0$.
 Construct continuous curves $c_j(t)$, $0 \leq t \leq 1$, $j = 1,\cdots,d-1$, such that $c_j(t)\neq 0$, $1 + e^{-ic_1(t)} + \cdots + e^{-ic_{d-1}(t)} \neq 0$ and $c_j(0) = e^{-i\zeta_j}$, $c_j(1)= e^{-ix_j}$. We can then construct a solution $w(t)$ of the equation 
(\ref{C1S4fd(z)quadratic}) with $z_j$ replaced by $c_j(t)$, continuous with respect to $t \in [0,1]$, such that 
$w(0) = e^{i\zeta_d}$, $w(1)=e^{ix_d}$.  Here, we use the fact that (\ref{C1S4fd(z)quadratic})  has a double root for $t=1$. This proves that there is a continuous curve in $\widetilde S_{a,d}(b_d)$ with end points $\zeta$ and $x$, hence the assertion (ii).

We prove (6). Note that each connected component of $S_{a}^{\bf C}(b_d)$ is a union of $\widetilde S_{a,j}(b_d), j\neq 0,$ and possibly a part of $\widetilde S_{a,0}(b_d)$. Take $z^{(0)} \in S_{a}^{\bf C}(b_d)$. If $z^{(0)} \in \widetilde S_{a,j}(b_d)$ for some $j \neq 0$, by (5 - ii), there is a continuous curve $c(t)$ such that $c(0) = z^{(0)}$ and $c(1) \in \widetilde S_{a,j}(b_d)\cap{\bf T}^d$.
 If $z^{(0)} \in \widetilde S_{a,0}(b_d)$, by (5 - i), one can find $\zeta \in \widetilde S_{a,j}(b_d)$ and a continuous curve with end points $z^{(0)}$ and $\zeta$. Then we can apply (5 - ii) again. Here, let us note that since $ a\neq -(d+1)/2$, we can avoid the singular points by (2).
 \qed


\section{Examples of periodic lattices}
We list up examples of periodic graphs, and study the algebraic properties of 
the matrix $H_0(x)$ associated with their Laplacians.
We denote by $\sigma_p(H_0(x))$  the set of all eigenvalues of $H_0(x)$, including the case of $s=1$, and let
\begin{equation}
\sigma(H_0) = {\mathop\cup_{x\in{\bf T}^d}}\sigma_p(H_0(x)),
\label{S2DefineSigma(H0)}
\end{equation}
\begin{equation}
p(x,\lambda) = {\rm det}\left(H_0(x) - \lambda\right),
\label{C1Definep(x,lambda)}
\end{equation}
\begin{equation}
M_{\lambda} = \{x \in {\bf T}^d\, ; \, p(x,\lambda)=0\},
\label{S2DefineMlambda}
\end{equation}
\begin{equation}
M_{\lambda}^{\bf C} = \{z \in {\bf T}^d_{\bf C}\, ; \, p(z,\lambda)=0\}.
\label{S2DefineMlambdaC}
\end{equation}
\begin{equation}
M_{\lambda,reg}^{\bf C} = \{z \in M_{\lambda}^{\bf C}\, ; \, \nabla_z p(z,\lambda)\neq 0\}.
\label{S2DefineMlambdaCreg}
\end{equation}
\begin{equation}
M_{\lambda,sng}^{\bf C} = \{z \in M_{\lambda}^{\bf C}\, ; \, \nabla_z p(z,\lambda)= 0\}.
\label{S2DefineMlambdaCsng}
\end{equation}
\begin{equation}
\widetilde{\mathcal T} = \{\lambda \in \sigma(H_0)\, ; \, M_{\lambda,sng}^{\bf C}\cap {\bf T}^d \neq \emptyset\}.
\label{S2TTilde}
\end{equation}
In this section, we drop the subscript $0$ from $\mathcal V_0$ and $\mathcal L_0$ for the notational convenience.


\subsection{Square lattice}\label{C1SquareLattice}
\begin{figure}[hbtp]
\centering
\includegraphics[width=12cm, bb=0 0 595 392]{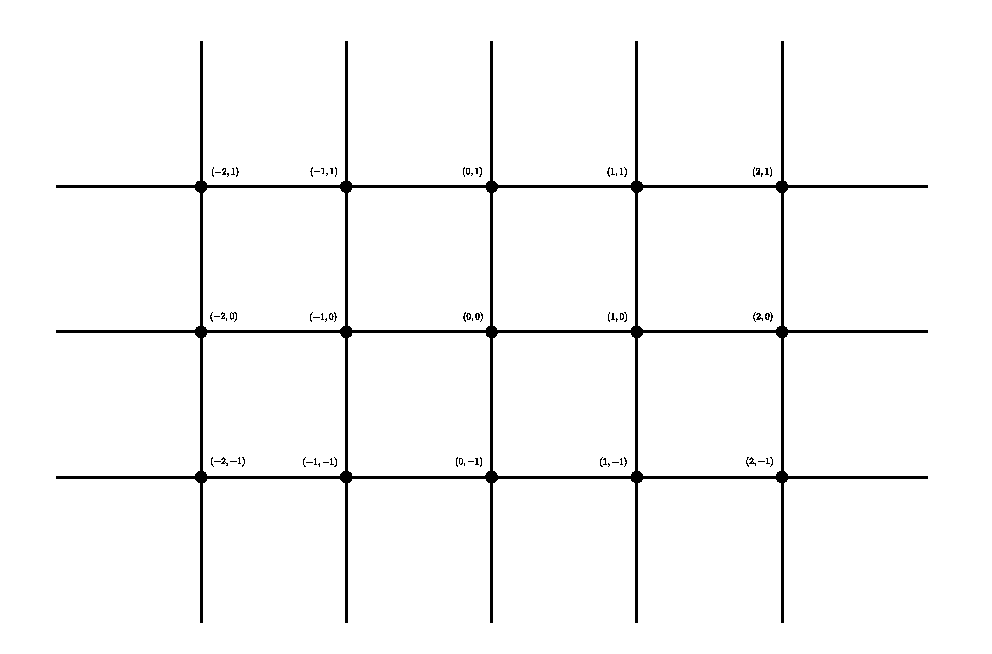}
\caption{Square lattice}
\label{C1SquareLattice}
\end{figure}
Let
$$
\mathcal V = \big\{{\bf v}(n)\, ; \, 
n \in {\bf Z}^d\},\quad
{\bf v}_1 = \big(1,0, \cdots,0\big), \cdots,
{\bf v}_d = \big(0,\cdots, 0,1\big).
$$
$$
\mathcal N_a = \{b \in {\mathcal V}\, ; \, |b-a|=1\} = \{a \pm {\bf v}_1,\cdots, a \pm {\bf v}_d\}, \quad a \in \mathcal V.
$$
Then, the Laplacian is defined by
\begin{equation}
\begin{split}
& \big(\widehat\Delta_{\Gamma} \widehat f\big)(n) =  \frac{1}{2d}
\big(\sum_{i=1}^d \widehat f(n + {\bf v}_i) + \widehat f(n-{\bf v}_i)\big).
\end{split}
\label{C1S3LaplacianSquareLattice}
\end{equation}
Passing to the Fourier series, this Laplacian is transformed into
\begin{equation}
 H_0(x)f(x) = 
-\frac{1}{d}\big(\sum_{i=1}^d\cos x_i\big)f(x).
\label{C1S3FourierLaplacianSquareLattice}
\end{equation}
Therefore, $p(x,\lambda) = -\frac{1}{d}\sum_{j=1}^d\cos x_j - \lambda$, and
Lemma \ref{LemmaSquareLattice} is rewritten as follows.


\begin{lemma}\label{C1S4PropertySquareLattice}
(1) $\ \sigma(H_0) = [-1,1]$. \\
\noindent
(2) $\ \widetilde{\mathcal T}  = \big\{n/d\, ; \, n = - d, - d+2,\cdots, d-2, d\big\}.$ \\
\noindent
(3) For $\lambda \in (-1,1)\setminus\widetilde{\mathcal T}$, $M_{\lambda}$ is a real analytic submanifold of ${\bf T}^d$, and $M_{\lambda}^{\bf C}$ is an analytic submanifold of ${\bf T}^d_{\bf C}$. \\
\noindent
(4) For $-1 < \lambda < 1$, $M_{\lambda,sng}^{\bf C} \subset \left(\pi{\bf Z}\right)^d\cap{\bf T}^d_{\bf C}$. \\
\noindent
(5) For $-1 < \lambda < 1$, each connected component of $M_{\lambda,reg}^{\bf C}$ intersects with ${\bf T}^d$ and the intersection is a $(d-1)$-dimensional real analytic submanifold of ${\bf T}^d$.
\end{lemma}


\subsection{Triangular lattice}
Let
$$
\mathcal V = \big\{{\bf v}(n)\, ; \, 
n \in {\bf Z}^2\}, \quad
{\bf v}_1 = \big(1,0\big), \quad 
{\bf v}_2 = \big(\frac{1}{2},\frac{\sqrt3}{2}\big),
$$
$$
\mathcal N_a = \{b \in {\mathcal V}\, ; \, |b-a|=1\} = \{
a \pm {\bf v_1}, \ a \pm {\bf v}_2, \ a \pm ({\bf v}_1 - {\bf v}_2)\}, \quad a \in \mathcal V.
$$
The Laplacian is defined by
\begin{equation}
\begin{split}
& \big(\widehat\Delta_{\Gamma} \widehat f\big)(n) = \frac{1}{6}
\big(\widehat f(n_1+1,n_2) + \widehat f(n_1-1,n_2) + 
 \widehat f(n_1,n_2+1) \\
& + \widehat f(n_1,n_2-1) + \widehat f(n_1+1,n_2-1)+ \widehat f(n_1-1,n_2+1)\big).
\end{split}
\end{equation}
Passing to the Fourier series, $-\widehat\Delta_{\Gamma}$ is rewritten as
\begin{equation}
H_0(x)f(x) = 
-\frac{1}{3}\big(\cos x_1 + \cos x_2 + \cos(x_1-x_2)\big)f(x).
\end{equation}

\begin{figure}[hbtp]
\centering
\includegraphics[width=12cm, bb=0 0 595 416]{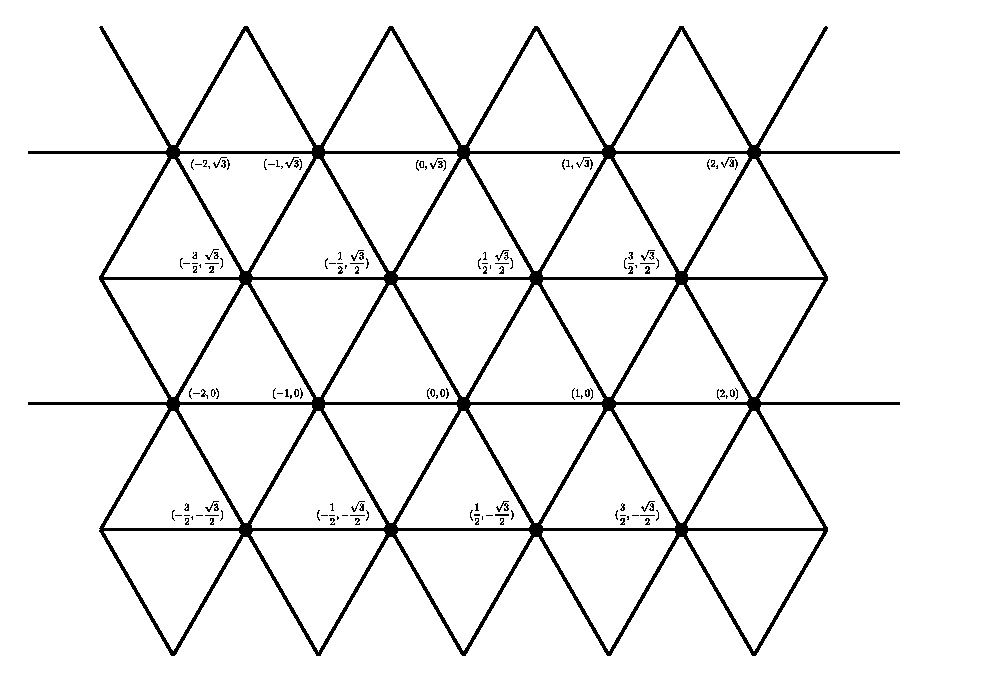}
\caption{Triangular lattice}
\label{C1TrangleLattice}
\end{figure}
Then $p(x,\lambda) =- \frac{1}{3}\big(\cos x_1 + \cos x_2 + \cos(x_1-x_2)\big) - \lambda$, and by Lemma \ref{LemmaDiamond}, we have the following


\begin{lemma}\label{C1S4TriangLatticeProperty}
(1)  $\ \sigma(H_0) = [-1,1/2].$ \\
\noindent
(2)$\ \widetilde{\mathcal T} = \big\{-1, 1/3, 1/2\big\}.$ \\
\noindent
(3) For $\lambda \in (-1,1/2)\setminus\widetilde {\mathcal T}$, $M_{\lambda}$ is a real analytic submanifold of ${\bf T}^2$, and $M_{\lambda}^{\bf C}$ is an analytic submanifold of ${\bf T}^2_{\bf C}$. \\
\noindent
(4) For $-1 < \lambda < 1/2$, $M_{\lambda,sng}^{\bf C} \subset
 \left(\pi{\bf Z}\right)^2\cap{\bf T}^2_{\bf C}$. \\
\noindent
(5) For $-1 < \lambda < 1/2$, each connected components of  $M_{\lambda,reg}^{\bf C}$ intersects with ${\bf T}^2$ and the intersection is a 1-dimensional real analytic submanifold of ${\bf T}^2$.
\end{lemma}


\subsection{Hexagonal lattice}
We put
\begin{equation}
\mathcal L = \big\{{\bf v}(n)\, ; \, n \in {\bf Z^2}\big\}, \quad
 {\bf v}_{1} = \Big( \frac{3}{2},\frac{\sqrt3}{2}\Big), \quad
 {\bf v}_{2} = \Big(0,\sqrt{3}\Big),
\label{C1S4vpm}
\end{equation}
\begin{equation}
p_1 = \big(\frac{1}{2},-\frac{\sqrt3}{2}\big), \quad p_2 = (1,0),
\end{equation}
and define the vertex set $\mathcal V$ by
\begin{equation}
\mathcal V = \mathcal V_1 \cup 
\mathcal V_2, \quad \mathcal V_i = p_i + \mathcal L.
\label{C1S4DefinemathcalV}
\end{equation}
Note that $\mathcal V_1 \cap 
\mathcal V_2 = \emptyset$.
The adjacent points of $a_1 \in \mathcal V_1$ and $a_2 \in \mathcal V_2$ are defined by 
\begin{equation}
\begin{split}
\mathcal N_{a_1}& = \{y \in {\bf R}^2\,;\, |a_1-y|=1\}\cap\mathcal V_2\\
&=
\Big\{
 a_1 + \frac{{\bf v}_1 + {\bf v}_2}{3}, \ a_1 + \frac{{\bf v}_1 - 2{\bf v}_2}{3}, \ a_1 - \frac{2{\bf v}_1 - {\bf v}_2}{3}\Big\},
\end{split}
\end{equation}
\begin{equation}
\begin{split} 
\mathcal N_{a_2} &= \{y \in {\bf R}^2\,;\, |a_2 -y|=1\}\cap \mathcal V_1 \\
& = \Big\{a_2 +\frac{2{\bf v}_1 - {\bf v}_2}{3}, \ a_2 - \frac{{\bf v}_1 -2{\bf v}_2}{3}, \ a_2 - \frac{{\bf v}_1+{\bf v}_2}{3}\Big\}.
\end{split}
\end{equation}

\begin{figure}[hbtp]
\centering
\includegraphics[width=12cm, bb=0 0 595 432]{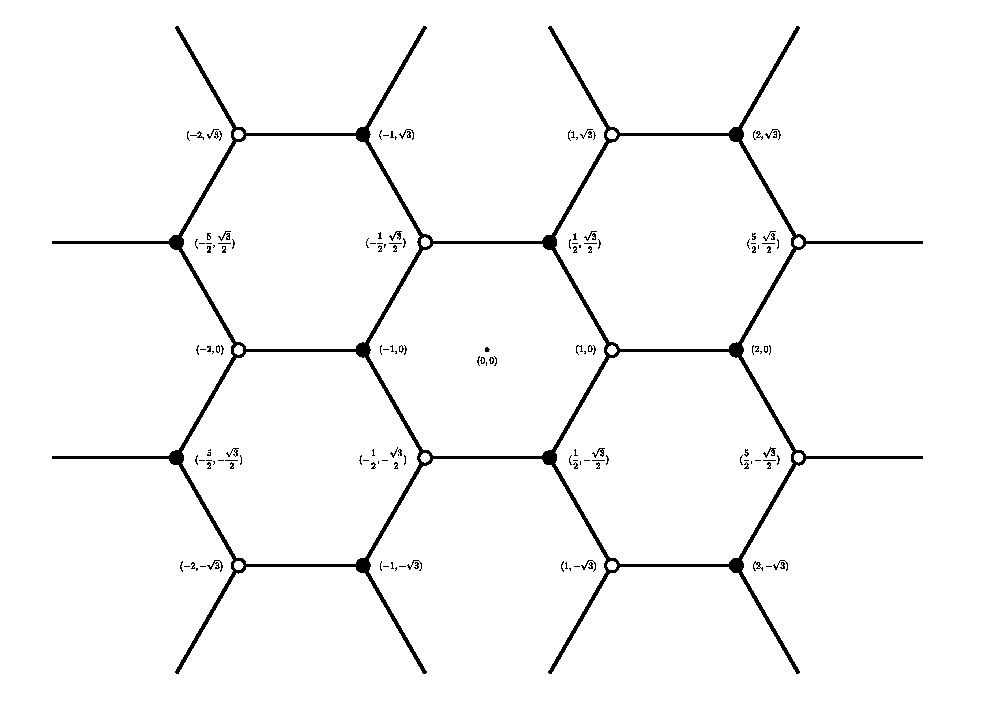}
\caption{Hexagonal lattice}
\label{C1HexagonalLattice}
\end{figure}

For a function $\widehat f(n) = (\widehat f_1(n),\widehat f_2(n))$, the Laplacian is defined by
\begin{equation}
\big(\widehat\Delta_{\Gamma}\widehat f\big)(n) = 
\frac{1}{3}
\left(
\begin{split}
\widehat f_2(n_1,n_2) + \widehat f_2(n_1-1,n_2) + 
\widehat f_2(n_1,n_2-1) \\
\widehat f_1(n_1,n_2) + \widehat f_1(n_1+1,n_2) + 
\widehat f_1(n_1,n_2+1)
\end{split}
\right)
\end{equation}
Passing to the Fourier series, $- \widehat\Delta_{\Gamma}$ is transformed to
\begin{equation}
H_0(x) = 
-\frac{1}{3}\left(
\begin{array}{cc}
0 & 1 + e^{ix_1} + e^{ix_2} \\
1 + e^{-ix_1} + e^{-ix_2}&0
\end{array}
\right).
\end{equation}
A direct computation yields
\begin{equation}
p(x,\lambda) = \det\big(H_0(x)- \lambda\big) = 
\lambda^2 - \frac{\alpha(x)}{9}.
\label{S3pxlambdahexa}
\end{equation}
\begin{equation}
\alpha(x) = 3 + 2\big(\cos x_1 + \cos x_2 + \cos(x_1-x_2)\big).
\label{C1HexagonalLattice_alpha(x)}
\end{equation}
 Lemma \ref{LemmaDiamond} implies the following


\begin{lemma}
(1) $\ \sigma(H_0) = [-1,1].$ \\
\noindent
(2) $\ \widetilde{\mathcal T} = \{-1, - 1/3,0,1/3, 1\}$. \\
\noindent
(3) For $\lambda \in (-1, 1) \setminus \widetilde{\mathcal T}$, $M_{\lambda}$ is a real analytic submanifold of ${\bf T}^2$ and $M_{\lambda}^{\bf C}$ is an analytic submanifold of ${\bf T}^2_{\bf C}$. \\
\noindent
(4) For $-1 < \lambda < 0$ and $0 < \lambda < 1$, $M_{\lambda,sng}^{\bf C}\subset 
\left(\pi{\bf Z}\right)^2\cap{\bf T}^2_{\bf C}$. \\
\noindent
(5) For $-1 < \lambda < 0$ and $0 < \lambda < 1$,  each connected component of $M_{\lambda}^{\bf C}$ intersects with ${\bf T}^2$ and the intersection is a 1-dimensional real analytic submanifold of ${\bf T}^2$ .
\end{lemma}


\subsection{Kagome lattice}\label{C1kagomeLattice}
Let
$$
\mathcal L = \big\{{\bf v}(n)\, ; \, n \in {\bf Z}^2\big\}, \quad
{\bf v}_1 = \big(\frac{1}{2},\frac{\sqrt3}{2}\big), \quad
{\bf v}_2 = \big(-\frac{1}{2},\frac{\sqrt3}{2}\big),
$$
$$
p_1 = (0,0), \quad p_2 = \big(\frac{1}{2},0\big), \quad
p_3 = \big(\frac{1}{4},\frac{\sqrt3}{4}\big).
$$
$$
\mathcal V = \mathcal V_1\cup\mathcal V_2\cup\mathcal V_3, \quad \mathcal V_j = p_j + \mathcal L.
$$
For $a_j \in \mathcal V_j$, the adjacent points are
$$
\mathcal N_{a_j} = \{y \in \mathcal V\, ; \, |a_j - y| = 1/2, y\not\in \mathcal V_j\},
$$
i.e.
\begin{equation}
\begin{split}
\mathcal N_{a_1} & = \Big\{a_1 \pm \frac{{\bf v}_1}{2}, \ 
a_1 \pm \frac{{\bf v}_1-{\bf v}_2}{2}\Big\}, \\
\mathcal N_{a_2} &= \Big\{a_2 \pm \frac{{\bf v}_2}{2}, \ 
a_2 \pm \frac{{\bf v}_1-{\bf v}_2}{2}\Big\}, \\
\mathcal N_{a_3} &= \Big\{a_3 \pm \frac{{\bf v}_1}{2}, \ 
a_3 \pm \frac{{\bf v}_2}{2}\Big\}.
\end{split}
\end{equation}

For a function $\widehat f(n) = (\widehat f_1(n),\widehat f_2(n),\widehat f_3(n))$, the Laplacian is defined by
$$
\big({\widehat\Delta}_{\Gamma}\widehat f\big)(n) = \frac{1}{4}\big(\widehat g_1(n),\widehat g_2(n),\widehat g_(n)\big),
$$
\begin{equation}
\begin{split}
\widehat g_1(n) &= \widehat f_2(n) + \widehat f_2(n_1-1,n_2+1) + 
\widehat f_3(n) + \widehat f_3(n_1-1,n_2),\\
\widehat g_2(n) &= \widehat f_1(n) + \widehat f_1(n_1+1,n_2-1) + 
\widehat f_3(n) + \widehat f_3(n_1,n_2-1),\\
\widehat g_3(n) &= \widehat f_1(n) + \widehat f_1(n_1+1,n_2) + 
\widehat f_2(n) + \widehat f_2(n_1,n_2+1).
\end{split}
\end{equation}
Passing to the Fourier series, $- \widehat\Delta_{\Gamma}$ becomes
\begin{equation}
H_0(x) =- \frac{1}{4} 
\left(
\begin{array}{ccc}
0& 1 + e^{ix_1}e^{-ix_2} &1 + e^{ix_1}\\
1+ e^{-ix_1}e^{ix_2}& 0& 1 + e^{ix_2}\\
1 + e^{-ix_1} &1 + e^{-ix_2}& 0
\end{array}
\right).
\end{equation}
A direct computation gives
\begin{equation}
p(x,\lambda) = \det\big(H_0(x) - \lambda\big) = -\big(\lambda - \frac{1}{2}\big)\big(\lambda^2 + \frac{\lambda}{2} - \frac{\beta(x)}{8}\big),
\label{pxlambdaKagome}
\end{equation}
\begin{equation}
\beta(x) = 1 +  \cos x_1 + \cos x_2 + \cos(x_1-x_2).
\end{equation}
Note that the case $\lambda = 1/2$ is exceptional in that $p(x,1/2) = 0$.
 Lemma \ref{LemmaDiamond} and a direct computation imply the following


\begin{lemma}\label{Kagomeproperty}
(1) $\sigma(H_0) = [-1,1/2]$. \\
\noindent
(2) $\widetilde{\mathcal T} = \{-1, -1/2, -1/4,0,1/2\}$. \\
\noindent
(3) For $(-1,1/2)\setminus\widetilde{\mathcal T}$, $M_{\lambda}$ is a real analytic submanifold of ${\bf T}^2$ and $M_{\lambda}^{\bf C}$ is an analytic submanifold of ${\bf T}^2_{\bf C}$. \\
\noindent
(4) For $-1 < \lambda < -1/4$ and $-1/4 < \lambda < 1/2$, $M_{\lambda,sng}^{\bf C} \subset 
\left(\pi{\bf Z}\right)^2\cap{\bf T}^2_{\bf C}$. \\
\noindent
(5) For $ -1 < \lambda < -1/4$ and $-1/4 < \lambda < 1/2$,  each connected component of $\mathcal M_{\lambda}^{\bf C}$ intersects with ${\bf T}^2$ and the intersection is a 1-dimensional real analytic submanifold of ${\bf T}^2$.\\
\noindent
(6) $H_0(x)$ has an eigenvalue $1/2$ with eigenvector $s(x)v(x)$, where $s(x)$ is an arbitrarily scalar function on ${\bf T}^2$ and
$$
v(x) = \big(- \frac{1}{2}(1 - e^{ix_1})(1 - e^{-ix_2}), 1 - \cos x_1, -\frac{1}{2}(1 - e^{ix_1})(1 + e^{-ix_2})\big).
$$
\end{lemma}

\begin{figure}[hbtp]
\centering
\includegraphics[width=12cm, bb=0 0 595 386]{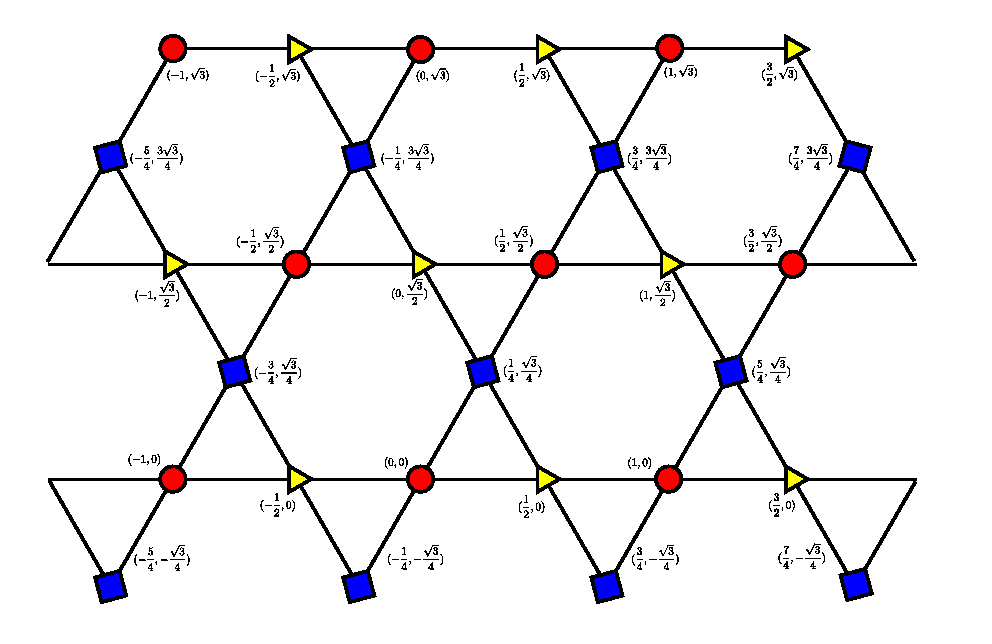}
\caption{Kagome lattice}
\label{C1KagomeLattice}
\end{figure}


\subsection{Diamond lattice}
We put 
\begin{equation}
\mathcal V_1 = \{\ell = (\ell_1,\ell_2,\ell_3) \in {\bf Z}^3\, ; \, \ell_1 + \ell_2 + \ell_3 \in 2{\bf Z}\},
\end{equation}
\begin{equation}
\mathcal V_2 = p + \mathcal V_1, \quad 
p = \big(\frac{1}{2},\frac{1}{2},\frac{1}{2}\big).
\end{equation}
\begin{equation}
\mathcal V = \mathcal V_1\cup{\mathcal V}_2.
\end{equation}
We want to define the adjacent points of $a \in \mathcal V$ as the nearest neighboring points. For this purpose, we prepare the following lemma.


\begin{lemma}\label{Diamondadjacent}
For $\ell \in \mathcal V_1, p + \ell' \in \mathcal V_2$, 
we have $|\ell - \ell'| \geq \sqrt{3}/2$ and 
$|\ell - (\ell' + p)| = \sqrt3/2$ if and only if
\begin{equation}
\ell - \ell' = (0,0,0), \ {\rm or} \ (0,1,1), \ {\rm or} \ (1,0,1), \ {\rm or} \ (1,1,0).
\nonumber
\end{equation}
\end{lemma}
Proof. Let $a = \ell - \ell'$ and $\rho$ be the distance of $\mathcal V_1$ and $\mathcal V_2$. Then we have
\begin{equation}
a_1^2 - a_1 + a_2^2 - a_2 + a_3^2 - a_3 \geq \rho^2 - 3/4.
\label{C1S3DiaLemmma}
\end{equation}
For $m \in {\bf Z}$, $m^2 \geq m$, and $m^2=m$ if and only if $m = 0,1$. Therefore, the left-hand side of (\ref{C1S3DiaLemmma}) is  non-negative, and vanishes for $a_1, a_2, a_3 = 0$ or $1$. The solutions $(1,0,0),(0,1,0),(0,0,1), (1,1,1)$ do not satisfy the condition $a_1 + a_2 + a_3 \in 2{\bf Z}$. On the other hand, $a = (0,0,0)$, $(0,1,1)$, $(1,0,1)$, $(1,1,0)$ meet the condition. This means that $\rho = \sqrt{3}/2$, and the equality of (\ref{C1S3DiaLemmma}) is attained by these  values. \qed

\medskip
By this lemma, the adjacent points of $a \in \mathcal V_1$ are
$$
p + a, \ p + (a_1,a_2-1,a_3-1), \ p + (a_1-1,a_2,a_3-1), \ p + (a_1-1,a_2-1,a_3),
$$
and the adjacent points of $a'+ p \in \mathcal V_2$ are
$$
a', \ (a'_1,a'_2+1,a'_3+1), \  (a'_1+1,a'_2,a'_3+1), \  (a'_1+1,a'_2+1,a_3).
$$


\begin{lemma}\label{C1Diabase} The vectors 
${\bf v}_1  = (0,1,1), \ {\bf v}_2 = (1,0,1), \ {\bf v}_3 = (1,1,0)$ form  a basis of the lattice $\mathcal V_1$.
\end{lemma}
Proof. We put
\begin{equation}
n_1 = \frac{-\ell_1 + \ell_2 + \ell_3}{2}, \quad 
n_2 = \frac{\ell_1 - \ell_2 + \ell_3}{2}, \quad 
n_3 = \frac{\ell_1 + \ell_2 - \ell_3}{2}. 
\end{equation}
Then, one can see that
$$
\ell \in {\bf Z}^3, \ \ell_1 + \ell_2 + \ell_3 \in 2{\bf Z}
\Longleftrightarrow n \in {\bf Z}^3.
$$
We also have
$n_1{\bf v}_1 + n_2{\bf v}_2 + n_3{\bf v}_3 = \ell$.
Therefore $\mathcal V_1 = \{n_1{\bf v}_1 + n_2{\bf v}_2 + n_3{\bf v}_3\, ; \, n \in {\bf Z}^3\}$. \qed

\medskip
In view of Lemma \ref{C1Diabase}, we have, letting
$\mathcal L = \big\{{\bf v}(n)\, ; \, n \in {\bf Z}^3\big\}$,
$$
\mathcal V_1 = \mathcal L, \quad
\mathcal V_2 = p + \mathcal L, \quad
p = \frac{{\bf v}_1 + {\bf v}_2 + {\bf v}_3}{4}.
$$
The edge set is rewritten as follows : 
$\mathcal N$ is the set of points $({\bf v}(n), p + {\bf v}(n')), \ (p + {\bf v}(n'),{\bf v}(n))$
with $n, n'$ satisfying
$$
n - n' = (0, 0, 0), \ (1,0,0), \ (0,1,0), \ (0,0,1).
$$

The Laplacian is defined by
\begin{equation}
\big(\widehat\Delta_{\Gamma}\widehat f)(n) = \frac{1}{4}(\widehat g_1,\widehat g_2),
\end{equation}
\begin{equation}
\begin{split}
\widehat g_1(n) = & \widehat f_2(n)+ \widehat f_2(n_1-1,n_2,n_3) \\
& + 
\widehat f_2(n_1,n_2-1,n_3)+\widehat f_2(n_1,n_2,n_3-1),\\
\widehat g_2(n) =& \widehat f_1(n)
 + \widehat f_1(n_1+1,n_2,n_3) \\
& + 
\widehat f_1(n_1,n_2+1,n_3)+\widehat f_1(n_1,n_2,n_3+1).
\end{split}
\end{equation}

Passing to the Fourier series, $- \widehat\Delta_{\Gamma}$ becomes
\begin{equation}
H_0(x) = -\frac{1}{4}
\left(
\begin{array}{cc}
0 & 1 + e^{ix_1} + e^{ix_2} + e^{ix_3}\\
 1 + e^{-ix_1} + e^{-ix_2} + e^{-ix_3} & 0
\end{array}
\right).
\end{equation}
We then have
\begin{equation}
p(x,\lambda) = \det(H_0(x)  -\lambda) = \lambda^2 - \gamma_3(x) ,
\end{equation}
\begin{equation}
\begin{split}
\gamma_3(x) & = \frac{1}{4} 
 + \frac{1}{8}\big(\cos x_1 + \cos x_2 + \cos x_3 \\
& + 
\cos(x_1-x_2) +\cos(x_2-x_3) + \cos(x_3-x_1)\big).
\end{split}
\end{equation}
 Lemma \ref{LemmaDiamond} implies the following


\begin{lemma}\label{C1S4DiaLattice}
(1) $\sigma(H_0) = [-1,1]$. \\
\noindent
(2) $\widetilde{\mathcal T} = \{-1, - 1/2, 0, 1/2,1\}$. \\
\noindent
(3) For $\lambda \in (-1,1)\setminus\widetilde{\mathcal T}$, $M_{\lambda}$ is a  real analytic submanifold of ${\bf T}^3$ and $M_{\lambda}^{\bf C}$ is an analytic submanifold of ${\bf T}^3_{\bf C}$. \\
\noindent
(4) For $-1 < \lambda < 0$ and $0 < \lambda < 1$, $M_{\lambda,sng}^{\bf C} \subset 
\left(\pi{\bf Z}\right)^3\cap{\bf T}^3_{\bf C}$. \\
\noindent
(5) For $-1 < \lambda < 0$ and $0 < \lambda < 1$,  each connected component of $M_{\lambda}^{\bf C}$ intersects with ${\bf T}^3$ and the intersection is a 2-dimensional analytic submanifold of ${\bf T}^3$.
\end{lemma}


\subsection{Higher-dimensional diamond lattice}\label{C1DiamondLattice}
There is a higher-dimensional analogue of diamond lattice. In fact, the hexagonal lattice and the 3-dimensional diamond lattice are just the cases for $d=2$ and $d=3$ of the lattice $A_d$ defined as follows.
\begin{equation}
A_d = \big\{x = (x_1,\cdots,x_{d+1}) \in {\bf Z}^{d+1}\, ; \, \sum_{i=1}^{d+1}x_i=0\big\}.
\label{C1DEfineAd}
\end{equation}
Let $e_1 = (1,0,\cdots,0)$, $\cdots$, $e_{d+1}=(0,\cdots,0,1)$ be the standard basis of ${\bf R}^{d+1}$, and put
$$
{\bf v}_i = e_{d+1} - e_i, \quad i = 1,\cdots,d.
$$
They satisfy
\begin{equation}
\begin{split}
& |{\bf v}_i|^2 =2, \quad i = 1,\cdots, d, \\
& {\bf v}_i\cdot{\bf v}_j = 1, \quad |{\bf v}_i-{\bf v}_j|^2 = 2, \quad {\rm if} \quad i\neq j.
\end{split}
\end{equation}


\begin{lemma}
Let ${\bf v}(n) = \sum_{i=1}^dn_i{\bf v}_i$, $n \in {\bf Z}^d$. Then
$A_d = \big\{{\bf v}(n)\, ; \, n \in {\bf Z}^d\}$,
and $\big\{{\bf v}_i\big\}_{i=1}^d$ is a basis of $A_d$.
\end{lemma}
Proof. We have an equivalent relation
\begin{equation}
(x_1,\cdots,x_{d+1}) = \sum_{i=1}^dy_i{\bf v}_i \Longleftrightarrow
x_i=-y_i, \ i =1,\cdots,d, \ x_{d+1} = \sum_{i=1}^dy_i.
\label{HidimDiaxandy}
\end{equation}
From this, the lemma follows immediately. \qed

\medskip
We put
\begin{equation}
\mathcal V = A_d\cup\big(p + A_d), \quad 
p = \frac{1}{d+1}\big({\bf v}_1 + \cdots + {\bf v}_d\big).
\label{C1HidimDia}
\end{equation}
This is the vertex set of $d$-dim. diamond lattice.


\begin{lemma}
For ${\bf v}(n), {\bf v}(n') \in A_d$, $|{\bf v}(n) - ({\bf v}(n')+p)| \geq \sqrt{d/(d+1)}$, and the equality occurs if and only if
$$
n - n' = (0,\cdots,0), \ (1,0,\cdots,0),\ \cdots, \ (0,\cdots,0,1).
$$
\end{lemma}
Proof. We put $a = {\bf v}(n)-{\bf v}(n')$. Then
$$
a-p = \Big(a_1+\frac{1}{d+1},\cdots,a_d+\frac{1}{d+1},
-\sum_{i=1}^da_i - \frac{d}{d+1}\Big).
$$
Then we have
\begin{equation}
|a-p|^2 = {\mathop\sum_{i=1}^d}a_i^2 -1 
+ \big({\mathop\sum_{i=1}^d}a_i + 1\big)^2 + 
\frac{d}{d+1}.
\nonumber
\end{equation}
This is always greater than or equal to $d/(d+1)$, and the equality occurs if and only if all $a_i=0$, or one of $a_i=-1$ and the others vanish. Taking into account of (\ref{HidimDiaxandy}), we obtain the lemma.  \qed

\medskip
Therefore, for ${\bf v}(n) \in A_d$, we define its adjacent points by
\begin{equation}
\begin{split}
\mathcal N_{{\bf v}(n)} = \big\{p + {\bf v}(n')\, ;& \, n-n' =  (0,\cdots,0),  \\
&(1,0,\cdots,0),
\cdots,  (0,\cdots,0,1)\big\},
\end{split}
\nonumber
\end{equation}
and for $p + {\bf v}(n') \in p + A_d$,
\begin{equation}
\begin{split}
\mathcal N_{p + {\bf v}(n')} = \big\{{\bf v}(n)\, ;& \, n-n' =  (0,\cdots,0),  \\
&(1,0,\cdots,0),
\cdots,  (0,\cdots,0,1)\big\}.
\end{split}
\nonumber
\end{equation}

The Laplacian is defined by
\begin{equation}
\big(\widehat{\Delta}_{\Gamma} \widehat{f}\big)(n) = 
\frac{1}{d+1}\big(\widehat g_1(n),\widehat g_2(n)\big),
\end{equation}
\begin{equation}
\widehat g_1(n) = \widehat f_2(n) + \widehat f_2(n-{\bf e}_1) + \cdots +\widehat f_2(n-{\bf e}_d),
\end{equation}
\begin{equation}
\widehat g_2(n) = \widehat f_1(n) + \widehat f_1(n+{\bf e}_1) + \cdots +\widehat f_1(n+{\bf e}_d),
\end{equation}
where $\{{\bf e}_i\}_{i=1}^d$ is the standard basis of ${\bf R}^d$. Passing to the Fourier series, $- \widehat\Delta_{\Gamma}$ is transformed to
\begin{equation}
H_0(x) =- \frac{1}{d+1}\left(
\begin{array}{cc}
0 & 1 + e^{ix_1} + \cdots + e^{ix_d} \\
1 + e^{-ix_1} + \cdots + e^{-ix_d} & 0
\end{array}
\right).
\end{equation}
We have
\begin{equation}
p(x,\lambda) = \det(H_0(x)  -  \lambda\big) = 
\lambda^2 -\gamma_d(x),
\label{S3ddimdiapxlambda}
\end{equation}
\begin{equation}
\gamma_d(x) = \frac{1}{d+1} + \frac{2}{(d+1)^2}\Big(\sum_{i=1}^d\cos x_i + \sum_{i<j}\cos(x_i-x_j)\Big).
\end{equation}
 Lemma \ref{LemmaDiamond} implies the following


\begin{lemma}
(1)  $\sigma(H_0) = [-1,1]$. \\
\noindent
(2) 
\begin{equation}
\widetilde{\mathcal T} = \left\{
\begin{split}
&\Big\{\pm \dfrac{\ell + 1}{d+1}\, ; \, \ell = d, d-2, \cdots, -d\Big\}\cup\big\{0\big\}, \quad {\rm if} \quad d = {\rm even}, \\
& \Big\{\pm \dfrac{\ell +1}{d+1}\, ; \, \ell = d, d-2,\cdots, -d\Big\}, \quad {\rm if} \quad d = {\rm odd}.
\end{split}
\right.
\nonumber
\end{equation}
(3) For $(-1,1)\setminus\widetilde{\mathcal T}$, $M_{\lambda}$ is a real analytic submanifold of ${\bf T}^d$ and $M_{\lambda}^{\bf C}$ is an analytic submanifold of ${\bf T}^d_{\bf C}$. \\
\noindent
(4) For $-1 < \lambda < 0$ and $0 < \lambda < 1$, $M_{\lambda,sng}^{\bf C} \subset \left(\pi{\bf Z}\right)^d\cap{\bf T}^d_{\bf C}$.\\
\noindent
(5) For $-1 < \lambda < 0$ and $0 < \lambda < 1$,  each connected component of $M_{\lambda}^{\bf C}$ intersects with ${\bf T}^d$ and the intersection is a $(d-1)$-dimensional real analytic submanifold of ${\bf T}^d$.
\end{lemma}


\subsection{Subdivision of $d$-dimensional Square Lattice ${\bf Z}^d$, $d\ge2$}\label{SubsectionSubdivision}

\begin{figure}[hbtp]
\centering
\includegraphics[width=10cm, bb=0 0 595 547]{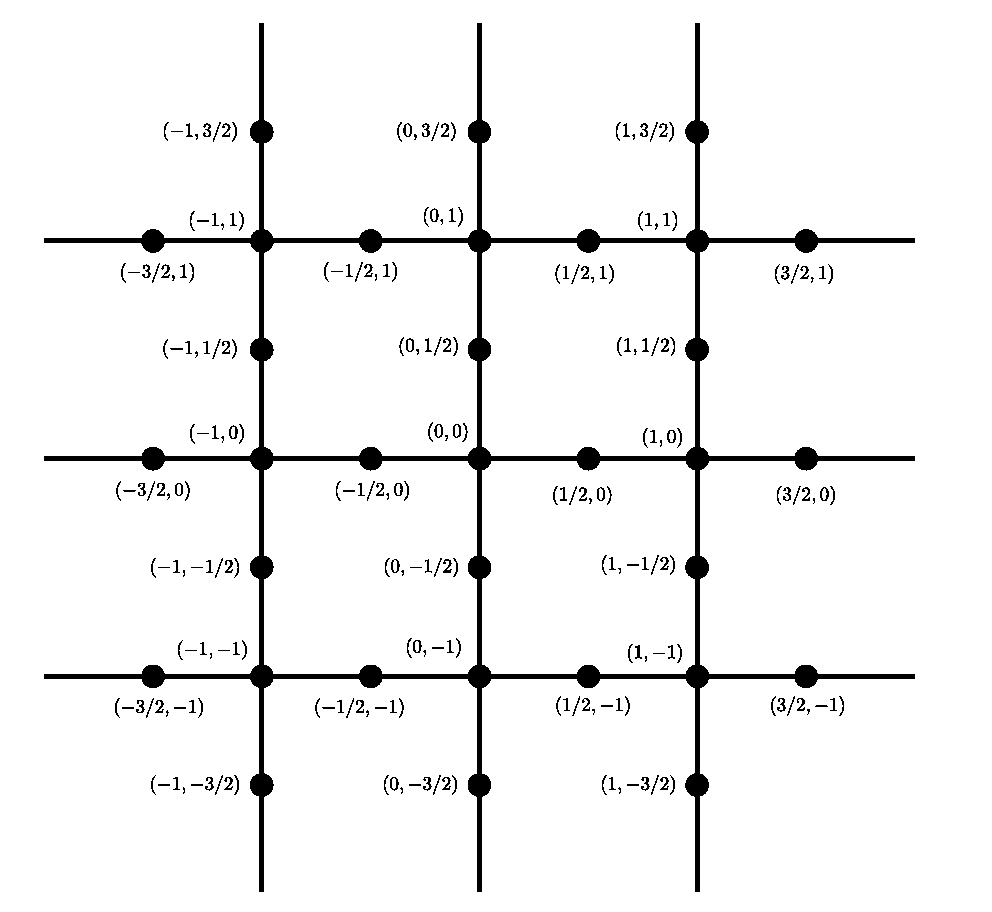}
\caption{Subdivision of $2$-dimensional square lattice}
\label{C1Subdivision}
\end{figure}

The next example is obtained by adding new vertices in the middle points of the edges of the square lattice ${\bf Z}^d$. Let us put
\begin{equation}
{\bf v}_1 = \big(1, 0, \cdots, 0\big),\ \cdots,\ {\bf v}_d = \big(0, \cdots, 0, 1\big),
\end{equation}
\begin{equation}
p_1 = \big(0, \cdots, 0\big),\ p_2 = \big(1/2, 0, \cdots ,0\big),\ \cdots,\ p_{d+1} = \big(0, \cdots, 0, 1/2\big),
\end{equation}
\begin{equation}
\mathcal L = \{{\bf v}(n)\, ; \, n \in {\bf Z}^d\},
\end{equation}
\begin{equation}
\mathcal V = {\mathop \cup_{j=1}^{d+1}} V_j,\ \mathcal V_j = p_j + \mathcal L.
\end{equation}
The edge relations are defined by
\begin{equation}
\mathcal N_{a_j} = \{y \in \mathcal V\, ; \, |y - a_j| = 1/2\}, \quad a_j \in \mathcal V_j.
\end{equation}
Then $a_1  \in \mathcal V_1$ has $2d$ adjacent points, while $a_j \in \mathcal V_j$, $j=2, \cdots, d+1$, have $2$ adjacent points. The Laplacian is then defined by
\begin{equation}
\big(\widehat\Delta_{\Gamma}\widehat f\big)(n) = \frac{1}{2}
\left(
\begin{array}{c}
\displaystyle\frac{1}{d}\sum_{j=1}^{d}\left(\widehat f_{j+1}(n) + \widehat f_{j+1}(n-{\bf e}_j)\right) \\
\widehat f_1(n) + \widehat f_1(n + {\bf e}_1) \\
\vdots \\
\widehat f_1(n) + \widehat f_1(n + {\bf e}_d)
\end{array}
\right).
\end{equation}
Passing to the Fourier series, $- \widehat\Delta_{\Gamma}$ becomes the following matrix
\begin{equation}
H_0(x) = -\frac{1}{2\sqrt{d}}
\left(
\begin{array}{cccc}
0 & 1 + e^{ix_1} & \cdots & 1+e^{ix_d}\\
1 + e^{-ix_1} & 0 & \cdots & 0 \\
\vdots & \vdots & \ddots & \vdots \\
1 + e^{-ix_d} & 0 & \cdots & 0
\end{array}
\right),
\end{equation}
whose determinant is computed as
\begin{equation}
p(x,\lambda) = \det(H_0(x) - \lambda) =  (-\lambda)^{d-1} \big(\lambda^2 - \frac{1}{2d}
(d + \sum_{j=1}^d\cos x_j)\big).
\label{Subdivisionpxlambda}
\end{equation}

Similarly to the case of Kagome lattice, the case $\lambda = 0$ is exceptional since
 $p(x,0) = 0$.  Lemma \ref{LemmaSquareLattice}  and a direct computation imply the following


\begin{lemma}
(1) $\sigma (H_0) = [-1,1]$.  \\
\noindent
(2) $\widetilde{\mathcal T} = \{0, \pm \sqrt{n/2d};\,n=1, 2, \cdots, 2d\}.$ \\
\noindent
(3) For $\lambda \in (-1,1)\setminus\widetilde{\mathcal T}$, $M_{\lambda}$ is a  real analytic submanifold of ${\bf T}^d$ and $M_{\lambda}^{\bf C}$ is an analytic submanifold of ${\bf T}^d_{\bf C}$. \\
\noindent
(4) For $-1 < \lambda <0$, $0 < \lambda < 1$, $M_{\lambda,sng}^{\bf C} \subset 
\left(\pi{\bf Z}\right)^d\cap{\bf T}^d_{\bf C}$.  \\
\noindent
(5) For $-1 < \lambda <0$, $0 < \lambda < 1$, $M_{\lambda,reg}^{\bf C}$ intersects with ${\bf T}^d$ and the intersection is a $d-1$-dimensional real analytic submanifold of ${\bf T}^d$. \\
\noindent
(5) $H_0(x)$ has an eigenvalue $0$, whose eigenvector is written as
$$
\sum_{j=1}^{d-1}s_j(x)v_j(x),
$$
where $s_1(x)$, $\cdots$, $s_{d-1}(x)$ are arbitrary scalar functions and
\begin{equation}
  \begin{split}
    v_1(x)&=\big(0, -(1+e^{ix_2}), 1+e^{ix_1}, 0, \cdots, 0\big), \\
    v_2(x)&=\big(0, -(1+e^{ix_3}), 0, 1+e^{ix_1}, \cdots, 0\big), \\
    &\cdots \\
    v_{d-1}(x)&=\big(0, -(1+e^{ix_d}), 0, 0, \cdots, 1+e^{ix_1}\big).
  \end{split}
  \nonumber
\end{equation}
\end{lemma}


\subsection{Ladder of $d$-dimensional Square Lattice in ${\bf R}^{d+1}$}\label{SubsectionLadder} 
The term ``ladder'' is named after the shape of the following graph ($d=1$):

\begin{figure}[hbtp]
\centering
\includegraphics[width=12cm, bb=0 0 595 183]{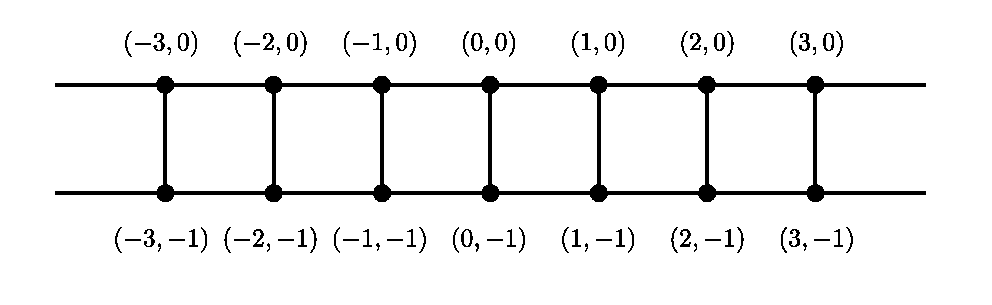}
\caption{2-dim. ladder}
\label{C2-dim.ladder}
\end{figure}

\begin{figure}[hbtp]
\centering
\includegraphics[width=12cm, bb=0 0 595 476]{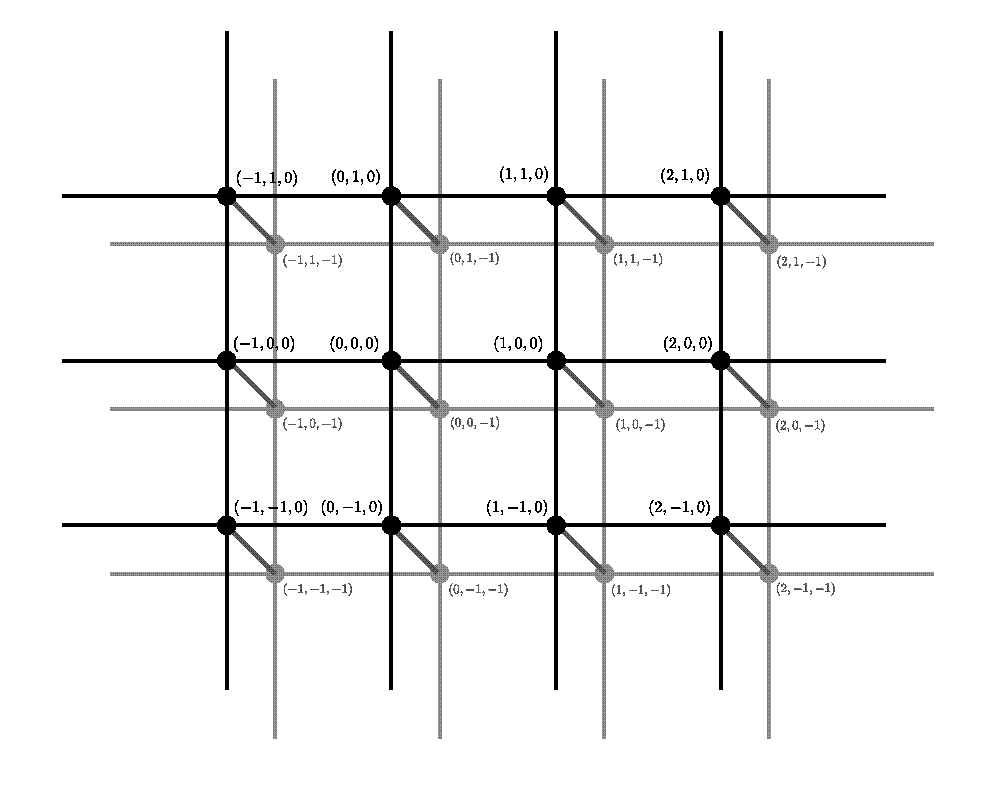}
\caption{3-dim. ladder}
\label{C3-dim.ladder}
\end{figure}

The $d$-dimensional ladder is defined as follows. 
Let $\mathcal L_d$ be the standard $d$-dim. square lattice realized in ${\bf R}^{d+1}$, i.e. $\mathcal L_d = \{(n_1,\cdots,n_d,0)\, ; \, n_i \in {\bf Z}\}$, and put
\begin{equation}
\mathcal V_1 = \mathcal L_d, \quad \mathcal V_2 = (0,\cdots,0,-1) + \mathcal L_d,
\end{equation}
$$
p_0 = (0,\cdots,0), \quad p_1 = (0,\cdots,0,-1),
$$
\begin{equation}
\mathcal V = \mathcal V_1 \cup \mathcal V_2, \quad \mathcal V_i = p_i + \mathcal L.
\end{equation}
The adjacent relation is defined by
\begin{equation}
\mathcal N_a = \{y \in \mathcal V\, ; \, |y-a| = 1\}, \quad a \in \mathcal V.
\end{equation}

Then the Laplacian is
\begin{equation}
(\widehat\Delta_{\Gamma}\widehat f\big)(n) = \frac{1}{2d+1}
\left(
\begin{array}{c}
\widehat f_2(n) + \sum_{j=1}^d\big(\widehat f_1(n + {\bf e}_j) + \widehat f_1(n - {\bf e}_j)\\
\widehat f_1(n) + \sum_{j=1}^d\big(\widehat f_2(n + {\bf e}_j) + \widehat f_2(n - {\bf e}_j)
\end{array}
\right).
\end{equation}
Passing to the Fourier series, $- \widehat\Delta_{\Gamma}$ is transformed to
\begin{equation}
H_0(x) = -\frac{1}{2d+1}
\left(
\begin{array}{cc}
2\sum_{j=1}^d \cos x_j & 1  \\
1 & 2\sum_{j=1}^d\cos x_j
\end{array}
\right).
\end{equation}
Then we have
$$
p(x,\lambda) = \det(H_0(x)-\lambda) = p_+(x,\lambda)p_-(x,\lambda),
$$
$$
p_{\pm}(x,\lambda) = \lambda +  \frac{1}{2d+1}\Big(2\sum_{j=1}^d\cos x_j \pm 1\Big).
$$
Then, $H_0(x)$ has two distinct eigenvalues  $\lambda_{\pm}(x) = 
(- 2\sum_{j=1}^d\cos x_j \pm 1)/(2d+1)$ with values 
$$
- 1 \leq \lambda_-(x) \leq \frac{2d-1}{2d+1}, \quad
\frac{-2d+1}{2d+1} \leq \lambda_+(x) \leq 1.
$$
Accordingly, $M_{\lambda}^{\bf C}$ is split into 2 parts : 
$$
M_{\lambda}^{\bf C} = M_{\lambda,+}^{\bf C}\cup M_{\lambda,-}^{\bf C}, \quad 
M_{\lambda,\pm}^{\bf C} = \{z \in {\bf T}_{\bf C}^{d}\, ; \, p_{\pm}(z,\lambda)=0\}.
$$
We put
\begin{equation}
\mathcal T_{\pm} = \{\lambda \, ; \, p_{\pm}(x,\lambda)=0,\ \nabla p_{\pm}(x,\lambda)=0 \ {\rm for} \ {\rm some} \ x \in {\bf T}^d\},
\end{equation}
which is equal to
$$
\widetilde{\mathcal T}_+ = \Big\{\frac{-2d+1}{2d+1}, \frac{-2d+5}{2d+1},\cdots,1\Big\},
$$
$$
\widetilde{\mathcal T}_- = \Big\{-1, \frac{-2d+3}{2d+1}, \cdots,\frac{2d-1}{2d+1}\Big\}.
$$

In view of Lemma \ref{LemmaSquareLattice}, we have the following


\begin{lemma}
(1) $\sigma(H_0) = [-1,1].$ \\
\noindent
(2) For $\lambda \in \big(-1,\frac{2d-1}{2d+1}\big)\setminus\widetilde{\mathcal T}_-$, $M_{\lambda,-}$ is a real analytic submanifold of ${\bf T}^d$ and $M_{\lambda,-}^{\bf C}$ is an analytic submanifold of ${\bf T}^d_{\bf C}$. \\
\noindent
(3)  For $\lambda \in \big(\frac{-2d+1}{2d+1},1\big)\setminus\widetilde{\mathcal T}_+$, $M_{\lambda,+}$ is a real analytic submanifold of ${\bf T}^d$ and $M_{\lambda,+}^{\bf C}$ is an analytic submanifold of ${\bf T}^d_{\bf C}$. \\
\noindent
(4) For $-1 < \lambda < \frac{2d-1}{2d+1}$, each connected component of $M_{\lambda,-}^{\bf C}$ intersects with ${\bf T}^d$ and the intersection is a $(d-1)$-dimensional real analytic submanifold of ${\bf T}^d$. \\
\noindent
(5) For $\frac{-2d+1}{2d+1} < \lambda < 1$,  each connected component of $M_{\lambda,+}^{\bf C}$ intersects with ${\bf T}^d$ and the intersection is a $(d-1)$-dimensional real analytic submanifold of ${\bf T}^d$. \\
\noindent
(6) For $-1< \lambda < \frac{-2d+1}{2d+1}$, $M_{\lambda,+}^{\bf C}\cap{\bf T}^d = \emptyset$. \\
\noindent
(7) For $ \frac{2d-1}{2d+1} < \lambda < 1$, $M_{\lambda,-}^{\bf C}\cap{\bf T}^d = \emptyset$.
\end{lemma}

\begin{figure}[hbtp]
\centering
\includegraphics[width=10cm, bb=0 0 483 179]{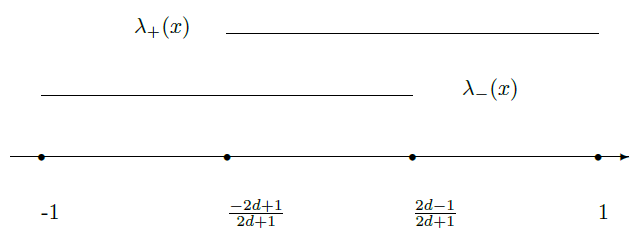}
  \caption{Eigenvalues for the Ladder}
\label{fig:eigenvalueLadder}
\end{figure}


\subsection{Graphite in ${\bf R}^3$}

\begin{figure}
\centering
\includegraphics[width=12cm, bb=0 0 595 439]{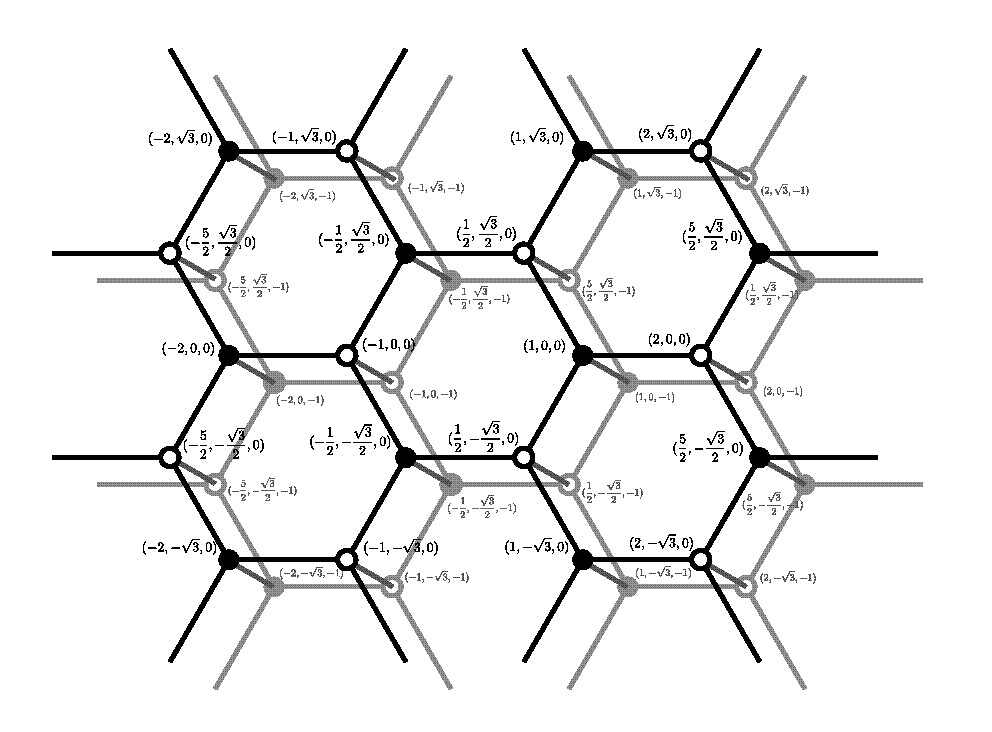}
\caption{Graphite}
\label{C1Graphite}
\end{figure}

The graphite has the same structure as above with the square lattice replaced by the hexagonal lattice. We put
\begin{equation}
\mathcal L_2 = \big\{{\bf v}(n) = n_1{\bf v}_1 + n_2{\bf v}_2\, ; \, n \in {\bf Z}^2\big\},
\end{equation}
\begin{equation}
{\bf v}_1 = \Big(\frac{3}{2}, \frac{\sqrt3}{2},0\Big), \quad {\bf v}_2 = \big(0,\sqrt 3,0\big),
\end{equation}
\begin{equation}
p_1 = \Big(\frac{1}{2}, - \frac{\sqrt 3}{2},0\Big), \quad p_2 = (1,0,0), 
\end{equation}
\begin{equation}
p_3 = p_1 + (0,0,-1), \quad p_4 = p_2 + (0,0,-1),
\end{equation}
and define the vertex set $\mathcal V$ by
\begin{equation}
\mathcal V = {\mathop\cup_{i=1}^4}\mathcal V_i, \quad 
\mathcal V_i = p_i + \mathcal L_2.
\end{equation}
The  adjacent relation is defined by
\begin{equation}
\mathcal N_a = \{y \in \mathcal V\, ; \, |y - a| = 1\}, \quad a \in \mathcal V.
\end{equation}
For a function $\widehat f(n) = (\widehat f_1(n),\widehat f_2(n),\widehat f_3(n),\widehat f_4(n))$, the Laplacian is defined by
\begin{equation}
 \big(\widehat \Delta_{\Gamma}\widehat f\big)(n) =
\frac{1}{4}\left(
\begin{split}
&\widehat f_3(n) +\widehat f_2(n) + \widehat f_2(n_1-1,n_2) + \widehat f_2(n_1,n_2-1) \\
&\widehat f_4(n) + \widehat f_1(n) + \widehat f_1(n_1+1,n_2) + \widehat f_1(n_1,n_2+1) \\
&\widehat f_1(n) + \widehat f_4(n) + \widehat f_4(n_1-1,n_2) + \widehat f_4(n_1,n_2-1) \\
&\widehat f_2(n) + \widehat f_3(n) + \widehat f_3(n_1+1,n_2) + \widehat f_3(n_1,n_2+1) 
\end{split}
\right).
\nonumber
\end{equation}
Passing to the Fourier series, $-\widehat\Delta_{\Gamma}$ is written as $H_0(x)$, where 
\begin{equation}
H_0(x) = 
-\frac{1}{4}
\left(
\begin{array}{cccc}
0 & \overline{c(x)}& 1&0 \\
c(x) & 0 & 0 & 1 \\
1 & 0 & 0 & \overline{c(x)} \\
0 & 1 & c(x) & 0
\end{array}
\right), \
\end{equation}
\begin{equation}
c(x) = 1 + e^{-ix_1} + e^{-ix_2}.
\label{C1S4Graphite_c(x)}
\end{equation}
Then we have, letting $|c|^2 = \alpha(x) = 3 + 2\big(\cos x_1 + \cos x_2 + \cos(x_1-x_2)\big)$,
\begin{equation}
p(x,\lambda) = \det(H_0(x)-\lambda) = 
\lambda^4 - \frac{\alpha+1}{8}\lambda^2 + \frac{(\alpha-1)^2}{4^4}.
\end{equation}
Therefore, $H_0(x)$ has 4 eigenvalues $\pm (\sqrt{\alpha(x)}\pm 1)/4$. We label them as
\begin{equation}
\begin{array}{l}
\lambda_1(x) = - \dfrac{1}{4} - \dfrac{\sqrt{\alpha(x)}}{4},\\
\lambda_2(x) =\left\{
\begin{split}
&- \frac{1}{4} + \frac{\sqrt{\alpha(x)}}{4}, \quad {\rm if} \quad {\alpha(x) \leq 1}\\
& \frac{1}{4} - \frac{\sqrt{\alpha(x)}}{4}, \quad {\rm if} \quad {\alpha(x) \geq 1},
\end{split}
\right.
 \\
\lambda_3(x) =\left\{
\begin{split}
&\frac{1}{4} - \frac{\sqrt{\alpha(x)}}{4}, \quad {\rm if} \quad {\alpha(x) \leq 1}\\
&- \frac{1}{4} + \frac{\sqrt{\alpha(x)}}{4}, \quad {\rm if} \quad {\alpha(x) \geq 1},
\end{split}
\right.\\
 \lambda_4(x) =  \dfrac{1}{4} + \dfrac{\sqrt{\alpha(x)}}{4}.
\end{array}
\label{C1S4GraphiteEigenvalues}
\end{equation}
Then we have
$$
\lambda_1(x) \leq \lambda_2(x) \leq \lambda_3(x) \leq \lambda_4(x),
$$
and they are distinct if $\alpha(x) \neq 0, 1$. 
Furthermore, $\nabla \lambda_j(x)\neq 0$ if 
$0 < \alpha(x) < 1$, $1 < \alpha(x) < 9$. Note that
\begin{equation}
\alpha = 0, 1 \quad 
{\rm on} \quad M_{\lambda} \Longleftrightarrow \lambda = 0, \pm\frac{1}{4}, \pm \frac{1}{2} .
\label{C1S4Graphitegradlambda(j)}
\end{equation}

\begin{figure}[hbtp]
\centering
\includegraphics[width=10cm, bb=0 0 483 274]{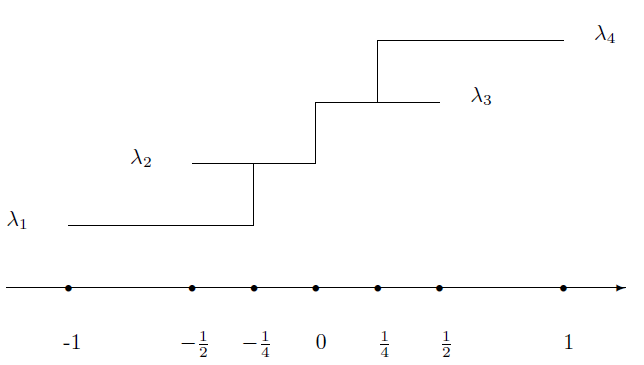}
  \caption{Eigenvalues for the Graphite}
\label{fig:eigenvalueGraphite}
\end{figure}

Letting
\begin{equation}
M_{\lambda}^{(j)} = \{x \in {\bf T}^2\, ; \, \lambda_j(x)=\lambda\},
\end{equation}
we have
$$
M_{\lambda} = {\mathop\cup_{j=1}^4}M_{\lambda}^{(j)}.
$$

We need another splitting of $M_{\lambda}$. Let
\begin{equation}
M_{\lambda,\pm}^{\bf C} = \{ z \in {\bf T}^2_{\bf C}\, ; \, 
\cos z_1 + \cos z_2 + \cos(z_1 - z_2) = 8\lambda^2 \pm 4\lambda - 1\}, 
\end{equation}
\begin{equation}
M_{\lambda,\pm} = M_{\lambda,\pm}^{\bf C}\cap{\bf T}^2.
\end{equation}
Then  $p=0 \Longleftrightarrow 16\lambda^2 = (\sqrt{\alpha} \pm 1)^2 
\Longleftrightarrow \alpha = 16\lambda^2 \pm 8\lambda +1$. This yields
\begin{equation}
M_{\lambda} = M_{\lambda,+}\cup M_{\lambda,-}, \quad 
M_{\lambda}^{\bf C} = M_{\lambda,+}^{\bf C}\cup M_{\lambda,-}^{\bf C}.
\end{equation}


\begin{lemma}\label{C1S4GraphiteMlambda}
(1) $\sigma (H_0) = [-1,1]$. \\
\noindent
(2) $\widetilde{\mathcal T} =\{0, \pm 1/4, \pm 1/2, \pm 1\}$. \\
\noindent
(3) For $\lambda \in (-1,1)\setminus\widetilde{\mathcal T}$ and $ 1 \leq j \leq 4$, $M_{\lambda}^{(j)}$ is a real analytic submanifold of ${\bf T}^2$. \\
\noindent
(4) For $-1 < \lambda < -1/4$, $-1/4 < \lambda < 1/4$ and $1/4 < \lambda < 1$, $M_{\lambda,sng}^{\bf C} \subset \left(\pi{\bf Z}\right)^2\cap{\bf T}^2_{\bf C}$. \\
\noindent
(5) For $-1 < \lambda < -1/4$ and $-1/4 < \lambda < 1/2$, each connected component of $M_{\lambda,+}^{\bf C}$ intersects with ${\bf T}^2$ and the intersection is a 1-dimensional real analytic submanifold of ${\bf T}^2$. \\
\noindent
(6) For $-1/2 < \lambda < 1/4$ and $1/4 < \lambda < 1$, each connected component of $M_{\lambda,-}^{\bf C}$ intersects with ${\bf T}^2$ and the intersection is a 1-dimensional real analytic submanifold of ${\bf T}^2$. \\
(7) For $-1 < \lambda < -1/2$, $M^{{\bf C}}_{\lambda, -} \cap {\bf T}^2 = \emptyset$. \\
\noindent
(8) For $\frac{1}{2} < \lambda < 1$, $M^{{\bf C}}_{\lambda, +} \cap {\bf T}^2 = \emptyset$. \\
\end{lemma}

Proof. The assertion (1) follows from Lemma \ref{LemmaDiamond} (3).  To prove (2), let
$$
a_i = \frac{\partial}{\partial x_i}\alpha = - 2\sin x_i - 2\sin(x_i-x_j), \quad i\neq j.
$$
Then we have
$$
\frac{\partial }{\partial x_i}p = a_i\Big(- \frac{\lambda^2}{8} + \frac{2}{4^2}(\alpha -1)\Big).
$$
A simple computation shows that
$$
- \frac{\lambda^2}{8} + \frac{2}{4^4}(\alpha - 1)^2 = 0, \quad 
\lambda^4 - \frac{\alpha + 1}{8}\lambda^2 + \frac{(\alpha - 1)^2}{4^4} = 0
$$
if and only if
$$
\alpha = 1, \quad \lambda = 0.
$$
Therefore, if $p(x,\lambda)= 0$, $\nabla_xp(x,\lambda)=0$, $\lambda \neq 0$, we have $a_1 = a_2=0$, which implies $(x_1,x_2) = (0,0), (0,\pi), (\pi,0), (\pi,\pi), (2\pi/3,4\pi/3), (4\pi/3,2\pi/3)$. For these values, $\alpha = 0, 1,9$ and $\lambda^2 = 0, 1/16, 1/4, 1$. This proves (2).

The assertion (3) follows from (\ref{C1S4Graphitegradlambda(j)}).

The assertion (4) follows from Lemma~\ref{LemmaDiamond} (2).

In view of Lemma~\ref{LemmaDiamond} (4) and (5), we have (5), (6), (7), and (8). \qed


\subsection{Graph-operation and characteristic polynomials}
Let us observe the above examples from a view point of graph-operation, which is a method of creating  new graphs from the given one.  It is worthwhile to note the general relations   (\ref{charpolynomlinegraph}), (\ref{Rrmark3.13pxlambda}), (\ref{char_pol_ladder}) between the characteristic polynomial of the Laplacian for the resulting graph and that of the original graph. We omit the proof of these formulas, although they follow from straightforward computation,  since we do not use them in this paper. The arguments in the following sub-subsections \ref{subsublinegraph},  \ref{subsubSubdivision} are based on  \cite{S99} and \cite{HiShi04}. 

We recall some notions in the graph theory. 

A graph $\Gamma$ is said to be {\it $k$-regular} if ${\rm deg}\,(v) = k$ for any $v \in \mathcal V(\Gamma)$. 

A {\it $(k_1, k_2)$-semiregular} graph $\Gamma$ is,  by definition, a {\it bipartite} graph with two partite sets $\mathcal V_1$ and $\mathcal V_2$, i.e. $\mathcal{V}(\Gamma)=\mathcal V_1 \cup \mathcal V_2$, $\mathcal V_1 \cap \mathcal V_2 = \emptyset$, and $\mathcal{N}_{v_1} \subset \mathcal{V}_2$ for any $v_1\in \mathcal V_1$,  $\mathcal{N}_{v_2} \subset \mathcal{V}_1$ for any $v_2\in \mathcal V_2$.   Furthermore, $\deg{(v_j)}=k_j$ for any $v_j\in \mathcal V_j$, $j=1, 2$.

 Any periodic graph $\Gamma$  can be viewed as an abelian covering graph of a finite graph $\Gamma_0$, which is called the {\it fundamental graph} of $\Gamma$.

\subsubsection{\bf Line graph}\label{subsublinegraph}
Given a graph $\Gamma = \{\mathcal V(\Gamma), \mathcal E(\Gamma)\}$, its {\it line graph} $L(\Gamma) = \{\mathcal V(L(\Gamma))$, $\mathcal E(L(\Gamma))\}$ is defined as follows : (1) The vertex set  $\mathcal V(L(\Gamma))$ is $\mathcal E(\Gamma)$. (2) $\mathcal E(L(\Gamma)) \ni (e_1, e_2)$, where $e_1, e_2 \in \mathcal E(\Gamma)$, if and only if $o(e_1) = o(e_2)$. 

The characteristic polynomials of $\Gamma$ and $L(\Gamma)$ are then related each other. In fact, let $\Gamma$ be a $k$-regular abelian covering graph of a finite graph $\Gamma_0$, $k\ge3$; $\mu$ and $\nu$ the numbers of the vertices and edges of $\Gamma_0$, respectively.
Let $\kappa=\nu-\mu$, which is a positive integer, and $p_{\Gamma}(x, \lambda)$ be the characteristic polynomial of $-\widehat \Delta_{\Gamma}$ on $\Gamma$.
Then $L(\Gamma)$ is a $2(k-1)$-regular abelian covering graph of $L(\Gamma_0)$ whose transformation group is also that of $\Gamma$, and the characteristic polynomial $p_{L(\Gamma)}(x, \lambda)$ of $-\widehat \Delta_{L(\Gamma )}$ is
\begin{equation}
  p_{L(\Gamma)}(x, \lambda)=\left(\frac{1}{k-1}-\lambda\right)^{\kappa}\left(\frac{k}{2k-2}\right)^{\mu}p_{\Gamma}\Big(x, \frac{2k-2}{k}(\lambda+\frac{k-2}{2k-2})\Big).
  \label{charpolynomlinegraph}
\end{equation}

 For instance, the Kagome lattice is a line graph of the hexagonal lattice.

We can also compute the characteristic polynomial of the line graph of a $(k_1, k_2)$-semiregular periodic graph, where $k_1 \ge k_2 \ge 3$ or $k_1 > k_2 = 2$. See \cite{S99} and \cite{HiShi04} for the details, where they study the spectrum of the discrete Laplacian on the line graph of $k$-regular or $(k_1, k_2)$-semiregular infinite graphs, which are not necessarily periodic.

\subsubsection{\bf Subdivision}\label{subsubSubdivision}
We can define the subdivisions of the triangular lattice, the hexagonal lattice, the Kagome lattice and the diamond lattice, and derive similar spectral properties for their Hamiltonians in the same way as in Subsection~\ref{SubsectionSubdivision}.
As in the case of line graph, the characteristic polynomials of a regular periodic graph and its subdivision are mutually related.
In fact, for a $k$-regular abelian covering graph $\Gamma$ of a finite graph $\Gamma_0$, $k \ge 3$, put $\mu$, $\nu$, and $\kappa$ in the same way as above, and let $p_{\Gamma}(x, \lambda)$ be the characteristic polynomial of $-\widehat \Delta_{\Gamma}$.
Then the subdivision $S(\Gamma)$ of $\Gamma$ is a $(k, 2)$-semiregular abelian covering graph of $S(\Gamma_0)$ whose transformation group is also that of $\Gamma$, and the characteristic polynomial $p_{S(\Gamma)}(x, \lambda)$ of $-\widehat \Delta_{S(\Gamma)}$ is
\begin{equation}
  p_{S(\Gamma)}(x, \lambda)=(-\lambda)^{\kappa}\left(-\frac{1}{2}\right)^{\mu}p_{\Gamma}(x, 1-2\lambda^2).
  \label{Rrmark3.13pxlambda}
\end{equation}
In \cite{S99} and \cite{HiShi04}, they also study the spectrum for the Laplacian on the subdivision of $k$-regular, $k\ge3$, infinite graphs, which are not necessarily periodic.

\subsubsection{\bf Ladder.}
The ladder structure is defined for any periodic graphs.
Let $\Gamma$ be a $k$-regular periodic graph, $\mu$ the number of vertices in the fundamental graph $\Gamma_0$, and $p_{\Gamma}(x, \lambda)$ the characteristic polynomial of $\widehat H_{\Gamma} = -\widehat \Delta_{\Gamma}$ on $\Gamma$.
The ladder $Lad(\Gamma)$ of $\Gamma$ is defined as a union of two copies of $\Gamma$ with additional edges joining the corresponding vertices, which is a $(k+1)$-regular periodic graph. Then $\widehat H_{Lad(\Gamma)} = -\widehat \Delta_{Lad(\Gamma)}$ has the following structure:
\begin{equation}
\widehat H_{Lad(\Gamma)}=\frac{k}{k+1}
\left(
  \begin{array}{cc}
    \widehat H_{\Gamma} & \dfrac{1}{k}I\\
    & \\
    \dfrac{1}{k}I & \widehat H_{\Gamma}
  \end{array}
\right): l^2(\Gamma) \oplus l^2(\Gamma) \to l^2(\Gamma) \oplus l^2(\Gamma).
\label{hamiltonian_ladder}
\end{equation}
Passing $\widehat H_{Lad(\Gamma)}$ to the Fourier transform, we have a multiplication operator by a $2\mu \times 2\mu$ symmetric matrix-valued function
\begin{equation}
H_{Lad(\Gamma)}(x)=\frac{k}{k+1}
\left(
  \begin{array}{cc}
    H_{\Gamma}(x) & \dfrac{1}{k}I\\
    & \\
    \dfrac{1}{k}I & H_{\Gamma}(x)
  \end{array}
\right) \ \text{ on } L^2(\mathbb{T}^d) \oplus L^2(\mathbb{T}^d),
\label{hamiltonian_ladder}
\end{equation}
where $H_{\Gamma}(x)$ is a $\mu \times \mu$ symmetric matrix-valued function which is obtained by passing $\widehat H_{\Gamma}$ to the Fourier transform.
The characteristic polynomial $p_{Lad(\Gamma)}(x, \lambda)$ of $\widehat H_{Lad(\Gamma)}$ is then computed as
\begin{equation}
  p_{Lad(\Gamma)}(x, \lambda)=\left(\frac{k}{k+1}\right)^{2\mu}p_{\Gamma}\Big(x, \frac{k+1}{k}\lambda+\frac{1}{k}\Big)\,p_{\Gamma}\Big(x, \frac{k+1}{k}\lambda-\frac{1}{k}\Big).
  \label{char_pol_ladder}
\end{equation}


\section{Distributions on the torus}


\subsection{Sobolev and Besov spaces on ${\bf R}^d$, lattice and torus}
Let $r_{-1}=0$, $r_j = 2^{j}$, $(j\geq0)$, and on ${\bf R}^d$, define the Besov space $\mathcal B({\bf R}^d)$ to be the set of all functions $f$ having the following norm :
\begin{equation}
\|f\|_{\mathcal B({\bf R}^d)} = \sum_{j=0}^{\infty}
r_j^{1/2}\Big(\int_{\Xi_j}|\widetilde f(\xi)|^2d\xi\Big)^{1/2},
\label{C2S1SpaceB}
\end{equation}
where   $\Xi_j = \{\xi \in {\bf R}^d\, ; \, 
r_{j-1} \leq |\xi| < r_j\}$, and $\widetilde f$ denotes the Fourier transform.
The (equivalent) norm of the dual space $\mathcal B^{\ast}({\bf R}^d)$ is
\begin{equation}
\|u\|_{\mathcal B^{\ast}({\bf R}^d)} = \Big(\mathop{\rm sup}_{R>1}
\frac{1}{R}\int_{|\xi|<R}|\widetilde u(\xi)|^2d\xi\Big)^{1/2}.
\label{S4BastnormRd}
\end{equation}
The space $\mathcal B_0^{\ast}({\bf R}^d)$ is defined as follows:
\begin{equation}
\mathcal B^{\ast}_0({\bf R}^d) = \Big\{u\in \mathcal B^{\ast}({\bf R}^d)\, ; \, 
\lim_{R\to\infty}\frac{1}{R}\int_{|\xi|<R}|\widetilde u(\xi)|^2d\xi = 0\Big\}.
\label{S4Bast0Def}
\end{equation}
The Sobolev space $H^{\sigma}({\bf R}^d)$ is defined as usual : 
\begin{equation}
H^{\sigma}({\bf R}^d) = \left\{u \in \mathcal S'({\bf R}^d)\, ; \, 
\|(1 + |\xi|^2)^{\sigma/2}\widetilde u(\xi)\|_{L^2({\bf R}^d)} < \infty\right\}, \quad \sigma \in {\bf R}.
\label{S4SobolevRdDef}
\end{equation}

The Besov spaces $\mathcal B$, $\mathcal B^{\ast}$ are also defined on ${\bf T}^d$. 
Take a $C^{\infty}$-partition of unity $\{\chi_{\ell}\}_{\ell=1}^N$ on ${\bf T}^d$ where the support of $\chi_{\ell}$ is sufficiently small, and define
\begin{equation}
\|f\|_{\mathcal B({\bf T}^d)} = \sum_{\ell=1}^N\|\chi_{\ell}f\|_{\mathcal B({\bf R}^d)},
\label{BonTd}
\end{equation}
\begin{equation}
\|u\|_{\mathcal B^{\ast}({\bf T}^d)} = \sum_{\ell=1}^N\|\chi_{\ell}u\|_{\mathcal B^{\ast}({\bf R}^d)}.
\label{BastonTd}
\end{equation}
The space $\mathcal B_0^{\ast}({\bf T}^d)$ is defined to be the set of $u \in \mathcal B^{\ast}({\bf T}^d)$ such that
\begin{equation}
\chi_{\ell } u \in \mathcal B_0^{\ast}({\bf R}^d), \quad 1 \leq \ell \leq N.
\label{Bast0Td}
\end{equation}
The Sobolev space $H^{\sigma}({\bf T}^d)$ is defined similarly.

The analogues of $\mathcal B$ and $\mathcal B^{\ast}$ on the lattice ${\bf Z}^d$ are defined  to be the
 Banach spaces endowed with norms
\begin{equation}
\|\widehat f\|_{\widehat{\mathcal B}({\bf Z}^d)} 
= 
\sum_{j=0}^{\infty}
r_j^{1/2}\Big(\sum_{r_{j-1} \leq |n| < r_j}|\widehat f(n)|^2\Big)^{1/2},
\end{equation}
\begin{equation}
\|\widehat u\|_{\widehat{\mathcal B}^{\ast}({\bf Z}^d)} = \Big(
\sup_{R>1}\frac{1}{R}\sum_{|n|<R}|\widehat u(n)|^2\Big)^{1/2}.
\end{equation}
The space $\widehat{\mathcal B}^{\ast}_0({\bf Z}^d)$ is then defined by
\begin{equation}
\widehat{\mathcal B}^{\ast}_0 ({\bf Z}^d)= \Big\{ \widehat u \in \widehat{\mathcal B}^{\ast}({\bf Z}^d)\, ; \, \lim_{R\to\infty}\frac{1}{R}\sum_{|n|<R}|\widehat u(n)|^2 = 0\Big\}.
\end{equation}
For $\sigma \in {\bf R}$, the weighted space $\ell^{2,\sigma}$  is the set of all $\widehat u$ satisfying
$$
\|\widehat u\|^2_{\ell^{2,\sigma}} = \sum_{n \in {\bf Z}^d}(1 + |n|^2)^{\sigma}|\widehat u(n)|^2 < \infty.
$$

These Besov spaces $\mathcal B({\bf T}^d)$, $\widehat{\mathcal B}({\bf Z}^d)$, $\mathcal B^{\ast}({\bf T}^d)$, $\widehat{\mathcal B}^{\ast}({\bf Z}^d)$ are related  by the  Fourier series.  For $\widehat f \in \mathcal S'({\bf Z}^d)$, let $f = \mathcal U\widehat f$, where $\mathcal U$ is defined by (\ref{S1Fourierseries}).
We define operators $\widehat N_j$ on the lattice ${\bf Z}^d$ and $N_j$ on the torus ${\bf T}^d$ by
$$
\big(\widehat N_j\widehat f)(n) = n_j\widehat f(n),
\quad
N_j = \mathcal U\widehat N_j\mathcal U^{\ast} = i\frac{\partial}{\partial x_j}.
$$
We put $N = (N_1,\cdots,N_d)$, and let $N^2$ be the self-adjont operator defined by
\begin{equation}
N^2 = \sum_{j=1}^dN_j^2 =  - \Delta, \quad {\rm on} \quad {\bf T}^d,
\nonumber
\end{equation}
where $\Delta$ is the Laplacian on ${\bf T}^d$ with periodic boundary condition. We put
\begin{equation}
|N| = \sqrt{N^2} = \sqrt{-\Delta}.
\label{S4Define|N|}
\end{equation}
For $\sigma \in {\bf R}$, $H^{\sigma}({\bf T}^d)$ coincides with the completion of $D(|N|^{\sigma})$ with respect to the norm
$\|u\|_{\sigma} = \|(1 + N^2)^{\sigma/2}u\|$ i.e.
\begin{equation}
H^{\sigma}({\bf T}^d) = \big\{u \in {\mathcal S'}({\bf T}^d)\, ; \, \|u\|_{\sigma} = \|(1 - \Delta)^{\sigma/2}u\| <
\infty \big\}.
\nonumber
\end{equation}

For a self-adjoint operator $T$, let $\chi(a \leq T < b)$ denote the operator $\chi_{I}(T)$, where $\chi_I(\lambda)$ is the characteristic function of the interval $I = [a,b)$. The operators $\chi(T < a)$ and $\chi(T \geq b)$ are defined similarly.  Then the Besov spaces $\mathcal B({\bf T}^d), \mathcal B^{\ast}({\bf T}^d)$ are rewritten by the 
equivalent norms :
\begin{equation}
\mathcal B({\bf T}^d) = \Big\{f \in L^2({\bf T}^d)\, ; \|f\|_{\mathcal B({\bf T}^d)} = \sum_{j=0}^{\infty}r_j^{1/2}\|\chi(r_{j-1} \leq |N| < r_j)f\| < \infty\Big\},
\label{S4BTdequivnorm}
\end{equation}
\begin{equation}
\mathcal B^{\ast}({\bf T}^d) = \Big\{ u \in \mathcal S'({\bf T}^d)\, ; \, \|u\|_{\mathcal B^{\ast}({\bf T}^d)}= \left(\sup_{R>1}\frac{1}{R}\|\chi(|N|<R)u\|^2\right)^{1/2}<\infty\Big\},
\label{S4BastTdequivnorm}
\end{equation}
and $\mathcal B^{\ast}_0({\bf T}^d)$ is rewritten as
\begin{equation}
\mathcal B^{\ast}_0({\bf T}^d) = \left\{u \in \mathcal B^{\ast}({\bf T}^d)\, ; \, 
\lim_{R\to\infty}\frac{1}{R}\|\chi(|N|<R)u\|^2=0\right\}.
\label{S4Bast0TdOpform}
\end{equation}

In fact, the equivalence of (\ref{BastonTd}) and (\ref{S4BastTdequivnorm}), and that of (\ref{Bast0Td}) and (\ref{S4Bast0TdOpform}) are proved in \cite{IsMo1}, Lemmas 3.1 and 3.2.
 The equivalence of (\ref{BonTd}) and (\ref{S4BTdequivnorm}) follows from this by duality.

Note also the following equivalence:
\begin{eqnarray*}
f \in \mathcal B({\bf T}^d) &\Longleftrightarrow&  \widehat f \in \widehat{\mathcal B}({\bf Z}^d), \\
u \in \mathcal B^{\ast}({\bf T}^d) &\Longleftrightarrow & \widehat u \in \widehat{\mathcal B}^{\ast}({\bf Z}^d), \\
u \in \mathcal B^{\ast}_0({\bf T}^d) &\Longleftrightarrow & \widehat u \in \widehat{\mathcal B}^{\ast}_0({\bf Z}^d). 
\end{eqnarray*}
In the sequel, we often write $\mathcal B$, $\widehat{\mathcal B}$, $\mathcal B^{\ast}$, $\widehat{\mathcal B}^{\ast}$, $\mathcal B^{\ast}_0$ and
 $\widehat{\mathcal B}^{\ast}_0$, omitting ${\bf R}^d$, ${\bf T}^d$, ${\bf Z}^d$.


\subsection{Basic Lemmas}
We use the following properties of Besov spaces. Let $S^m_{1,0}$ be the standard H{\"o}rmander class of $\Psi$DO on ${\bf R}^d$.


\begin{lemma}\label{PsDOandBspace}
(1) If $f \in {\mathcal B}({\bf R}^d)$, we have
$$
\int_{-\infty}^{\infty}\|f(x_1,\cdot)\|_{L^2({\bf R}^{d-1})}dx_1 \leq 
\sqrt2 \|f\|_{\mathcal B({\bf R}^d)}.
$$
(2) If $P \in S^0_{1,0}$, we have
$$
P \in {\bf B}(\mathcal B;\mathcal B)\cap{\bf B}({\mathcal B}^{\ast};{\mathcal B}^{\ast})\cap{\bf B}(\mathcal B^{\ast}_0;{\mathcal B}^{\ast}_0).
$$
\end{lemma}

Proof. The assertion (1) is proven in \cite{HoVol2}, Theorem 14.1.2. 
For $\chi \in C_0^{\infty}({\bf R}^d)$, let $\chi_R$ be the $\Psi$DO with symbol $\chi(\xi/R)$. Then we have
$$
u \in \mathcal B^{\ast} \Longleftrightarrow \sup_{R>1}\frac{1}{\sqrt{R}}
\|\chi_Ru\|_{L^2({\bf R}^d)} < \infty, \quad 
\forall \chi \in C_0^{\infty}({\bf R}^d).
$$
This and the symbolic calculus of $\Psi$DO imply $P \in {\bf B}({\mathcal B}^{\ast};{\mathcal B}^{\ast})$. Taking the adjoint, we have $P \in {\bf B}({\mathcal B};{\mathcal B})$. The fact $P \in {\bf B}({\mathcal B}^{\ast}_0;{\mathcal B}^{\ast}_0)$ is proven similarly. \qed


\begin{lemma}\label{C3S3u0estimatedbylinsup}
Suppose $u \in \mathcal S'({\bf R}^d)$ satisfies $\widetilde{u} \in L^2_{loc}({\bf R}^d)$ and
$$
\limsup_{R\to\infty}\frac{1}{R}\int_{|\xi|<R}|\widetilde{u} (\xi)|^2d\xi < \infty. 
$$
If there is a submanifold $M$ of codimension 1 in ${\bf R}^d$ such that ${\rm supp}\,u \subset M$, then there exists $u_0 \in L^2(M)$ such that
\begin{equation}
\langle u,\varphi\rangle = \int_Mu_0\varphi\,  dM, \quad \forall \varphi \in \mathcal S({\bf R}^d),
\label{C3S2u=u0onM}
\end{equation}
\begin{equation}
\int_M|u_0|^2dM \leq C\limsup_{R\to\infty}\frac{1}{R}\int_{|\xi|<R}|\widetilde{u} (\xi)|^2d\xi < \infty.
\label{C3S2u0estimate}
\end{equation} 
\end{lemma}

For the proof, see \cite{HoVol1}, Theorem 7.1.27. 


\subsection{Distribution $(x_1 \mp i0)^{-1}$ and its wavefront set}
We need a division theorem and its micro-local consequences.
 Let us begin with the case of ${\bf R}^d$.


\begin{lemma}\label{C3S3u+lemma}
For $f \in \mathcal B({\bf R}^d)$ and $\epsilon > 0$, let $u_{\epsilon}(x) = f(x)/(x_1 - i\epsilon)$. Then,
$\lim_{\epsilon\to0}u_{\epsilon} = f(x)/(x_1-i0) =:u_+$ exists in the weak $\ast$ sense, i.e. 
\begin{equation}
(u_{\epsilon},g) \to (u_+,g), \quad \forall g \in \mathcal B({\bf R}^d), 
\label{C2S4uepsiontou+}
\end{equation}
with the following estimate
\begin{equation}
\|u_+\|_{\mathcal B^{\ast}({\bf R}^d)} \leq 2\|f\|_{\mathcal B({\bf R}^d)}.
\label{C4S3u+isdominatedbyf}
\end{equation}
Moreover, $\widetilde u_+(\xi)$ is an $L^2({\bf R}^{d-1})$-valued bounded function of $\xi_1$ and 
\begin{equation}
\|\widetilde u_+(\xi_1,\cdot)\|_{L^2({\bf R}^{d-1})} \to 0, \quad {\rm as}\quad \xi_1 \to \infty,
\label{Lemma4.3xi1toinfty}
\end{equation}
\begin{equation}
\Big\|\widetilde u_+(\xi_1,\cdot) - i\int_{-\infty}^{\infty}\widetilde f(\eta_1,\cdot)d\eta_1\Big\|_{L^2({\bf R}^{d-1})} \to 0 , \quad {\rm as} \quad \xi_1\to-\infty.
\label{C4S3u+limitxitoinfty}
\end{equation}
\end{lemma}
Proof. 
Letting $\theta$ be the Heaviside function, we have
\begin{equation}
(u_{\epsilon},g) = i\int _{{\bf R}^{d+1}}\theta(\eta_1-\xi_1)e^{\epsilon(\xi_1-\eta_1)}\widetilde f(\eta_1,\xi')\overline{\widetilde g(\xi)}\,d\eta_1d\xi.
\label{S4Lemma4.1Formula}
\end{equation}
By the Schwarz inequality,
$$
\int_{{\bf R}^{d-1}}|\widetilde f(\eta_1,\xi')\overline{\widetilde g(\xi_1,\xi')}|d\xi' \leq 
\|\widetilde f(\eta_1,\cdot)\|_{L^2({\bf R}^{d-1})}\|\widetilde g(\xi_1,\cdot)\|_{L^2({\bf R}^{d-1})}.
$$
 Using Lemma \ref{PsDOandBspace} (1), we have
\begin{equation}
|(u_{\epsilon},g)| \leq 2\|f\|_{\mathcal B}\|g\|_{\mathcal B},
\label{C4S3uepsilonlessthanfg}
\end{equation}
which implies
$\|u_{\epsilon}\|_{{\mathcal B}^{\ast}} \leq 2\|f\|_{\mathcal B}$. 
If $f, g \in C_0^{\infty}({\bf R}^d)$, $u_{\epsilon} \to u_+$ pointwise, and 
$(u_{\epsilon},g) \to (u_+,g)$. By (\ref{C4S3uepsilonlessthanfg}), we have 
$|(u_+,g)|\leq 2\|f\|_{\mathcal B}\|g\|_{\mathcal B}$, hence $u_+ \in \mathcal B^{\ast}$. Approximating $f$ and $g$ by elements of $C_0^{\infty}({\bf R}^d)$, we see that (\ref{C2S4uepsiontou+}) and (\ref{C4S3u+isdominatedbyf}) hold for $f, g \in \mathcal B$.
The equality (\ref{S4Lemma4.1Formula}) shows
\begin{equation}
\widetilde u_+(\xi) = i\int_{\xi_1}^{\infty}\widetilde f(\eta_1,\xi')d\eta_1.
\label{Lemma4.3u+intf(xi)}
\end{equation}
Lemma \ref{PsDOandBspace} (1) then  implies that $\widetilde u_+(\xi)$ is an $L^2({\bf R}^{d-1})$-valued bounded function of $\xi_1$, and (\ref{Lemma4.3xi1toinfty}).
Moreover, we have
$$
\|\widetilde u_+(\xi_1,\cdot) - i \int_{-\infty}^{\infty}
\widetilde f(\eta_1,\cdot)d\eta_1\|_{L^2({\bf R}^{d-1})} 
\leq \int_{-\infty}^{\xi_1}\|\widetilde f(\eta_1,\cdot)\|_{L^2({\bf R}^{d-1})}d\eta_1,
$$
which yields (\ref{C4S3u+limitxitoinfty}). \qed


\begin{definition}
For $u \in \mathcal S'({\bf R}^d)$, the {\it  wave front set } $WF^{\ast}(u)$ is defined as follows. For $(x_0,\omega) \in {\bf R}^d\times S^{d-1}$,  
$(x_0,\omega) \not\in WF^{\ast}(u)$, if there exist $0 < \delta < 1$ and $\chi \in C^{\infty}_0({\bf R}^d)$ satisfying $\chi(x_0)=1$ such that
\begin{equation}
\lim_{R\to\infty}\frac{1}{R}\int_{|\xi|<R}\big|C_{\omega,\delta}(\xi)(\widetilde{\chi u})(\xi)\big|^2d\xi = 0,
\label{C4S3ScatteringWFdefine}
\end{equation}
 where 
$C_{\omega,\delta}(\xi)$ is the characteristic function of the cone 
$\big\{\xi \in {\bf R}^d\, ; \,  \omega\cdot\xi> \delta|\xi|\big\}$.
\end{definition}


\begin{lemma}\label{C3S4u+wavefrontLemma}
For $u_+$ in Lemma \ref{C3S3u+lemma}, we have
\begin{equation}
WF^{\ast}(u_+) \subset \left\{\big((0,x'),(-1,0)\big)\, ; \, x' \in {\bf R}^{d-1}\right\}.
\label{C3S4u+wavefront}
\end{equation}
\begin{equation}
u_+(x) - \frac{1}{x_1-i0}\otimes f(0,x') \in \mathcal B^{\ast}_0({\bf R}^d).
\label{C3S4u+expansion}
\end{equation}
\end{lemma}

Proof. Take $h\in C_0^{\infty}({\bf R}^d)$ and put $w(x) = h(x)/(x_1-i0)^{-1}$. Then by (\ref{C4S3u+isdominatedbyf}), 
\begin{equation}
\frac{1}{R}\int_{|\xi|<R}|C_{\omega,\delta}(\xi)\widetilde{\chi u_+}(\xi)|^2d\xi
\leq  \frac{1}{R}\int_{|\xi|<R}|C_{\omega,\delta}(\xi)\widetilde{\chi w}(\xi)|^2d\xi + C{\|f-h\|^2}_{\mathcal B},
\nonumber
\end{equation}
where the constant $C$ is independent of $R>1$. Therefore, we have only to prove the lemma when $f \in C_0^{\infty}({\bf R}^d)$. Obviously, $\big((y_1,y'),\omega\big) \not\in WF^{\ast}(u_+)$ if $y_1 \neq 0$. Take $\big((0,y'),\omega\big)$ such that $(-1,0) \neq \omega \in S^{d-1}$. Take $\chi(x) \in C_0^{\infty}({\bf R}^d)$ such that $\chi((0,y'))=1$, and put $v(x) = \chi(x)u_+(x), \; g(x) = \chi(x)f(x)$. Then $v(x) = g(x)/(x_1-i0)$, hence by passing to the Fourier transform,
$$
\widetilde v(\xi_1,\xi') = i\int_{\xi_1}^{\infty}\widetilde g(\eta_1,\xi')d\eta_1.
$$
This implies
$$
|\widetilde v(\xi)| \leq C_N\int_{\xi_1}^{\infty}(1 + |\eta_1| + |\xi'|)^{-N}d\eta_1, \quad \forall N>0.
$$
Since $\omega \neq (-1, 0, \cdots,0)$, by taking $0<\delta<1$ sufficiently close to 1, 
on the region $\{\omega\cdot \xi >\delta|\xi|\}$, we have either $C|\xi_1|\leq |\xi'|$, or $\xi_1 > 0, \ C|\xi'| \leq \xi_1$, where $C>0$. In both case, we obtain $|\widetilde v(\xi)|\leq C_N(1 + |\xi|)^{-N}$, which proves (\ref{C3S4u+wavefront}).

Passing to the Fourier transform, $u_+(x) = f(x)/(x_1-i0)$ becomes $\widetilde u_+(\xi) = i\int_{\xi_1}^{\infty}\widetilde f(\eta_1,\xi')d\eta_1$. Letting $\theta(t)$ be the Heaviside function, we then have
\begin{equation}
 \widetilde u_+(\xi) - i\theta(-\xi_1)\int_{-\infty}^{\infty}\widetilde f(\eta_1,\xi')d\eta_1
	  = \left\{
\begin{split}
i\int_{\xi_1}^{\infty}\widetilde f(\eta_1,\xi')d\eta_1 \quad {\rm if} \quad \xi_1 > 0, \\
-i\int_{-\infty}^{\xi_1}\widetilde f(\eta_1,\xi')d\eta_1 \quad {\rm if} \quad \xi_1 < 0.
\end{split}
\right.
\nonumber
\end{equation}
We then have
$$
\|\widetilde u_+(\xi) - i\theta(-\xi_1)\int_{-\infty}^{\infty}\widetilde f(\eta_1,\xi')d\eta_1\|_{L^2({\bf R}^{d-1})} \to 0, \quad 
{\rm as} \quad |\xi_1| \to \infty.
$$
This and the inequality
\begin{equation}
\begin{split}
& \frac{1}{R}\int_{|\xi|<R}\big|\widetilde u_+(\xi) - i\theta(-\xi_1)\int_{-\infty}^{\infty}\widetilde f(\eta_1,\xi')d\eta_1\big|^2d\xi \\
& \leq \frac{1}{R}\int_{|\xi_1|<R}\|\widetilde u_+(\xi_1,\cdot) - i\theta(-\xi_1)\int_{-\infty}^{\infty}\widetilde f(\eta_1,\cdot)d\eta_1\|^2_{L^2({\bf R}^{d-1})}d\xi_1 
\end{split}
\nonumber
\end{equation}
yield (\ref{C3S4u+expansion}). \qed

\medskip
Consider the equation
\begin{equation}
x_1u = f, \quad f \in \mathcal B.
\label{C3S3Eqxiu=f}
\end{equation}
A solution $u_+ \in \mathcal B^{\ast}$ ($u_- \in \mathcal B^{\ast}$, respectively), of the equation (\ref{C3S3Eqxiu=f}) is said to be {\it outgoing} ({\it incoming}) if it satisfies
\begin{equation}
WF^{\ast}(u_{+}) \subset \{((0,x'),(- 1, 0))\, ; \, x' \in {\bf R}^{d-1}\},
\end{equation}
\begin{equation}
WF^{\ast}(u_{-}) \subset \{((0,x'),(1, 0))\, ; \, x' \in {\bf R}^{d-1}\}.
\end{equation}


\begin{lemma}\label{LemmaIm(uf)onRd}
Let $u \in \mathcal B^{\ast}$ be a solution to the equation (\ref{C3S3Eqxiu=f}). Then $u$ is outgoing (incoming) if and only if
$$
u = \frac{f(x)}{x_1 - i0},  \quad \left(u = \frac{f(x)}{x_1 + i0}\right).
$$
For outgoing (incoming) solution $u_+$ ($u_-$), we have
\begin{equation}
{\rm Im}\, (u_{\pm},f) = \pm \pi \|f(0,\cdot)\|^2_{L^2({\bf R}^{d-1})}.
\label{S4Impmuf}
\end{equation}
\end{lemma}

Proof. The "if" part is proven in Lemma \ref{C3S4u+wavefrontLemma}. To prove the "only if" part, let $u$ be an outgoing solution and 
$v = u - f(x)/(x_1-i0)$. Then $v \in \mathcal B^{\ast}$ is an outgoing solution 
to $x_1v=0$. Passing to the Fourier transform, this implies that $\widetilde v(\xi)$ depends only on $\xi'$ : $\widetilde v(\xi)= w(\xi')$. Integrating over the region 
 $\Omega_R=\{(\xi_1,\xi') \; 0 < \xi_1<R, |\xi'|<R/2\}$, we have
 $$
 \int_{|\xi'|<R/2}|w(\xi')|^2d\xi \leq \frac{1}{R}\int_{\Omega_R}
 |\widetilde v(\xi)|^2d\xi.
 $$
Letting $R \to \infty$, we have $v=0$, which shows that $u = f(x)/(x_1-i0)$.
The well-known formula
$$
\frac{1}{t\mp i0} = \pm i\pi\delta(t) + {\rm p.v.}\frac{1}{t}
$$
implies (\ref{S4Impmuf}). \qed

\subsection{Distribution $(h(x) \mp i0)^{-1}$ on ${\bf T}^d$}
We now consider the equation
\begin{equation}
\big(h(x)  - z\big)u(x) = f(x), \quad {\rm on} \quad 
{\bf T}^d,
\label{C4S3equationonTd}
\end{equation}
where $h(x)$ is a real-valued $C^{\infty}$-function on ${\bf T}^d$. We put
$$
M = \{x \in {\bf T}^d\, ; \, h(x) = 0\},
$$
and  assume that 

\medskip
\noindent
(C-1) $\nabla h(x) \neq 0$ on $M$.

\medskip
Take $x_0 \in M$, $\chi \in C^{\infty}({\bf T}^d)$ such that $\chi(x_0)=1$ and the support of $\chi$ is sufficiently small. We make a change of variable $ x \to y$ around $x_0$, where $y_1 = h(x)$. Letting $v(y) = \chi(x)u(x)$, $F(y) = \chi(x)f(x)$, we then have
$$
(y_1 - z)v(y) = F(y), \quad z \not\in {\bf R}.
$$
By $T_{x_0}(M)^{\perp}$, we mean the orthogonal compliment of $T_{x_0}(M)$ in $T_{x_0}({\bf R}^d)$. 
Applying Lemmas \ref{C3S3u+lemma}, \ref{C3S4u+wavefrontLemma}, we obtain the following lemma. 


\begin{lemma}\label{C3S4h(x)u=flemma}
Let $u_{z} = f(x)/(h(x) - z)$, $z \not\in {\bf R}$. Then, there exists a limit $\lim_{\epsilon\to0} u_{\pm i\epsilon}=:u_{\pm}$ in  the
weak $\ast$ sense, i.e.
\begin{equation}
(u_{\pm i \epsilon},g) \to (u_{\pm},g), \quad \forall g \in \mathcal B.
\label{C3S4limiupmonTd}
\end{equation}
Moreover,
\begin{equation}
\|u_{\pm}\|_{\mathcal B^{\ast}} \leq C\|f\|_{\mathcal B},
\label{C3S4upmastonTd}
\end{equation}
\begin{equation}
WF^{\ast}(u_{\pm}) \subset \{ (x, \pm \omega_x)\, ; \, x \in M\},
\label{C3S4WFupmonTd}
\end{equation}
where $\omega_x \in S^{d-1}\cap T_x(M)^{\perp}$, and $ \omega\cdot \nabla h(x) < 0$,
\begin{equation}
u_{\pm}(x) - \frac{1}{h(x)\mp i0}\otimes \left(f\big|_M\right) \in 
\mathcal B^{\ast}_0,
\label{C3S4upmsingexpansion}
\end{equation}
where $f\big|_M$ means the restriction of $f$ to $M$.
\end{lemma}

A solution $u \in \mathcal B^{\ast}$ of the equation 
\begin{equation}
h(x)u= f(x) \in \mathcal B
\label{C3S4hxu=f}
\end{equation}
is said to be outgoing (incoming) if it satisfies
\begin{equation}
WF^{\ast}(u) \subset \{ (x,  \omega_x)\, ; \, x \in M\}, \quad 
\Big(WF^{\ast}(u) \subset \{ (x,  -\omega_x)\, ; \, x \in M\}\Big),
\end{equation}
$\omega_x$ being as a above. The following lemma is a direct consequence of Lemma \ref{LemmaIm(uf)onRd}.


\begin{lemma}\label{S4Lemma47}
Let $u \in \mathcal B^{\ast}$ be a solution to the equation (\ref{C3S4hxu=f}). Then $u$ is outgoing (incoming) if and only if
$$
u = \frac{f(x)}{h(x) - i0},  \quad \left(u = \frac{f(x)}{h(x) + i0}\right).
$$
For outgoing (incoming) solution $u_+$ ($u_-$), we have
$$
{\rm Im}\, (u_{\pm},f) = \pm \pi \big\|f\big|_{M}\big\|^2_{L^2(M)}.
$$
In particular, we have
$$
\frac{1}{2\pi i}(u_+ - u_-,f) =  \big\|f\big|_{M}\big\|^2_{L^2(M)}.
$$
\end{lemma}


\section{Rellich type theorem on the torus}


\subsection{Rellich type theorem}
Let $\mathcal H_0 = L^2\left({\bf T}^d\right)^s$  equipped with the inner product (\ref{S2L2bfTdinnerproduct}), and $H_0(x)$ an $s \times s$ hermitian matrix whose entries are polynomials of $e^{\pm ix_1}$, $\cdots$, $e^{\pm ix_d}$. Then, the operator of multiplication by $H_0(x)$ is a bounded self-adjoint operator on $\mathcal H_0$, which is denoted by $H_0$.
The spectrum of $H_0$ is given by (\ref{S2DefineSigma(H0)}).

 We are going to study the following theorem for the  Hamiltonian on the lattice (see \cite{IsMo1}) 
We define $\widehat H_0 = \mathcal U^*_{\mathcal{L}_0} H_0\, \mathcal U_{\mathcal{L}_0} = - \widehat \Delta_{\Gamma_0}$, and suppose $\widehat u$ satisfies for some $R_0>0$  and $\lambda \in \sigma(\widehat H_0)$  \textit{except for some exceptional points} to be defined below
$$
(\widehat H_0 - \lambda)\widehat u = 0 \ \ {\rm for} \ \ |n| > R_0, \qquad 
\lim_{R\to\infty}\frac{1}{R}{\mathop\sum_{R_0 < |n| < R}}|\widehat u(n)|^2 = 0.
$$
Then, there exists $R > R_0$ such that $\widehat u(n) = 0$ for $|n| > R$.

\medskip
 This theorem is proven by passing to the torus ${\bf T}^d$, and requires the following property on the regular part of the \textit{complex Fermi surface}
 $M_{\lambda}^{\bf C} = \{z \in {\bf T}^d_{\bf C}\, ; \, p(z,\lambda) = 0\}$.
Recall the notation (\ref{S2DefineMlambdaCreg}), (\ref{S2DefineMlambdaCsng}). We assume 

\bigskip
\noindent
{\bf (A-1)} {\it There exists a  subset $\mathcal T_1 \subset \sigma(H_0)$ such that for $\lambda \in \sigma(H_0)\setminus\mathcal T_1$ : }

\smallskip
{\bf (A-1-1)}  $M_{\lambda,sng}^{\bf C}${\it is discrete.}

\smallskip
{\bf (A-1-2)} {\it Each connected component of $M_{\lambda,reg}^{\bf C}$ intersects with ${\bf T}^d$ and the intersection is a $(d-1)$-dimensional real analytic submanifold of ${\bf T}^d$.}

\bigskip
We reformulate the Rellich type theorem in the following way. A {\it trigonometric polynomial} is a vector function, each component of which has the form
$\sum_{|\alpha| \leq N}c_{\alpha}e^{i\alpha\cdot x}$,
 $c_{\alpha}$ being a constant. 


\begin{theorem}\label{C4S2RellichTh}
Assume {\bf (A-1)}, and let  $\lambda \in \sigma(H_0)\setminus{\mathcal T_1}$. 
Then if $u \in \mathcal B^{\ast}_0({\bf T}^d)$ satisfies
\begin{equation}
(H_0 - \lambda)u = f(x) \quad {\rm on} \quad {\bf T}^d,
\label{C4S2EquationH0-lambdau=f}
\end{equation}
for some trigonometric polynomial $f(x)$, $u(x)$ is also a trigonometric polynomial.
\end{theorem}

The 1st step of the proof of Theorem \ref{C4S2RellichTh} is to multiply the equation (\ref{C4S2EquationH0-lambdau=f}) by the cofactor matrix of $H_0(x)-\lambda$ and transform it as
\begin{equation}
p(x,\lambda) u(x) = g(x),
\label{C4S2pu=g}
\end{equation}
where $g(x)$ is a trigonometric polynomial.
In the following, we pick up one of the components of $u$ and $g$, and denote them by $u$ and $g$ again.


\begin{lemma}
Let $\lambda$ and $u$ be as in Theorem \ref{C4S2RellichTh}.
Then $u \in C^{\infty}({\bf T}^d\setminus M_{\lambda,sng}^{\bf C})$. In particular, we have
\begin{equation}
g(x)=0 \quad {\rm  on} \quad M_{\lambda,reg}^{\bf C}\cap {\bf T}^d.
\label{C4S2g(x)=0}
\end{equation} 
\end{lemma}

Proof. Take $x^{(0)}\in M_{\lambda,reg}^{\bf C}\cap{\bf T}^d$, and let $U$ be a sufficiently small neighborhood of $x^{(0)}$ in ${\bf T}^d$ such that 
$U \cap M_{\lambda,sng}^{\bf C} = \emptyset$. Take $\chi \in C^{\infty}({\bf T}^d)$ satisfying ${\rm supp}\,\chi \subset U$ and $\chi(x^{(0)})=1$.
Since $\nabla p(x_0,\lambda)\neq0$,we can make the change of variables on $U$ : $x \to y=(y_1,y')$ so that $y_1 = p(x,\lambda)$. Let $ v =\chi(x)u(x) $, $h=\chi (x)g(x)$,
 Since $u \in \mathcal B^{\ast}_0$, passing to the Fourier transform,
\begin{equation}
\lim_{R \to \infty} \frac{1}{R} \int_{|\eta |<R} |\widetilde v (\eta )|^2 d\eta =0.
\label{rellich_asymptoticzero2} 
\end{equation}
By (\ref{C4S2pu=g}),  $\frac{\partial}{\partial \eta _1} \widetilde{v} (\eta)= i \widetilde{h} (\eta )$. 
Integrating this equation, we have 
$$
\widetilde v(\eta) = i\int_0^{\eta_1}\widetilde h(s, \eta' )ds + \widetilde v(0, \eta' ).
$$
Since $\widetilde h(\eta)$ is rapidly decreasing, there exists the limit
$$
\lim_{\eta_1\to\infty}\widetilde v(\eta) = i \int_0^{\infty}\widetilde h(s, \eta' )ds + \widetilde v(0, \eta' ).
$$
We show that this limit vanishes. 
Let $D_R $ be the slab such that 
\begin{equation*}
D_R = \left\{ \eta \ ; \ |\eta '|<\delta R , \ \frac{R}{3 } < \eta_1 < \frac{2R}{3} \right\}. 
\end{equation*}
Then we have $D_R \subset \{ |\eta |<R \} $ for a sufficiently small $\delta >0 $.
We then see that 
\begin{equation*}
\frac{1}{R} \int_{ D_R } | \widetilde{v} (\eta )|^2 d\eta =\frac{1}{R} \int_{|\eta '|<\delta R} \int_{R/3}^{2R/3} |\widetilde{v} (\eta_1,\eta' )|^2 d\eta_1 d\eta'  \leq \frac{1}{R} \int_{ |\eta |<R} |\widetilde{v} ( \eta )|^2 d\eta. 
\end{equation*}
As $R\rightarrow \infty $, the right-hand side tends to zero by (\ref{rellich_asymptoticzero2}), hence so does the left-hand side, which proves that 
$\lim_{\eta_1\to\infty}\widetilde v(\eta) = 0$.
We have, therefore, 
\begin{equation}
\widetilde{v} (\eta )= -i\int_{\eta_1 } ^{\infty } \widetilde{h} (s ,\eta ')ds.
\nonumber
\end{equation}
Then $\widetilde v(\eta)$ is rapidly decreasing as $\eta_1 \to \infty$. 
Similarly, $\widetilde v(\eta)$ is rapidly decreasing as $\eta_1 \to - \infty$.
Therefore, $v=\chi u \in C^{\infty }({\bf T}^d)$.

It is easy to see that $\chi u$ is smooth outside $M_{\lambda }$.
We have thus proven  $u \in C^{\infty}(M_{\lambda,reg}^{\bf C}\cap {\bf T}^d)$. In particular, $g(x)= 0$ on $M_{\lambda,reg}^{\bf C}\cap {\bf T}^d$.
\qed



\begin{lemma}\label{gz=0onMlambdaC}
$g(z) = 0$ on $M_{\lambda,reg}^{\bf C}$.
\end{lemma}

Proof. Near any point in $M_{\lambda,reg}^{\bf C}\cap {\bf T}^d$, one can take local coordinates $(\zeta_1,\cdots,\zeta_d)$ so that $M_{\lambda,reg}^{\bf C}$ is represented as $\zeta_d=0$. Let $\zeta_j = s_j + it_j$. We expand $g\big|_{M_{\lambda,reg}^{\bf C}}$ into a Taylor series: 
$$
g\big|_{M_{\lambda,reg}^{\bf C}} = 
\sum c_{n_1\cdots n_{d-1}}\zeta_1^{n_1}\cdots\zeta_{d-1}^{n_{d-1}},
$$
which vanishes for $t_1 = \cdots=t_{d-1}=0$. We then have $c_{n_1\cdots n_{d-1}} =0$, hence $g(z)$ vanishes in a neighborhood of $M_{\lambda,reg}^{\bf C}\cap{\bf T}^d$. By virtue of {\bf (A-1-2)} and the analytic continuation, $g(z)$ vanishes on $M_{\lambda,reg}^{\bf C}$  (see e.g. Corollary 7 of \cite{KuVa}).\qed


\begin{lemma}\label{C4S2g(x)/p(z)analytic}
The meromorphic function $g(z)/p(z,\lambda)$ is analytic on ${\bf T}^d_{\bf C}$.\end{lemma}

Proof. By Lemma \ref{gz=0onMlambdaC}, $g(z)/p(z,\lambda)$ is analytic near $M_{\lambda,reg}^{\bf C}$. This can be proven by taking $\zeta_d = p(z,\lambda)$ as one of local coordinates, and expand $g(z)$ into a power series. 
Then the singularities of $g(z)/p(z,\lambda)$ are on $M_{\lambda,sng}^{\bf C}$.
However, since we have assumed $d \geq 2$ and (A-1-1), the singularities are removable (for the proof, see e.g. Corollary 7.3.2 of \cite{Kr82}). \qed

\medskip
We pass to the variables $w_j = e^{iz_j}$, $j = 1,\cdots, d$, and let ${\bf C}[w_1,\cdots,w_d]$ be the ring of polynomials of $w_1,\cdots,w_d$ with coefficients in $ {\bf C}$. The map
$$
{\bf T}^d_{\bf C} \ni z \to w \in {\bf C}^d\setminus
{\mathop\cup_{j=1}^d} A_j, 
\quad A_j = \{w \in {\bf C}^d\, ; \, w_j=0\}
$$
is biholomorphic, i.e. both of the mappings $z \to w$, $w \to z$ are holomorphic.
Let us note that $p(z,\lambda)$ has the form
\begin{equation}
p(z,\lambda) =\sum_{\alpha \in {\bf Z}^d, |\alpha|\leq N}
c_{\alpha}(\lambda)e^{i\alpha\cdot z}, \quad
\overline{c_{\alpha}(\lambda)} =c_{-\alpha}(\lambda).
\nonumber
\end{equation}
Letting $\gamma_j = {\rm max}_{|\alpha|\leq N}\alpha_j$, we factorize $p(z,\lambda)$ as
\begin{equation}
p(z,\lambda) = P(w,\lambda)\prod_{j=1}^dw_j^{-\gamma_j},
\label{C4S2pzlambda=Pwjpower}
\end{equation}
where $P(w,\lambda) \in {\bf C}[w_1,\cdots,w_d]$. Note that $\gamma_j \geq 0$, and this is the minimum choice of $\gamma_j$ for which the factorization (\ref{C4S2pzlambda=Pwjpower}) is possible. Similarly, we factorize $g(z)$ as
$$
g(z) = G(w)\prod_{j=1}^{d}w_j^{-\beta_j},
$$ 
where $\beta_j$ is a non-negative integer and $G(w) \in {\bf C}[w_1,\cdots.w_d]$.

 Let us recall some basic facts for \textit{analytic set} (see e.g. \cite{KuVa}, or Chapter 1 of \cite{Chi}).
A subset $E\subset {\bf C}^d $ is called an analytic set if $E$ is, in a neighborhood of each point $\in E$, the set of common zeros of a certain finite family of holomorphic functions.
An analytic set $E = \cap_{1 \leq j \leq N}f_j^{-1}(0)$, $f_j$ being analytic, splits into several parts : the set of regular points, which is an analytic submanifold with complex dimension $p =: \mathrm{dim}_{{\bf C}} E \leq d-1$, and the singular locus, which is a union of the set of singular points and submanifolds with complex dimension $ <p $.
Note that the set of regular points of $E$ is dense in $E$, and the singular locus of $E$ is nowhere dense in $E$ (see Lemma 6 of \cite{KuVa}).


\begin{lemma}\label{C4S2AjZlambdadim}
Let $Z_{\lambda} = \{w \in {\bf C}^d\, ; \, P(w,\lambda)=0\}$.
For $1 \leq j \leq d$,  $\dim_{\bf C}A_j\cap Z_{\lambda} \leq d-2$.
\end{lemma}

Proof. Assume $j=1$, and let $w'=(w_2,\cdots,w_d)$. We rewrite $P(w,\lambda)$ as 
$$
P(w,\lambda) = P_0(w',\lambda) + P_1(w',\lambda)w_1 + \cdots + P_{m}(w',\lambda)w_1^m,
$$
where $P_{\ell}(w',\lambda) \in {\bf C}[w_2,\cdots,w_d]$. 
If $\dim_{\bf C} A_1\cap Z_{\lambda} = d-1$, $P_0(w',\lambda)=0$ on an open set in $A_1$, hence it vanishes identically. Therefore, $P(w,\lambda) = w_1Q(w,\lambda)$ with $Q(w,\lambda) \in {\bf C}[w_1,\cdots, w_d]$. This implies
$$
p(z,\lambda) = w_1Q(w,\lambda)\prod_{j=1}^{d}w_j^{-\gamma_j} 
= Q(w,\lambda)\Big(w_1^{-(\gamma_1-1)}\prod_{j=2}^dw_j^{-\gamma_j}\Big).
$$
This contradicts the minimum choice of $\gamma_1$ in (\ref{C4S2pzlambda=Pwjpower}). \qed

\medskip
We now observe
\begin{equation}
\frac{g(z)}{p(z,\lambda)} = \frac{G(w)}{P(w,\lambda)}\prod_{j=1}^dw_j^{\gamma_j-\beta_j}.
\label{C4S2variableschange}
\end{equation}
Since $g(z)/p(z,\lambda)$ is analytic, $G(w)/P(w,\lambda)$ is analytic except possibly on hyperplanes, $A_j$, $j=1,\cdots,d$. In view of Lemma \ref{C4S2AjZlambdadim}, we then see that $G(w)/P(w,\lambda)$ is analytic except on some set of complex dimension at most $d-2$, which are removable singularities. For the proof, see e.g. Corollary 7.3.2 of \cite{Kr82}. We have, therefore,

\begin{itemize}
\item 
$G(w)/P(w,\lambda)$ is an entire function.
\end{itemize}
In particular,
\begin{itemize}
\item 
$G(w)=0$ on the set $\{w\in {\bf C}^d\, ; \, P(w,\lambda)=0\}$.
\end{itemize}

We factorize $P(w,\lambda)$ so that
$$
P(w,\lambda)= P^{(1)}(w,\lambda) \cdots P^{(N)}(w,\lambda),
$$
where each $P^{(j)} (w,\lambda) $ is an irreducible polynomial.
We prove inductively that 
$$
G(w)/\big(P^{(1)}(w,\lambda)\cdots P^{(n)}(w,\lambda)\big)\ {\rm is}\ {\rm a}\ {\rm polynomial}\ {\rm for}\  1 \leq n \leq N.
$$
Note that, since we know already that $G(w)/P(w,\lambda)$ is entire, 
\begin{itemize}
\item $G(w)/P^{(1)}(w,\lambda)\cdots P^{(n)}(w,\lambda)$ is also entire, 
\item $G(w) = 0$ on the zeros of $P^{(1)}(w,\lambda)\cdots P^{(n)}(w,\lambda)$.
\end{itemize}

\medskip
We make us of the Hilbert Nullstellensatz (see e.g \cite{Shafa}).


\begin{lemma}\label{C4S2HilbertNullestellensatz}
Suppose $f, g \in {\bf C}[w_1,\cdots,w_d]$ and $f$ is irreducible. If $g = 0$ on all zeros of $f$, there exists $h \in {\bf C}[w_1,\cdots,w_d]$ such that $g = fh$.
\end{lemma}

\medskip

Consider the case $n=1$. Since  $G(w) = 0$ on the zeros of $P^{(1)}(w,\lambda)$,
 Lemma \ref{C4S2HilbertNullestellensatz} implies that  $G(w) / P^{(1)} (w,\lambda)$ is a polynomial.
 
Assuming the case $n \leq \ell -1$, we consider the case $n=\ell$.
By the induction hypothesis, there exists a polynomial $P_{\ell-1}(w,\lambda)$ such that 
\begin{equation*} 
\frac{G(w)}{P^{(1)} (w,\lambda) \cdots P^{(\ell-1)} (w,\lambda) }=P_{\ell-1}(w, \lambda) . 
\end{equation*}
Then we have $G(w)/ (P^{(1) } (w,\lambda) \cdots P^{(\ell)} (w,\lambda) ) =P_{\ell-1}(w,\lambda)/P^{(\ell)}(w,\lambda)$. 
This is entire. Therefore, $P_{\ell-1}(w,\lambda) = 0$ on the zeros of $P^{(\ell)}(w,\lambda)$.
By Lemma \ref{C4S2HilbertNullestellensatz},
 there exists a polynomial $Q^{(\ell)}(w,\lambda) $ such that 
\begin{equation*}
 \frac{P_{\ell-1}(w,\lambda)}{P^{(\ell)} (w,\lambda)}=Q^{(\ell)}(w,\lambda). 
\end{equation*}
Therefore, $G(w)/\big(P^{(1)} (w,\lambda)\cdots P^{(n)}(w,\lambda)\big)$ is a polynomial for $1 \leq n \leq N$. 
Taking $n = N$, we have that $G(w)/P(w,\lambda)$ is a polynomial of $w$, hence $g(z)/p(z,\lambda)$ is a polynomial of $e^{iz_j}$ by (\ref{C4S2variableschange}).
This implies that $u(z)$ is a trigonometric polynomial. We have thus completed the proof of Theorem \ref{C4S2RellichTh}.
\qed


\subsection{Thresholds}\label{SubsectionThresholds}
Let $\lambda_1(x) \leq \lambda_2(x) \leq \cdots \leq \lambda_s(x)$ be the eigenvalues of $H_0(x)$, and 
\begin{equation}
M_{\lambda, j} = \{x \in {\bf T}^d\, ; \, \lambda_j(x) = \lambda\}.
\label{C4S1MlambdajDefine}
\end{equation}
Then we have
\begin{equation}
p(x,\lambda) = \prod_{j=1}^s(\lambda_j(x) - \lambda),
\label{C4S1pxlambda=prodpj}
\end{equation}
\begin{equation}
M_{\lambda} = {\mathop\cup_{j=1}^s} M_{\lambda,j}.
\label{C4S1Mlambda=cupMlambdaj}
\end{equation}
To study the spectral properties of $H_0$, we need another series of assumptions : 

\bigskip
\noindent
 {\it There is a finite set $\mathcal T_0 \subset \sigma(H_0)$ such that }

\medskip
\noindent
(A-2) $\ 
M_{\lambda, i} \cap M_{\lambda,j} = \emptyset, \quad {\rm if} \quad i \neq j, \quad \lambda \in \sigma(H_0)\setminus\mathcal T_0.$ 

\medskip
\noindent
(A-3) $\ 
\nabla_xp(x,\lambda) \neq 0, \quad {\rm on} \quad M_{\lambda}, \quad \lambda \in \sigma(H_0)\setminus\mathcal T_0.$

\bigskip
The assumption (A-2) implies that the eigenvalues are distinct in a neighborhood of $M_{\lambda}$ if $\lambda \in \sigma(H_0)\setminus{\mathcal T}_0$, moreover $H_0(x)$ is smoothly diagonalizable. 
By (A-2),  (\ref{C4S1pxlambda=prodpj}) and (\ref{C4S1Mlambda=cupMlambdaj}), the assumption (A-3) is equivalent to 
$$
\nabla_x\lambda_j(x) \neq 0, \quad {\rm on} \quad M_{\lambda,j}, \quad \lambda \in \sigma(H_0)\setminus\mathcal T_0.
$$
We use (A-2) and (A-3) in the proof of limiting absorption principle for the resolvent (see \S 6).

Let us examine the assumptions (A-1), (A-2), (A-3) for the examples in \S 3.
Let $a_d(x)$, $b_d(x)$ be defined by (\ref{C1S4Definead(z)}) and (\ref{C1S4Definebd(z)}).
In view of Lemma \ref{LemmaSquareLattice}, for $H_0(x) = a_d(x)$, we can take
$$
\mathcal T_0 = SV(a_d),\qquad 
\mathcal T_1 = \{-d, d\}.
$$
By Lemma \ref{LemmaDiamond}, for $H_0(x) = b_d(x)$, we can take
$$
\mathcal T_0 = SV(b_d), \qquad
\mathcal T_1 = \{- (d+1)/2, d(d+1)/2\}.
$$
Therefore, we have

\begin{itemize}
\item 
for the $d$-dim. square lattice
$$
\mathcal T_0 = \{n/d\, ; \, -d \leq n \leq d\}, \quad 
\mathcal T_1 = \{-1, 1\},
$$
\item 
for the triangular lattice
$$
\mathcal T_0 = \{-1, 1/3, 1/2\}, \qquad
\mathcal T_1 = \{-1, 1/2\},
$$
\item
for the hexagonal lattice
$$
\mathcal T_0 = \{-1, - 1/3, 0, 1/3, 1\}, \qquad 
\mathcal T_1 = \{-1, 0, 1\},
$$
\item
for the Kagome lattice
$$
\mathcal T_0 = \{-1, - 1/4, -1/2, 0, 1/2\}, \qquad 
\mathcal T_1 = \{-1, -1/4, 1/2\},
$$
\item
for the $d$-dim. diamond lattice,
$$
\mathcal T_0 = \left\{
\begin{split}
&\{\pm (\ell + 1)/(d+1)\, ; \, \ell = d, d-2, \cdots, -d\}\cup\{0\}, \quad {\rm if} \quad d = {\rm even}, \\
& \{\pm (\ell +1)/(d+1)\, ; \, \ell = d, d-2,\cdots, -d\}, \quad {\rm if} \quad d = {\rm odd},
\end{split}
\right.
$$
$$
\mathcal T_1=
\{-1, 0, 1\},
$$
\item
for the subdivision of $d$-dim. square lattice ${\bf Z}^d$,
$$
\mathcal T_0 = \{0, \pm n/d; \, n=1, 2, \cdots, d\}, \qquad
\mathcal T_1 = \{0,\pm 1\}.
$$
\end{itemize}

Attention must be payed to the ladder of the square lattice ${\bf Z}^d$ in ${\bf R}^{d+1}$ and the graphite. For both the cases, some connected component of $M_{\lambda}^{\bf C}$ has no intersection with ${\bf T}^d$. We have
\begin{itemize}
\item
for the ladder of the square lattice ${\bf Z}^d$ in ${\bf R}^{d+1}$,
$$
\mathcal T_0 = \Big\{-1, \frac{-2d+1}{2d+1}, \frac{-2d+3}{2d+1},\cdots,\frac{2d-1}{2d+1}, 1\Big\}, \quad
\mathcal T_1 = \Big\{\frac{2d-1}{2d+1} \leq |\lambda| \leq 1\Big\},
$$
\item 
for the graphite
$$
\mathcal T_0 = \{0, \pm 1/4, \pm 1/2, \pm 1\}, \quad
\mathcal T_1 = \{1/2 \leq |\lambda| \leq 1\}.
$$
\end{itemize} 


\subsection{Unique continuation property}
 As in the case of  differential operators, the unique continuation property is a delicate issue for Laplacians on periodic lattices. In this paper, we do not pursue the general condition for it, but assume it in the following form.

\medskip
\noindent
(A-4) \  Suppose $\widehat u$ satisfies $(\widehat H_0 - \lambda)\widehat u = 0$ on $\mathcal V_0$ for some constant $\lambda \in  {\bf C}$. If 
 there exists $R_0>0$ such that $\widehat u = 0$ for $|n| > R_0$, 
 then $\widehat u =0$ on $\mathcal V_0$.

\medskip
  Theorem \ref{C4S2RellichTh} then implies the following theorem.


\begin{theorem} \label{UniqueConti}
Assume  (A-1) and (A-4).
Suppose $\widehat u(n)$ satisfies 
$(\widehat H_0 - \lambda)\widehat u = 0$ on $\mathcal V_0$ for some $\lambda \in \sigma(\widehat H_0)\setminus \mathcal T_1$. If
$$
\lim_{R\to\infty}\frac{1}{R}\sum_{|n| < R}|\widehat u(n)|^2 = 0, 
 $$
then,  $\widehat u = 0$ identically on $\mathcal V_0$. 
In particulr,  $\widehat H_0$ has no eigenvalue in $\sigma(\widehat H_0)\setminus \mathcal T_1$.
\end{theorem}

\medskip
The validity of the unique continuation property (A-4) depends largely on the geometry of the lattice.  We give here an algebraic condition to guarantee (A-4).


\begin{lemma}\label{p(x,lambda)=f(ad)}
Suppose there exists a polynomial $f(z,\lambda)$ such that $p(x,\lambda)$ is written as 
$p(x,\lambda) = f(a_d(x),\lambda)$ or $p(x,\lambda) = f(b_d(x), \lambda)$. If $p(x,\lambda)$ 
is a non-zero polynomial, $\widehat H_0- \lambda$ has the unique continuation property (A-4).
\end{lemma}

Proof.  Passing to the torus, and multilplying the cofactor matrix of $H_0(x)-\lambda$, we  obtain the equation
$p(x,\lambda)u = 0$. Returning to the lattice, $\widehat u$ satisfies
\begin{equation}
\widehat P(\lambda)\widehat u =  0,
\label{C4S3Pv=Ev}
 \end{equation} 
 where $\widehat P(\lambda)$ is defined by $p(x,\lambda)$ with $e^{-ix_j}$, $e^{ix_j}$ replaced by shift operators $\widehat S_j$, $\widehat S_j^{\ast}$, respectively.

 We let 
 $$
\widehat S = \sum_{j=1}^d\left(\widehat S_j + \widehat S_j^{\ast}\right), \quad
\widehat T = \sum_{1\leq j<k\leq d}\left(\widehat S_j\widehat   S_k^{\ast} + \widehat S_k\widehat S_j^{\ast}\right).
$$
Then $2 a_d(x)$ and $2b_d(x)$ correspond to 
$\widehat S$ and $- d-1 + \widehat S + \widehat T$, respectively.
By the assumption of Lemma \ref{p(x,lambda)=f(ad)},  the equation (\ref{C4S3Pv=Ev}) is rewritten as
either
\begin{equation}
\sum_{p =0}^Nc_{p}\widehat S^{p}\,\widehat u =0,
\label{C4S3Fourierequation1}
\end{equation}
or 
\begin{equation}
\sum_{p =0}^Nc_{p}\big(\widehat S + \widehat T\big)^{p}\widehat u =0,
\label{C4S3Fourierequation2}
\end{equation}  
where $c_{p}$ is a constant and $c_N=1$.

\medskip
\noindent
{\it The case (\ref{C4S3Fourierequation1})}.
For $n \in {\bf Z}^d$, we put
\begin{equation}
\begin{split}
D_n & = \{n - {\bf e}_1, \widehat S_j(n - {\bf e}_1), \widehat S_j^{\ast}(n - {\bf e}_1)\, ; \, j= 1,\cdots,d\}\\
& = \{n-{\bf e}_1 \pm {\bf e}_j\, ; \, j=1,\cdots, d\},
\end{split}
\label{C4S3DefineDnCase1}
\end{equation}
\begin{equation}
\Omega_n = \{\ell \in {\bf Z}^d\, ; \, 
\sum_{j=2}^d|n_j-\ell_j| \leq n_1-\ell_1\}.
\label{C4S3OmegapCase1}
\end{equation}
Geometrically, $\Omega_n$ is a cone with vertex $n$, and related with $D_n$ as follows. We define 
$$
D_k \prec D_{\ell} \Longleftrightarrow 
\ell \in D_k \setminus \{k\}.
$$
Then starting from $D_n$, one can construct a chain of $D_k$'s satisfying
$$
D_n \prec D_k \prec D_{k'} \prec \cdots.
$$
$\Omega_n$ is the union of such chains.

\begin{figure}[hbtp] 
\centering
\includegraphics[width=6cm, bb=0 0 410 406]{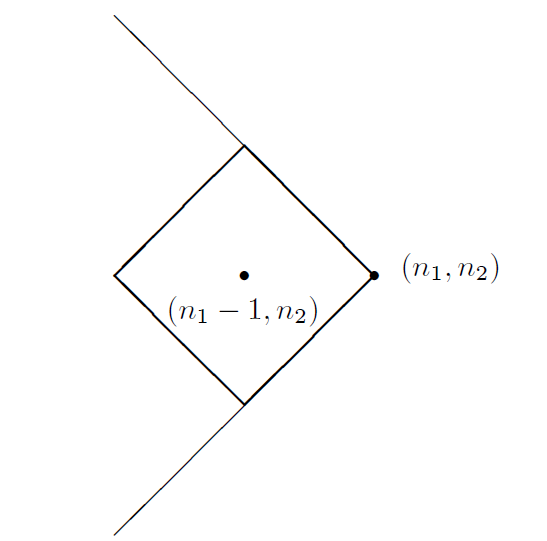}
\caption{$D_n$ and $\Omega_n$}
\label{DnOnforad}
\end{figure}

Note that, although $\widehat u(n)$ is defined on the multiple lattice $\mathcal V_0$, each component $\widehat u_i(n)$ is a function on a single lattice ${\bf Z}^d$.
Evaluating the equation (\ref{C4S3Fourierequation1}) at $k \in {\bf Z}^d$, we see that $\widehat u_i(k + N{\bf e}_1)$ is written by a linear combination of $\widehat u_i(\ell)$ for $\ell \in \Omega_{k + N{\bf e}_1}$, and $\widehat u_i(\ell)$ is written as a linear combination of $\widehat u_i(m)$, where $m \in \Omega_{\ell}$. This procedure can be repeated as long as possible. Now suppose $\widehat u_i(k) = 0$ near infinity.
By the above procedure, we then see that $\widehat u_i(n) =0$. Hence $\widehat u=0$ identically. 

\medskip
\noindent
{\it The case  (\ref{C4S3Fourierequation2})}. We define
$$
D_n = \big\{\widehat S_j(n-{\bf e}_1), \widehat S_j^{\ast}(n-{\bf e}_1),  
\widehat S_i\widehat S_j^{\ast}(n-{\bf e}_1), \widehat S_i^{\ast}\widehat S_j(n-{\bf e}_1)\, ; \,
1\leq i, j \leq d\big\},
$$
$$
\Omega_n = \{\ell \in {\bf Z}^d\, ; \, 
\ell_1 \leq n_1, \ell_i + \ell_j \leq n_i + n_j, \ 
1 \leq i,j \leq d\big\}.
$$
Then, the same argument works also for this case. \qed

\begin{figure}[hbtp] 
\centering
\includegraphics[width=6cm, bb=0 0 427 410]{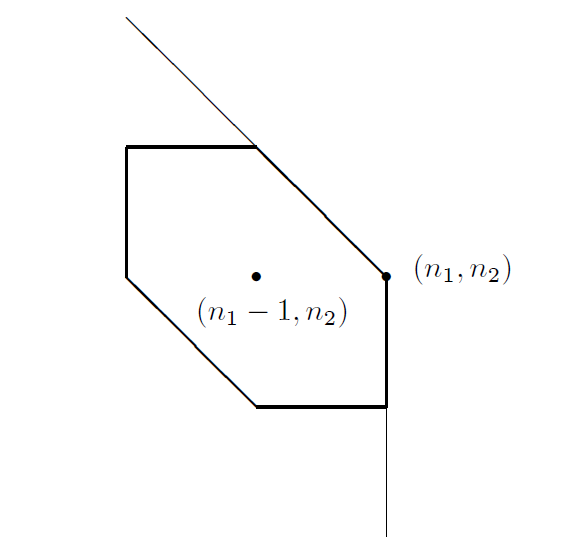}
\caption{$D_n$ and $\Omega_n$}
\label{DnOnforad}
\end{figure} 

\bigskip
As has been computed in \S 3, all of our examples satisfy the assumption of Lemma \ref{p(x,lambda)=f(ad)}, hence have the 
unique continuation property. Note that $p(x,1/2)=0$ for the Kagome lattice by (\ref{pxlambdaKagome}), and $p(x,0)=0$ for the subdivision by (\ref{Subdivisionpxlambda}). Later, we show that 
$\sigma_p(\widehat H_0) = \{1/2\}$ for the Kagome lattice, and $\sigma_p(\widehat H_0) = \{0\}$ for the
subdivision. See (\ref{sigmapH0Kagome}) and  (\ref{sigmapSubdivision}).

\bigskip

We can also add perturbations by scalar potentials. Here, by a scalar potential on a graph $\Gamma = \{\mathcal V, \mathcal E\}$, we mean an operator $\widehat V$ such that
$$
(\widehat V\widehat f)(v) = \widehat V(v)\widehat f(v), \quad \forall v \in \mathcal V, 
$$
where $\widehat V(v) \in {\bf C}$. First let us consider the diamond lattice.


\begin{lemma}
\label{Lemmaddimdia} Let $\widehat H_0$ be the Laplacian on the $d$-dimensional diamond lattice, where $d\geq 2$, and $\widehat V$ a compactly supported scalar potential. 
Then $\widehat H = \widehat H_0 + \widehat V$ has the unique continuation property (A-4). 
\end{lemma}

Proof. Suppose $(\widehat H_0 -\lambda)\widehat u = - \widehat V\widehat u$. Taking account of (\ref{S3ddimdiapxlambda}) 
and 
mutiplying this equation by the matrix 
\begin{equation}
\widehat C_0 = \left(
\begin{array}{cc}
0 & 1 + \widehat S_1^{\ast} + \cdots + \widehat S_d^{\ast} \\
1 + \widehat S_1 + \cdots + \widehat S_d & 0
\end{array}
\right),
\nonumber
\end{equation}
we have
\begin{equation}
\begin{split}
& \sum_{i}\big(\widehat u(n + {\bf e}_i) + \widehat u(n - {\bf e}_i)\big)+ \sum_{ i<j }
\big(\widehat u(n + {\bf e}_i - {\bf e}_j) + \widehat u(n + {\bf e}_j - {\bf e}_i)\big) \\ 
& = c_1\widehat u(n) + c_2\left(\widehat C_0\widehat V\widehat u\right)(n), 
\end{split}
\label{S6Ddimdiaeq}
\end{equation}
where $c_i$ is a constant. The right-hand side is rewritten as
\begin{equation}
c_1\left(
\begin{split}
\widehat u_1(n) \\
 \widehat u_2(n)
\end{split}
\right) 
+ 
c_2
\left(
\begin{split}
\widehat u_2(n) + \sum_{j=1}^d b(n+{\bf e}_j)\widehat u_2(n + {\bf e}_j) \\
\widehat u_1(n) + \sum_{j=1}^d a(n-{\bf e}_j)\widehat u_1(n - {\bf e}_j)
\end{split}
\right),
\label{S6Ddimdiaeq1}
\end{equation}
where $a(n), b(n) \in {\bf C}$. To prove the lemma, we have only to show that if $\widehat u(n) = 0$ for all $n= (n_1,n')$ such that $n_1 < k_1$, then $\widehat u(k_1,n')=0$ for all $n' \in {\bf Z}^{d-1}$. Let $n_1= k_1-1$ in (\ref{S6Ddimdiaeq}), and take the lower component. Since the right-hand side vanishes, we have
$$
\widehat u_2(k_1,n') + \sum_{j=2}^d\widehat u_2(k_1,n'-{\bf e}_j) = 0, \quad \forall n'.
$$
In ${\bf Z}^{d-1}$, consider a cone with vertex $n'$: 
$$
C(n') = \{m' = (m'_2,\cdots,m'_{d})\, ; \, m_i' \leq n_i', \ 2 \leq i \leq d\}.
$$
Since $\widehat u_2(k_1,m') = 0$ near infinity of $C(n')$, one can show inductively that $\widehat u_2(k_1,n')=0$.  Taking the upper component, we then have
$$
\widehat u_1(k_1,n') + \sum_{j=2}^d\widehat u_1(k_1,n'-{\bf e}_j) = 0, \quad \forall n'.
$$
Arguing as above, we have $\widehat u_1(k_1,n')=0$. 
This proves the lemma. \qed

\bigskip
We consider the other examples.


\begin{theorem} \label{UniqueConti2} 
Let $\widehat H_0$ be the Laplacian of one of the
following lattices: square lattice, 
  triangular lattice, $d$-dimensional  diamond lattice ($d \geq 2$), 
  ladder of $d$-dimensional square lattice, 
  graphite in ${\bf R}^3$.
Let $\widehat V$ be a complex-valued compactly supported scalar potential. Then $\widehat H = \widehat H_0 + \widehat V$ has the unique continuation property (A-4).
In particulr,  $\widehat H$ has no eigenvalue in $\sigma(\widehat H_0)\setminus \mathcal T_1$.
\end{theorem}

\begin{figure}
\centering
\includegraphics[width=12cm, bb=0 0 595 204]{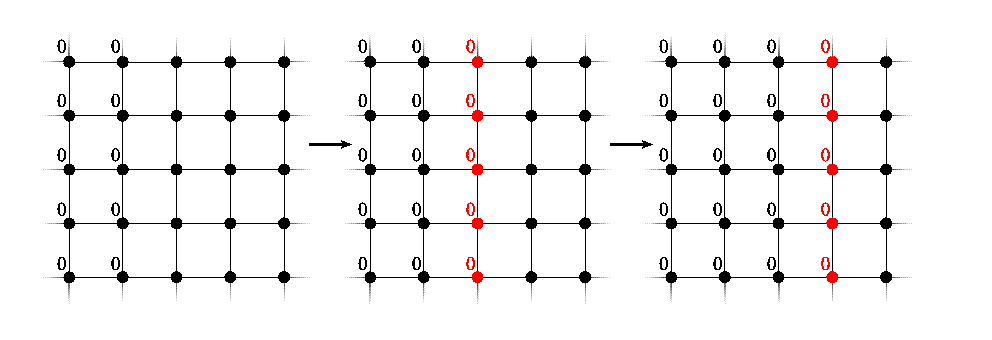}
\caption{The unique continuation on the square lattice.}
\label{UniqContSqrLattice}
\end{figure}

\begin{figure}
\centering
\includegraphics[width=13cm, bb=0 0 595 205]{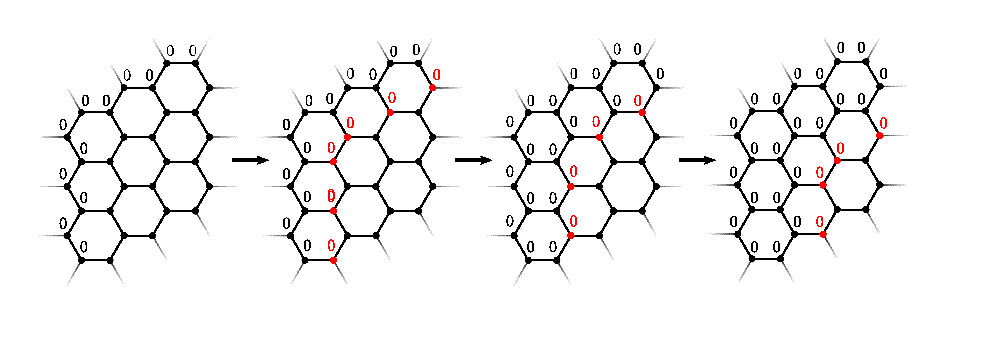}
\caption{The unique continuation on the hexagonal lattice.}
\label{UniqContHexLattice}
\end{figure}

\medskip
This theorem is proven by observing the figure of the graph, and the idea of the proof is similar to the case of square lattice given in \cite{IsMo1}.
(See the Figure~\ref{UniqContSqrLattice}.)
We enlarge the region on which $\widehat u=0$ step by step by using the equation $(\widehat H - \lambda)\widehat u=0$. 
Let us illustrate it for the hexagonal lattice, although this is the case of $d=2$ in Lemma \ref{Lemmaddimdia}. From the equation $- \widehat \Delta_{\Gamma} \widehat u + (\widehat V - \lambda)\widehat u = 0$, one obtains
\begin{equation}
  \widehat u_2(n_1,n_2-1) = - \widehat u_2(n_1-1,n_2) - \widehat u_2(n_1,n_2)  + 3\left(\widehat V(n_1,n_2) - \lambda\right)\widehat u_1(n_1,n_2),
\nonumber
\end{equation}
\begin{equation}
\widehat u_1(n_1+1,n_2) = - \widehat u_1(n_1,n_2+1) - \widehat u_1(n_1,n_2)  + 3\left(\widehat V(n_1,n_2) - \lambda\right)\widehat u_2(n_1,n_2).
\nonumber
\end{equation}
The left-hand side vanishes, if so does each term of the right-hand side, and this occurs by the assumption 
that $\widehat u(n)=0$ near infinity and the induction procedure. (See the Figure~\ref{UniqContHexLattice}.)

\begin{figure}
\centering
\includegraphics[width=13cm, bb=0 0 595 163]{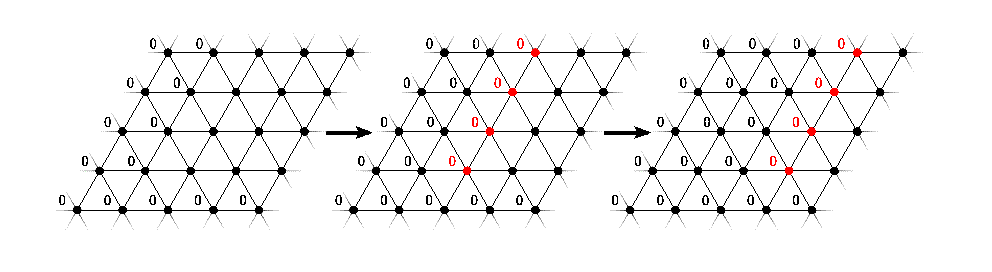}
\caption{The unique continuation on the triangular lattice.}
\label{UniqContTriangLattice}
\end{figure}

\begin{figure}
\centering
\includegraphics[width=13cm, bb=0 0 595 214]{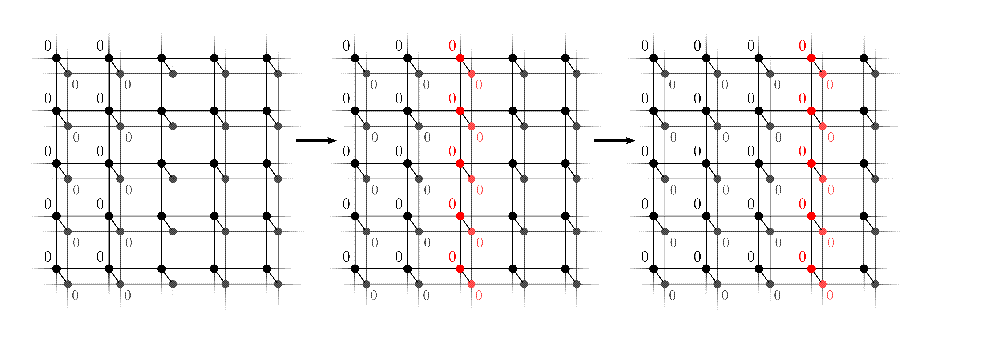}
\caption{The unique continuation on the square ladder.}
\label{UniqContSquareLadder}
\end{figure}

In the case of the triangular lattice and the square ladder, the unique continuation procedure is illustrated in  
 Figures~\ref{UniqContTriangLattice} and ~\ref{UniqContSquareLadder}, respectively. 
We omit the proof for the graphite, which can be easily imagined by comparing  Figures \ref{UniqContHexLattice} and \ref{C1Graphite}.

\bigskip

\begin{figure}
\centering
\includegraphics[width=12cm, bb=0 0 595 370]{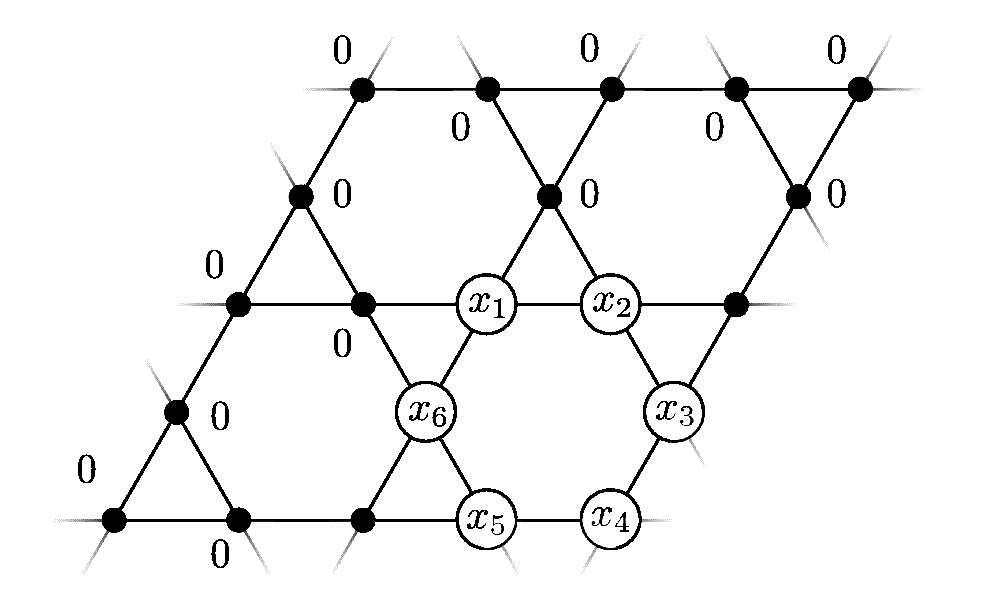}
\caption{The unique continuation fails on the Kagome lattice.}
\label{UniqContKagomeLattice}
\end{figure}

The Figure~\ref{UniqContKagomeLattice} shows the reason why this procedure does not work for the case of the Kagome lattice. In fact, for the potential
\[
\widehat{V} (n) =
\begin{cases}
  v,&\ \text{ if } \ \  n = x_j, \ \  j = 1, 2, \cdots, 6,\\
  0,&\ \text{ otherwise},
\end{cases}
\]
$v$ being an arbitrary constant,
we have an eigenvalue $\lambda_v = v + 1/2$ with the compactly supported eigenvector
\[
\widehat{u} (n) =
\begin{cases}
  (-1)^j,&\ \text{ if } \ \ n = x_j,\ \  j = 1, 2, \cdots, 6\\
  0,&\ \text{otherwise}.
\end{cases}
\]
 Leting $v=0$, we get that
\begin{equation}
\sigma_p(\widehat H_0) = \{1/2\}
\label{sigmapH0Kagome}
\end{equation}
for the Kagome lattice. 
Note that if $-3/2 < v < 0$, $\lambda_v = v + 1/2$ is an embedded eigenvalue for the non-zero potential $\widehat V$. See \cite{S99} for the case of $\widehat{V} \equiv 0$.

\begin{figure}
\centering
\includegraphics[width=8cm, bb=0 0 595 568]{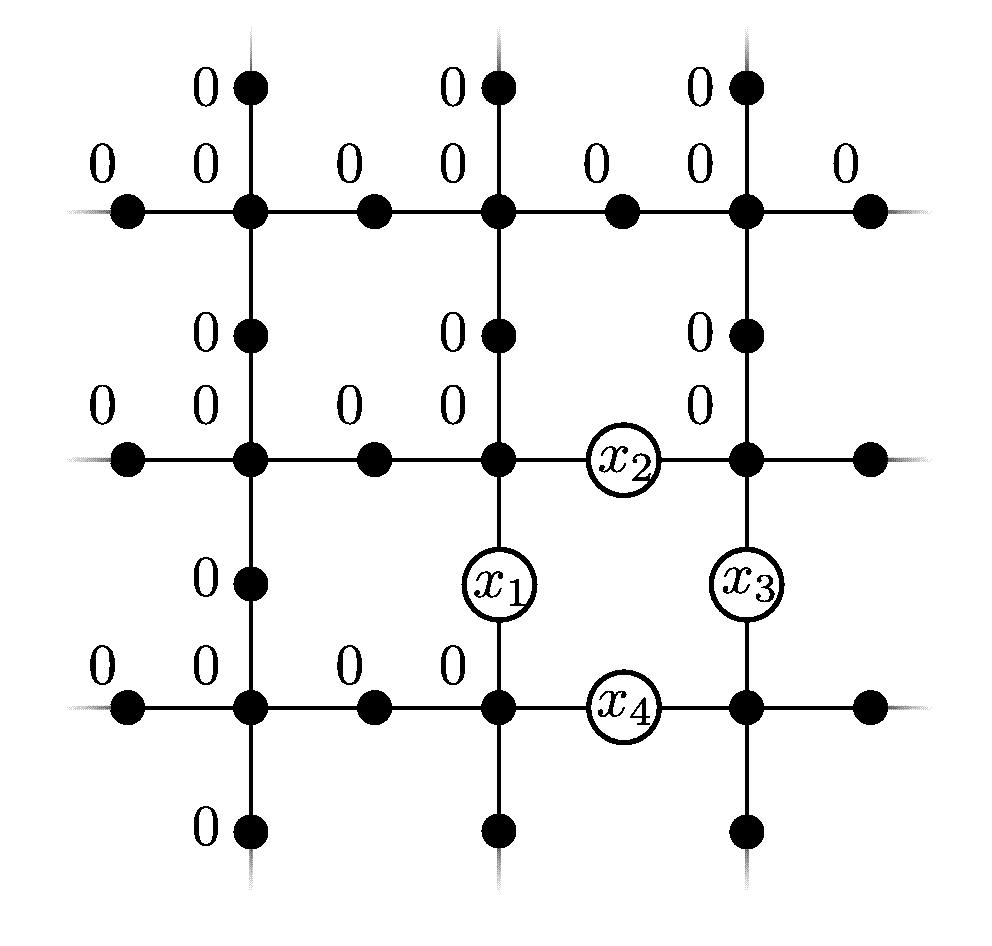}
\caption{The unique continuation fails on the subdivision of $2$-dimensional square lattice.}
\label{UniqContSubdivision}
\end{figure}

The unique continuation also fails for the subdivision, which is illustrated in the Figure~\ref{UniqContSubdivision} for the  $2$-dimensional square lattice case.
 In fact,  the potential
\[
\widehat{V} (n) =
\begin{cases}
  v,& \ \text{ if } \ \ n = x_j, \ \ j = 1, 2, 3, 4, \\
  0,& \ \text{ otherwise},
\end{cases}
\]
$v$ being an arbitrary constant, has an eigenvalue $\lambda_v = v$ with the compactly supported eigenvector
\[
\widehat{u} (n) =
\begin{cases}
  (-1)^j,& \ \text{ if }\ \  n = x_j, \ \ j = 1, 2, 3, 4, \\
  0, & \ \text{ otherwise}.
\end{cases}
\]
Letting $v=0$, we have
\begin{equation}
\sigma_p(\widehat H_0) = \{0\}
\label{sigmapSubdivision}
\end{equation}
for the subdivision.

\subsection{A counter example}
The assumption {\bf (A-1-2)} is characteristic in its topological feature. For the case of differential operators, H{\"o}rmander \cite{Hor73} imposed a similar assumption and also proved that it is necessary for the Rellich type theorem. In the discrete case, we can construct examples showing  the necessity of excluding the set $\mathcal T_1$.

Consider the Hamiltonian $\widehat H_0$ for the ladder of the square lattice in ${\bf R}^{d+1}$ (Subsection \ref{SubsectionLadder}). We add a scalar potential $\widehat V$ to $\widehat H_0$, where
$\widehat V(0) = c$, $\widehat V(n) =0$ for $n \neq 0$, $c$ being a real constant yet to be determined. Then the eigenvalue problem is rewritten as
$\big(\widehat H_0(n) - \lambda\big)\widehat u(n) = - c\widehat u(0)\delta_{n0}$, where $\delta_{n0} = 1 $ for $n=0$, $\delta_{n0} = 0$ for $n \neq 0$, i.e.
$$
\big(-\frac{2}{2d+1}\sum_{j=1}^d\cos x_j - \lambda\big)u_1(x) - \frac{1}{2d+1}
u_2(x) =  -\frac{c}{2\pi}\widehat u_1(0),
$$
$$
\big(-\frac{2}{2d+1}\sum_{j=1}^d\cos x_j - \lambda\big)u_2(x) - \frac{1}{2d+1}
u_1(x) =  -\frac{c}{2\pi}\widehat u_2(0).
$$
We seek the solution in the form $u_1(x)=\pm u_2(x)=:v_{\pm}(x)$, and obtain
$$
v_{\pm}(x) =  c \,\frac{\widehat v_{\pm}(0)}{2\pi}\frac{1}{\dfrac{1}{2d+1}(\pm1 + 2\sum_{j=1}^d\cos x_j) + \lambda}.
$$
If $\pm\lambda > (2d-1)/(2d+1)$, we take $v_{\pm}(x) = \big(\frac{1}{2 d+1}(\pm 1+2\sum_{j=1}^d\cos x_j)+\lambda\big)^{-1}$, and $c = 2\pi/\widehat v_{\pm}(0)$ to get the desired potential and eigenvector.
Note that the associated eigenvector $\widehat v_{\pm}(n)$ is not compactly supported.

\medskip

In the case of graphite, we consider the eigenvalue problem $\big(\widehat H_0 + \widehat{V}\big)\widehat u = \lambda \widehat u$, where the scalar potential is $\widehat V(0)=(c_1, c_2, c_1, c_2)$, $\widehat V(n)=0$ for $n\neq0$.
We rewrite it as
$$
\left\{
-\frac{1}{4}
\left(
  \begin{array}{cc}
    0 & \overline{c(x)} \\
    c(x) & 0
  \end{array}
\right)-\lambda
\right\}
\left(
  \begin{array}{c}
    u_1(x) \\
    u_2(x)
  \end{array}
\right)
-\frac{1}{4}
\left(
  \begin{array}{c}
    u_3(x) \\
    u_4(x)
  \end{array}
\right)
=\frac{1}{2\pi}
\left(
  \begin{array}{c}
    c_1 {\widehat u}_1(0) \\
    c_2 {\widehat u}_2(0)
  \end{array}
\right),
$$
$$
\left\{
-\frac{1}{4}
\left(
  \begin{array}{cc}
    0 & \overline{c(x)} \\
    c(x) & 0
  \end{array}
\right)-\lambda
\right\}
\left(
  \begin{array}{c}
   u_3(x) \\
    u_4(x)
  \end{array}
\right)
-\frac{1}{4}
\left(
  \begin{array}{c}
    u_1(x) \\
    u_2(x)
  \end{array}
\right)
=\frac{1}{2\pi}
\left(
  \begin{array}{c}
    c_1 {\widehat u}_3(0) \\
    c_2 {\widehat u}_4(0)
  \end{array}
\right),
$$
where $c(x)$ is defined by (\ref{C1S4Graphite_c(x)}).
Putting $u_1(x) = \pm u_3(x) = v_{\pm,1}(x)$ and $u_2(x) = \pm u_4(x) = v_{\pm,2}(x)$, we obtain
$$
-\frac{1}{4}
\left(
  \begin{array}{cc}
    4 \lambda \pm 1 & \overline{c(x)} \\
    c(x) & 4 \lambda \pm 1
  \end{array}
\right)
\left(
  \begin{array}{c}
    v_{\pm,1}(x) \\
    v_{\pm,2}(x)
  \end{array}
\right)
=\frac{1}{2\pi}
\left(
  \begin{array}{c}
    c_1 v_{\pm,1}(0) \\
    c_2 v_{\pm,2}(0)
  \end{array}
\right),
$$
which is solved as
\begin{eqnarray*}
& &\left(
  \begin{array}{c}
    v_{\pm,1}(x) \\
    v_{\pm,2}(x)
  \end{array}
\right) \\
& =& \frac{1}{\left(\frac{3}{4}\right)^2\left\{\left(\frac{4}{3}\lambda \pm \frac{1}{3}\right)^2 - \frac{\alpha(x)}{9}\right\}}
\left(-\frac{1}{8 \pi}\right)
\left(
  \begin{array}{cc}
    4 \lambda \pm 1 & -\overline{c(x)} \\
    -c(x) & 4 \lambda \pm 1
  \end{array}
\right)
\left(
  \begin{array}{c}
    c_1 v_{\pm,1}(0) \\
    c_2 v_{\pm,2}(0)
  \end{array}
\right),
\end{eqnarray*}
where $\alpha(x)$ is defined by (\ref{C1HexagonalLattice_alpha(x)}). Note that the characteristic polynomial of the Hamiltonian on the hexagonal lattice (\ref{S3pxlambdahexa}) appears in the dominator with $\lambda$ replaced by $4\lambda/3\pm1/3$.
If $\pm \lambda > 1/2$, choosing $c_{j}=-8 \pi / v_{\pm,j}(0)$, $j = 1, 2$, 
we have the desired potential and eigenvectors which are not compactly supported.

\medskip
Due to the formula (\ref{hamiltonian_ladder}),
 the same argument as above works for the ladder $Lad(\Gamma)$ of $k$-regular periodic graph $\Gamma$. The characteristic polynomial of the associated Hamiltonian is given by (\ref{char_pol_ladder}) using that on $\Gamma$.

\medskip
Summing up the arguments in this subsection, we have


\begin{theorem}
  Let $\Gamma$ be a $k$-regular periodic graph, and $Lad(\Gamma)$ the ladder of $\Gamma$. Then we have \\
  (1)  For any $\lambda \in \big( -1, -\frac{k-1}{k+1} \big) \cap \sigma(\widehat H_0)$, 
  there exists $0 \neq \widehat v \in L^2({\bf T}^d) \subset {\mathcal B}^*_0$ 
  satisfying $(\widehat H_0 - \lambda) \widehat v = 0$ in $| n | > R_0$ for sufficiently large $R_0$. \\
  \noindent
  (2) For  any
  $\lambda \in \big( -1, -\frac{k-1}{k+1} \big) \cap \sigma(\widehat H_0)$, there exists a compactly supported potential $\widehat V$ such that $\lambda$ is an eigenvalue of $\widehat H_0 + \widehat V$. \\
  \noindent
(3) Let $\sigma_{max} := \max{\left( \sigma(\widehat H_0) \setminus \sigma_p(\widehat H_0) \right)}$. 
Then, the assertions (1) and (2) hold with $\big( -1, -\frac{k-1}{k+1} \big)$ replaced by  $\big( \frac{k \sigma_{max}-1}{k+1}, \sigma_{max} \big)$.
\end{theorem}


\begin{remark}
Let  $\lambda_1(x) \le \lambda_2(x) \le \cdots \le \lambda_s(x)$ be the eigenvalues of 
 $H_0(x)$
as in Subsection \ref{SubsectionThresholds}.
Then $\lambda_1(x)$ attains its minimum at $x = 0 \in {\bf T}^d$ and $\lambda_1(0)=-1$, which is a simple eigenvalue of $H_0(0)$ with a constant eigenfunction.
Moreover, the Hessian of $\lambda_1(x)$ at $x = 0 \in {\bf T}^d$ is positive definite.  Therefore, $\inf\sigma(\widehat H_0) = -1$ and $[-1,-1+\epsilon] \subset \sigma(\widehat H_0)$ for some $\epsilon > 0$.  See e.g. \cite{KoShiSu98}.
\end{remark}


\begin{remark}
  We cannot take $\sigma_{max}$ as $\max{\sigma(\widehat H_0)}$, since it could be an isolated eigenvalue with $\infty$-multiplicity and the interval in (3) might be empty. 
  Such an example of periodic graph is $L(L(\Gamma))$, i.e. the line graph of the line graph of a $k$-regular periodic graph $\Gamma$ with $k\ge3$.
  See \cite{S99}.
\end{remark}

Let us further remark here that \cite{Ship} obtains the eigenfunctions for embedded eigenvalues with
 unbounded support in a way very similar to above on periodic combinatorial graphs.
Note also the construction of eigenfunctions with such properties on periodic metric graphs given there.


\section{Spectral properties of  the free Hamiltonian}


\subsection{Resolvent estimates}
We now study $H_0$ in $\left(L^2({\bf T}^d)\right)^s$  under the assumptions (A-1), (A-2), (A-3). 
The eigenvalues of $H_0(x)$ are arranged so that
$\lambda_1(x) \leq \cdots \leq \lambda_s(x)$, 
and  assumed to be distinct on $M_{\lambda}$. Let $P_j(x)$ be the eigenprojection associated with $\lambda_j(x)$.
Take a compact interval $J \subset \sigma(H_0)\setminus\mathcal T_0$ and let $U^{\epsilon_0}$ be an $\epsilon_0$-neighborhood of $\cup_{\lambda\in J}M_{\lambda}$. Then $P_j(x) \in C^{\infty}(U^{\epsilon_0})$ for a sufficiently small $\epsilon_0 > 0$.
We fix $\lambda \in J$, and take $x_0 \in M_{\lambda}$. Let $\chi \in C_0^{\infty}(U^{\epsilon_0})$ such that $\chi(x)=1$ on a small neighborhood of $x_0$. Then, letting $u = R_0(z)f$, we have
$$
(\lambda_j(x) - z)\chi(x)P_ju(x) = \chi(x)P_jf(x).
$$
We can then apply Lemma \ref{C3S4h(x)u=flemma} to obtain the following theorem.


\begin{theorem}\label{C4S4TheoremR0lambda)}
(1) For $f \in \mathcal B$ and $\lambda \in \sigma(H_0)\setminus{\mathcal T}_0$, there  exists a weak $\ast$ limit, $\lim_{\epsilon \to 0}R_0(\lambda \pm i\epsilon)f = : R_0(\lambda\pm i0)f$, i.e.
$$
(R_0(\lambda\pm i\epsilon)f,g) \to (R_0(\lambda\pm i0)f,g), \quad \forall g \in \mathcal B.
$$
(2) Moreover,
$$
\|R_0(\lambda\pm i0)f\|_{\mathcal B^{\ast}} \leq C\|f\|_{\mathcal B},
$$
where the constant $C$ is independent of $\lambda$ when $\lambda$ varies over a compact 
set in $\sigma(H_0)\setminus\mathcal T_0$. \\
\noindent
(3)  $R_{0}(\lambda \pm i0)f$ is weakly $\ast$ continuous in the sense that
$$
\sigma(H_0)\setminus\mathcal T_0 \ni \lambda \to  (R_0(\lambda\pm i0)f,g)
$$
is continuous for $f, g \in \mathcal B$. \\
\noindent
(4) Letting $u_{\pm} = R_0(\lambda\pm i0)f$, we have
$$
WF^{\ast}(P_ju_{\pm}) \subset \{(x,\pm\omega_x)\, ; \, x \in M_{\lambda,j}\},
$$
where $\omega_x \in S^{d-1}\cap T_x(M_{\lambda,j})^{\perp}$, and $\omega_x\cdot
\nabla \lambda_j(x) < 0$,
$$
P_ju_{\pm} \mp \frac{1}{\lambda_j(x)- \lambda \mp i0}\otimes
\left((P_jf)\big|_{M_{\lambda,j}}\right) \in \mathcal B^{\ast}_0.
$$
\end{theorem}

Proof. We take $h(x) = \lambda_j(x)-\lambda$ in  Lemma \ref{C3S4h(x)u=flemma}. The independence of the constant $C$ in (2) can be seen by examining the proof. The local uniformity of the convergence with respect to $\lambda$ as $\epsilon \to 0$ implies the continuity in (3). The other assertions are direct consequences of  Lemma \ref{C3S4h(x)u=flemma}.  \qed

\medskip
A solution $u \in \mathcal B^{\ast}$ of the equation
\begin{equation}
\left(H_0(x) - \lambda\right)u = f \in \mathcal B
\label{C4S4H0x-lambdau=f}
\end{equation}
is outgoing (incoming) if it satisfies
\begin{equation}
\begin{split}
WF^{\ast}(P_ju) \subset \{(x, \omega_x)\, ; \, x \in M_{\lambda,j}\}, \\
\Big(WF^{\ast}(P_ju) \subset \{(x, -\omega_x)\, ; \, x \in M_{\lambda,j}\}\Big)
\end{split}
\label{C4S4OutInRadCond}
\end{equation}
with $\omega_x$ satisfying $\omega_x \in S^{d-1}\cap T_x(M_{\lambda,j})^{\perp}$, and $\omega_x\cdot
\nabla \lambda_j(x) < 0$.

The next lemma follows from  Lemma \ref{S4Lemma47}.


\begin{lemma}\label{C4S4Im(upm,f)=|f|}
Let $u \in \mathcal B^{\ast}$ be a solution to (\ref{C4S4H0x-lambdau=f}). Then $u$ is outgoing (incoming) if and only if
\begin{equation}
P_j(x)u(x) = \frac{P_j(x)f(x)}{\lambda_j(x)-\lambda - i0}, \quad 
\Big(P_j(x)u(x) = \frac{P_j(x)f(x)}{\lambda_j(x)-\lambda + i0}\Big).
\label{S6Lemma62Pju=Pjfoverlambda}
\end{equation}
For the outgoing (incoming) solution $u_+$ ($u_-$), we have
$$
{\rm Im}\,(u_{\pm},f) = \pm \pi \big\|f\big|_{M_{\lambda}}\big\|^2_{L^2(M_{\lambda})}.
$$
\end{lemma}


\begin{lemma}\label{C4S4RadcondinBast0}
Suppose $u \in \mathcal B^{\ast}$ satisfies $(H_0 - \lambda)u = f\in \mathcal B$ and one of the radiation conditions. If ${\rm Im}(u,f) = 0$, then $u \in \mathcal B^{\ast}_0$.
\end{lemma}

This lemma follows from Theorem \ref{C4S4TheoremR0lambda)} (4) and Lemma \ref{C4S4Im(upm,f)=|f|}. The following lemma is a direct consequence of of Theorem \ref{C4S2RellichTh} and Lemma \ref{C4S4RadcondinBast0}.


\begin{lemma}\label{RadCondUniqueFree}
Let $\lambda \in \sigma(H_0)\setminus\mathcal T_1$, and suppose $u \in \mathcal B^{\ast}$ satisfies $(H_0 - \lambda)u=0$ in ${\bf T}^d$. Then, $u=0$ if $u$ satisfies the outgoing or incoming radiation condition.
\end{lemma}


\subsection{Spectral representation}
We define a spectral representation of $H_0$, which is essentially a diagonalization of $H_0$. Let 
$$
{\bf h}_{\lambda,j} = L^2\big(M_{\lambda,j},{\bf C}, dM_{\lambda,j}/|\nabla\lambda_j(x)|\big)
$$ 
be the Hilbert space of $\bf C$-valued functions on $M_{\lambda,j}$ equipped with the inner product
$$
(\phi,\psi) = \int_{M_{\lambda,j}}\phi(x)\overline{\psi(x)}\frac{dM_{\lambda,j}}{|\nabla\lambda_j(x)|}.
$$
We let 
$$
\widetilde{\bf H}_j = \left\{P_j(x)f(x)\, ; \, f(x) \in L^2({\bf T}^d,{\bf C}^s,dx)\right\}.
$$
By using an eigenvector $a_j(x) \in {\bf C}^s$  satisfying $H_0(x)a_j(x) = \lambda_j(x)a_j(x)$, $ |a_j(x)|=1$, we can rewrite
$$
P_j(x)f(x) = \big(f(x)\cdot\overline{ a_j(x)}\big)a_j(x),
$$
hence $\widetilde{\bf H}_j = \{\alpha(x)a_j(x)\, ; \, \alpha(x) \in L^2({\bf T}^d,{\bf C},dx)\}$.
One must note, however,  that the eigenvector cannot be chosen uniquely.  Therefore, we introduce an equivalent relation $\sim $ for $(\alpha(x), a_j(x)), (\beta(x), b_j(x))$, where 
$\alpha(x), \beta(x) \in  L^2({\bf T}^d,{\bf C},dx)$ and $a_j(x), b_j(x)$ are normalized eigenvectors for $H_0(x)$ associated with eigenvalue $\lambda_j(x)$,
$$
(\alpha(x),a_j(x)) \sim (\beta(x),b_j(x)) \Longleftrightarrow 
(\alpha(x)\cdot\overline{a_j(x)})a_j(x) = (\beta(x)\cdot\overline{b_j(x)})b_j(x).
$$
Then $\widetilde{\bf H}_j$ is the resulting equivalence class. Noting that
$$
dx = \frac{dM_{\lambda,j}d\lambda}{|\nabla \lambda_j(x)|},
$$
we identify $\widetilde{\bf H}_j$ with
\begin{equation}
{\bf H}_j^{(0)} = L^2(I_j^{(0)},{\bf h}_{\lambda,j}a_j,d\lambda),
\label{C4S3DefinebfHj}
\end{equation}
\begin{equation}
I_j^{(0)} = \lambda_j({\bf T}^d)\setminus{\mathcal T}_0,
\nonumber
\end{equation}
 where ${\bf h}_{\lambda,j}a_j = \{\alpha a_j\big|_{M_{\lambda,j}}\, ; \, \alpha \in {\bf h}_{\lambda,j}\}$, $a_j(x)$ being the normalized eigenvector of $H_0(x)$.
Finally, letting
$$
I^{(0)} = {\mathop\cup_{j=1}^s} I_j^{(0)} = \sigma(H_0)\setminus{\mathcal T_0},
$$ 
we define ${\bf h}_{\lambda,j}$ and ${\bf H}_j^{(0)}$  to be $\{0\}$ for $\lambda \in I^{(0)} \setminus I_j^{(0)}$, and put
\begin{equation}
{\bf h}_{\lambda} = {\bf h}_{\lambda,1}a_1\oplus \cdots \oplus {\bf h}_{\lambda,s}a_s, 
\end{equation}
\begin{equation}
{\bf H}^{(0)} = {\bf H}_1^{(0)} \oplus \cdots \oplus {\bf H}_s^{(0)} = L^2(I^{(0)},{\bf h}_{\lambda},d\lambda).
\label{C4S4SpacebfhbfH}
\end{equation}
In the above definition of ${\bf H}_j^{(0)}$, ${\bf h}_{\lambda,j}$ depends also on $\lambda$. Since $M_{\lambda,j}$ is defined by $\lambda_j(x) = \lambda$, by splitting the interval $I_j^{(0)}$ into subintervals $I_{j,\ell}^{(0)}$ so that on each $I_{j,\ell}^{(0)}$, $M_{\lambda,j}$ is diffeomorphic to a fixed manifold $\mathcal N_{j,\ell}$, then  
${\bf H}_{j}^{(0)}$ should be written as a direct sum $\oplus_{\ell}L^2(I_{j,\ell}^{(0)}, L^2(\mathcal N_{j,\ell}), d\mu_{j,\ell}(\lambda))$ for a suitable measure 
$d\mu_{j,\ell}(\lambda)$. For the sake of simplicity, however, we use the definition (\ref{C4S3DefinebfHj}).

\medskip
By Lemma \ref{C4S4Im(upm,f)=|f|}, we have
$$
P_j\big(R_0(\lambda + i0) - R_0(\lambda - i0)\big) = \big( (\lambda_j(x) - \lambda - i0)^{-1} - (\lambda_j(x)- \lambda + i0)^{-1}\big)P_j.
$$
 In particular, we have the following Parseval formula : 
\begin{equation}
 \frac{1}{2\pi i}\big((R_0(\lambda + i0) - R_0(\lambda - i0))f,g\big) \\
= \sum_{j=1}^s\int_{M_{\lambda,j}}(P_jf)(x)\cdot\overline{(P_jg)(x)}\, \frac{dM_{\lambda,j}}{|\nabla \lambda_j(x)|}.
\label{C4S5ParesevalWholespace}
\end{equation}
Using the distribution
\begin{equation}
\mathcal E'({\bf R}^d) \ni \delta_{M_{\lambda,j}} :  \varphi \to \langle \delta_{M_{\lambda,j}},\varphi\rangle = \int_{M_{\lambda,j}}\varphi(x)\frac{dM_{\lambda,j}}{|\nabla\lambda_j(x)|},
\label{S6DefinedeltaMlambdaj}
\end{equation}
we can rewrite (\ref{C4S5ParesevalWholespace}) as
\begin{equation}
\frac{1}{2\pi i}\left(R_0(\lambda + i0) - R_0(\lambda - i0)\right) = 
\sum_{j=1}^s\delta_{M_{\lambda,j}}P_j.
\end{equation}

With the above formulas in mind, we define
\begin{equation}
\mathcal F_0 : {\mathcal H}_0 \ni f \to (P_1f,\cdots,P_sf)\in {\bf H}^{(0)}.
\end{equation}


\begin{theorem}\label{C4S4F0unitarydiagonalize}
(1) $\ \mathcal F_0$ is a unitary operator from $(L^2({\bf T}^d))^s$ to ${\bf H}^{(0)}$. \\
\noindent
(2) \ For any $f \in {\mathcal H}_0$, $(\mathcal F_0H_0f)(\lambda) = \lambda (\mathcal F_0f)(\lambda)$. 
\end{theorem}

\medskip
Let us further study
\begin{equation}
\mathcal F_0(\lambda)f = \left(\mathcal F_0f\right)(\lambda) = \big( \mathcal F_{0,1}(\lambda)f, \cdots, \mathcal F_{0,s}(\lambda)f\big)\in {\bf h}_{\lambda}.
\end{equation}
More precisely,
\begin{equation}
\mathcal F_{0,j}(\lambda)f  = 
\left\{
\begin{split}
& P_j(x)f(x)\Big|_{M_{\lambda,j}}, \quad {\rm if} \quad \lambda \in I_j^{(0)}, \\
&0, \quad {\rm otherwise},
\end{split}
\right.
\end{equation}
\begin{equation}
\mathcal F_{0,j}(\lambda)^{\ast}\phi = 
\left\{
\begin{split}
& P_j(x)\phi(x)\otimes\delta_{M_{\lambda,j}}, \quad {\rm if} \quad \lambda \in I_j^{(0)}, \\
& 0, \quad {\rm otherwise}.
\end{split}
\right.
\label{C4S5F0lambdaast}
\end{equation}
Then, by (\ref{C4S5ParesevalWholespace})
\begin{equation}
\frac{1}{2\pi i}\left((R_0(\lambda + i0) - R_0(\lambda - i0))f,g\right) = 
\left(\mathcal F_0(\lambda)f,\mathcal F_0(\lambda)g\right)_{{\bf h}_{\lambda}}.
\label{C4S3E07lambda=f0lambda}
\end{equation}
This implies
\begin{equation}
\|\mathcal F_0(\lambda)f\|^2_{\bf h_{\lambda}} = \|f\|^2_{L^2(M_{\lambda})} = 
\sum_{j=1}^s\|P_jf\|^2_{L^2(M_{\lambda,j})}.
\end{equation}


\begin{lemma}\label{C4S5F0lambdaestmateabove}
For any compact set $J \subset \sigma(H_0) \setminus \mathcal T_0$, there exists a constant $C > 0$ such that
$$
\|\mathcal F_0(\lambda)f\|_{{\bf h}_{\lambda}} \leq C\|f\|_{\mathcal B}, \quad \lambda \in J.
$$
\end{lemma}

Proof. This follows from (\ref{C4S5ParesevalWholespace}) and Theorem \ref {C4S4TheoremR0lambda)}. \qed

\medskip
As a direct consequence of this lemma, we have
$$
\mathcal F_0(\lambda)^{\ast} \in {\bf B}({\bf h}_{\lambda}\, ; \, \mathcal B^{\ast}).
$$
By Theorem \ref{C4S4F0unitarydiagonalize} (2), we have $\mathcal F_0(\lambda)(H_0 -\lambda)f=0$, $\forall f \in \mathcal B$. Therefore
\begin{equation}
(H_0-\lambda)\mathcal F_0(\lambda)^{\ast}=0.
\label{C4S4F0lambdastaeigenequation}
\end{equation}
Multiplying the cofactor matrix from the left, we then obtain
$$
p(x,\lambda)\mathcal F_0(\lambda)^{\ast} = 0,
$$
which implies that ${\rm supp}\,\mathcal F_0(\lambda)^{\ast}\phi \subset M_{\lambda}$ for any $\phi  \in {\bf h}_{\lambda}$.


\begin{lemma}\label{C4S5F0ast1to1}
For any compact set $J \subset \sigma(H_0) \setminus \mathcal T_0$, 
there exists a constant $C > 0$ such that
$$
C^{-1}\|\phi\|_{{\bf h}_{\lambda}} \leq  \|\mathcal F_0(\lambda)^{\ast}\phi\|_{\mathcal B^{\ast}} \leq C\|\phi\|_{{\bf h}_{\lambda}}, \quad \lambda \in J.
$$
\end{lemma}

Proof. The estimate from above follows from Lemma \ref{C4S5F0lambdaestmateabove}. To prove the estimate from below, we have only to take $\phi$ in a dense set of ${\bf h}_{\lambda}$.
Therefore, without loss of generality, we can assume that $T(x)^{-1}\phi$ is continuous where $ T(x) := \oplus_{j=1}^s | \nabla \lambda_j (x)|$. 
We take a partition of unity $\{\chi_j\}$ on $M_{\lambda}$ so that $\phi_j = \chi_j\phi$ has a small support. 
 Let $u_j  = \mathcal F_0(\lambda)^{\ast}\phi_j$ on the fundamental domain, and extend it to be 0 outside. Then, $u_j \in \mathcal B^{\ast}({\bf R}^d)$ and one can apply Lemma \ref{C3S3u0estimatedbylinsup}
to get the estimate 
$$
\|\phi_j\|_{{\bf h}_{\lambda}} \leq  C\|\mathcal F_0(\lambda)^{\ast}\phi_j\|_{\mathcal B^{\ast}}.
$$
We can take a smooth extension $\widetilde \chi_j \in C^{\infty}_0({\bf R}^d)$ of $\chi_j$ so that $\mathcal F_0(\lambda)^{\ast}\phi_j = \widetilde{\chi_j}\mathcal F_0(\lambda)^{\ast}\phi$.
Then we have
$$
\|\widetilde{\chi_j}\mathcal F_0(\lambda)^{\ast}\phi\|_{\mathcal B^{\ast}} \leq C\|\mathcal F_0(\lambda)^{\ast}\phi\|_{\mathcal B^{\ast}},
$$
hence
$$
\|\phi\|_{{\bf h}_{\lambda}} \leq \sum_j \|\phi_j\|_{{\bf h}_{\lambda}} \leq \sum_jC\|\widetilde{\chi_j}\mathcal F_0(\lambda)^{\ast}\phi\|_{\mathcal B^{\ast}} \leq C\|\mathcal F_0(\lambda)^{\ast}\phi\|_{\mathcal B^{\ast}},
$$
which completes the proof of the lemma.  \qed


\begin{lemma}\label{C4S5OntoTrace}
 For $\lambda \in \sigma(H_0) \setminus \mathcal T_0$, $\mathcal F_0(\lambda){\mathcal B} = {\bf h}_{\lambda}$.
\end{lemma}

Proof. We use the following Banach's closed range theorem (see \cite{Yosida66}, p. 205). 


\begin{theorem}\label{C4S5Banachclosedrange}
Let $X, Y$ be Banach spaces, and $T$ a bounded operator from $X$ to $Y$. Let $R(T) = \{Rx\, ; \, x \in X\}$, $N(T) = \{x \in X\, ; \,  Tx=0\}$, and denote the pairing between $X$ and its dual space $X^{\ast}$ by $\langle\, , \rangle$. Then the following 4 assertions are equivalent. \\
\noindent
(1) $\ R(T)$ is closed. \\
\noindent
(2) $\ R(T^{\ast})$ is closed. \\
\noindent
(3) $\ R(T) = N(T^{\ast})^{\perp} = \{y \in Y\, ; \, \langle y,y^{\ast}\rangle = 0, \ \forall y^{\ast} \in N(T^{\ast})\}$.\\
\noindent
(4) $\ R(T^{\ast}) = N(T)^{\perp} = \{x^{\ast} \in X^{\ast}\, ; \, \langle x,x^{\ast}\rangle = 0, \ \forall x \in N(T)\}$.
\end{theorem}

\medskip
We take $X = {\mathcal B}, Y = {\bf h}_{\lambda}$ and $T = \mathcal F_0(\lambda)$. Then by 
Lemma \ref{C4S5F0ast1to1}, $T^{\ast}$ is 1 to 1, and $R(T^{\ast})$ is closed.
 Theorem \ref{C4S5Banachclosedrange} then implies that $R(T)$ is dense and closed. 
$\qquad \qed$


\begin{lemma}\label{C4S4N(lambda)=F0lambdaasth}
For $\lambda \in \sigma(H_0) \setminus \mathcal T_0$,
$\{u \in \mathcal B^{\ast}\, ; \, (H_0-\lambda)u=0\} =  \mathcal F_0(\lambda)^{\ast}{\bf h}_{\lambda}$.
\end{lemma}

Proof. The inclusion relation $\{u \in \mathcal B^{\ast}\, ; \, (H_0-\lambda)u=0\} \supset  \mathcal F_0(\lambda)^{\ast}{\bf h}_{\lambda}$ is proven in (\ref{C4S4F0lambdastaeigenequation}). To prove the converse relation, in view of (4) of Theorem \ref{C4S5Banachclosedrange}, we have only to prove
$$
u\in \mathcal B^{\ast},\ f \in \mathcal B, \ (H_0-\lambda)u=0, \ \mathcal F_0(\lambda)f =0 \Longrightarrow (u,f)=0.
$$
However, as has been seen before Lemma \ref{C4S5F0ast1to1}, $(H_0-\lambda)u=0$ implies ${\rm supp}\, u \subset M_{\lambda}$. By Lemma \ref{C3S3u0estimatedbylinsup}, $u$ is an $L^2$-density on $M_{\lambda}$.  Then
$(u,f)$ is an integral of $u(x)\overline{f(x)}$ on $M_{\lambda}$, However, the restriction of $f(x)$ to $M_{\lambda}$ is 0. This proves the lemma. $\qquad \qed$

\medskip
Let us also mention the singularity expansion of the 
solution of the Helmholtz equation $(H_0-\lambda)u=0$. Let 
\begin{equation}
A_{\pm}(\lambda) = \frac{1}{2\pi i}\sum_{j=1}^s
\dfrac{1}{\lambda_j(x) - \lambda \mp i0}\otimes P_j(x)\Big|_{x\in M_{\lambda,j}}.
\label{C4S4DefineApmlambda}
\end{equation}
Then by (\ref{C4S5F0lambdaast}),
\begin{equation}
\mathcal F_0(\lambda)^{\ast} = A_+(\lambda) - A_-(\lambda).
\label{C4S4F0lambda=A+-A-}
\end{equation}
In view of Lemma \ref{C4S4N(lambda)=F0lambdaasth}, we obtain the following lemma. 


\begin{lemma}
Let $\lambda \in \sigma(H_0)\setminus{\mathcal T}_0$. Then for any solution $u \in \mathcal B^{\ast}$ of the equation $(H_0-\lambda)u=0$, there exists 
$\phi \in {\bf h}_{\lambda}$ such that
$$
u = A_+(\lambda) \phi - A_-(\lambda)\phi.
$$
\end{lemma}


\section{Hamiltonian on the perturbed lattice}


\subsection{Perturbed lattice}
Let $\Gamma_0 = \{\mathcal L_0, \mathcal V_0, \mathcal E_0\}$ be the periodic graph given in Subsection \ref{subsectionperiodicgraph}. To {\it remove an edge} $e$ from $\Gamma_0$, we mean to remove $e$ from the edge set $\mathcal E_0$, but not to remove end points $o(e)$ and $t(e)$ of $e$ from the vertex set $\mathcal V_0$. 
To {\it remove a vertex} $v$ from $\mathcal V_0$, we mean to remove $v$ from $\mathcal V_0$ as well as all edges having $v$ as an end point.
We can also perturb a lattice by adding vertexes and edges.

Now let us perturb a finite part of $\Gamma_0$, and denote the resulting graph by $\Gamma = \{\mathcal V, \mathcal E\}$. 
Take a large {\bf integer} $a > 0$, and put
\begin{eqnarray}
{\bf Z}^d_{ext} &= & {\bf Z}^d\setminus\{n\, ; \, |n_i| \leq a , 1\leq i \leq d\}, \\
\partial{\bf Z}^d_{ext} &=& {\mathop\cup_{i=1}^d}\{n \in {\bf Z}^d_{ext}\, ; \, |n_i|=a\}, \\
\mathcal V_{ext} &=& {\mathop \cup_{j=1}^s}\{p_j + {\bf v}(n)\, ; \, n \in {\bf Z}^d_{ext}\}, \\
\partial \mathcal V_{ext} &=& {\mathop \cup_{j=1}^s}\{p_j + {\bf v}(n)\, ; \, n \in \partial{\bf Z}^d_{ext}\}, \\
\stackrel{\circ}{\mathcal V_{ext} } &=& \mathcal V_{ext}\setminus\partial\mathcal V_{ext}, \\
\mathcal V_{int} &= & \mathcal V\setminus \stackrel{\circ}{\mathcal V_{ext}}, \\
\partial\mathcal V_{int} &=&\partial\mathcal V_{ext}, \\
\stackrel{\circ}{\mathcal V_{int}} &=& \mathcal V_{int}\setminus\partial\mathcal V_{int}.
\end{eqnarray}
Then, $\mathcal V$ consists of a disjoint union of two parts :
$$
\mathcal V = \stackrel{\circ}{\mathcal V_{ext}} \cup\, \mathcal V_{int}, \quad \; {\sharp}\mathcal V_{int} < \infty.
$$
Accordingly, the Hilbert space $\ell^2({\mathcal V})$ admits an orthogonal decomposition
$$
\ell^2(\mathcal V) = \ell^2(\stackrel{\circ}{\mathcal V_{ext}}) \oplus \ell^2(\mathcal V_{int}).
$$
The elements $\widehat u$ of $\ell^2(\stackrel{\circ}{{\mathcal V}_{ext}})$ are written as vectors of $s$-components : $\widehat u(n) = (\widehat u_1(n),\cdots,\widehat u_s(n)), \ \widehat u_i(n) \in \ell^2({\bf Z}^d_{ext})$, while $\widehat w \in \ell^2(\mathcal V_{int})$ is simply  a finite dimensional vector.

Let 
\begin{equation}
 \widehat P_{ext} : \ell^2(\mathcal V) \to \ell^2(\stackrel{\circ}{\mathcal V_{ext}}), \quad
 \widehat P_{int} : \ell^2(\mathcal V) \to \ell^2(\mathcal V_{int})
\nonumber
\end{equation}
be the associated orthogonal projections, which are naturally extended to $\ell^2_{loc}(\mathcal V)$.
Let $\widehat \Delta_{\Gamma}$ be the Laplacian on the graph $\Gamma$, which is self-adjoint on $\ell^2({\mathcal V})$.


\subsection{Spectral properties for $\widehat H$}
We put 
\begin{equation}
\widehat H = - \widehat\Delta_{\Gamma} + \widehat V,
\end{equation}
and study its spectral properties. Let us repeat our assumptions (A-1) $\sim$ (A-4) 
on the unperturbed lattice and add new assumptions (A-5), (A-6) on perturbations.

\medskip
\noindent
{\bf (A-1)}  {\it There exists a subset $\mathcal T_1 \subset \sigma(H_0)$ such that for $\lambda \in \sigma(H_0)\setminus\mathcal T_1$ : }

\smallskip
{\bf (A-1-1)}  $M_{\lambda,sng}^{\bf C}$ {\it  is discrete.}

\smallskip
{\bf (A-1-2)} {\it Each connected component of $M_{\lambda,reg}^{\bf C}$ intersects with ${\bf T}^d$ and the intersection is a $(d-1)$-dimensional real analytic submanifold of ${\bf T}^d$.}

\medskip
\noindent
{\bf (A-2)} 
{\it There exists a finite set $\mathcal T_0 \subset \sigma(H_0)$ such that }
$$
\ 
M_{\lambda, i} \cap M_{\lambda,j} = \emptyset, \quad { if} \quad i \neq j, \quad \lambda \in \sigma(H_0)\setminus\mathcal T_0.
$$

\medskip
\noindent
{\bf (A-3)} $\ 
\nabla_xp(x,\lambda) \neq 0, \quad {on} \quad M_{\lambda}, \quad \lambda \in \sigma(H_0)\setminus\mathcal T_0.$

\medskip
\noindent
{\bf (A-4)}   {\it The unique continuation property holds for $\widehat H_0$ in $\mathcal V_0$.}

\medskip
\noindent
{\bf (A-5)} \ {\it $\widehat V$ is bounded self-adjoint on $\ell^2(\mathcal V)$ and has  support in $\mathcal V_{int}$, i.e. $\widehat V\widehat u = 0$ on  $\mathcal V_{ext}$, $\forall \widehat u \in \ell^2(\mathcal V)$.}

\medskip
\noindent
{\bf (A-6)}   {\it The unique continuation property holds for $\widehat H_0$ in ${\mathcal V_{ext}}$.}

\medskip
Here, (A-6) is defined in the same way as in (A-4) with $\mathcal V_0$ replaced by $\mathcal V_{ext}$.

\bigskip
The spaces $\widehat{\mathcal B}$, $\widehat{\mathcal B}^{\ast}$, $\widehat{\mathcal B}^{\ast}_0$, $\ell^{2,t}$ are naturally defined on $\mathcal V$.
For $z \not\in {\bf R}$, let
\begin{equation}
\widehat R(z) = (\widehat H- z)^{-1}, 
\end{equation}
\begin{equation}
\widehat R_0(z) = (\widehat H_0 - z)^{-1}, 
\end{equation}
where $\widehat H_0$ is the Hamiltonian on $\mathcal V_0$ defined in Section 4. 
Given a subset $S \subset {\bf Z}^d$ and its characteristic function $\widehat\chi_S$, we use $\widehat\chi_S$ to denote either the operator of restriction to $S$ : $\ell^2_{loc}({\bf Z}^d) \ni \widehat u \to \widehat u\big|_S \in \ell^2_{loc}(S)$, or the operator of extension  : $\ell^2_{loc}(S) \ni \widehat u \to \widehat v \in \ell^2_{loc}({\bf Z}^d)$, where $\widehat v = \widehat u$ on $S$, $\widehat v = 0$ on ${\bf Z}^d\setminus S$. This will not confuse our arguments.
 The spaces $\ell^2(\mathcal V)$ and $\ell^2(\mathcal V_0)$, on which $\widehat H$ and $\widehat H_0$ live, differ only by a finite dimensional space. Therefore by the well-known theorem on the compact perturbation of self-adjoint operators, we have the following theorem.


\begin{theorem}
$\ \sigma_{e}(\widehat H) = \sigma(\widehat H_0)$. 
\end{theorem}


\begin{lemma}\label{Lemma7.2}
(1)  The eigenvalues of $\widehat H$ in $\sigma_{e}(\widehat H)\setminus \mathcal T_1$ is finite with finite multiplicities. \\
\noindent
(2)  There is no eigenvalue in $\sigma_{e}(\widehat H)\setminus\mathcal T_1$, provided ${\mathcal V_{int}}$ has the unique continuation property.
\end{lemma}

Here, the unique continuation property on $\mathcal V_{int}$ is the following assertion.
{\it Suppose for some $\lambda \in {\bf C}$, $(- \widehat\Delta_{\Gamma}+ \widehat V - \lambda)\widehat u =0$ holds on $\mathcal V$ and $\widehat u(n) = 0$ on $\stackrel{\circ}{\mathcal V_{ext}}$. 
Then $\widehat u =0$ on whole ${\mathcal V}$.}

\medskip

Proof. Let $\widehat u$ be the eigenvector of $\widehat H$ with eigenvalue in $\sigma(\widehat H_0)\setminus\mathcal T_1$. 
Then, by Theorem \ref{C4S2RellichTh},
 it vanishes near infinity. By (A-6), we then have $\widehat u = 0$ on ${\mathcal V_{ext}}$. 
Therefore, all eigenvectors are supported in $\mathcal V_{int}$, hence are finite-dimensional. The assertion (2) is obvious. \qed

\bigskip
 So far, we have studied the operator $\widehat H$ under the assumptions (A-1) $\sim$ (A-6). It is because we are interested in the Rellich-Vekua theorem and the absence of embedded eigenvalues, both of which play important roles in the application to the inverse problem. However, by adopting the  already established perturbation method in scattering theory, we can study the spectral properties of $\widehat H$, including $\mathcal T_1$,  under the assumptions (A-2), (A-3) and (A-5)  only. The trade-off is the weaker result for embedded eigenvalues.

For $\widehat u \in \widehat{\mathcal B}^{\ast}$, its wave front set is defined by
$$
WF^{\ast}(\widehat u) = WF^{\ast}(\mathcal U\widehat P_{ext}\widehat u).
$$
A solution $\widehat u \in \widehat{\mathcal B}^{\ast}$ of the equation
$$
(\widehat H - \lambda)\widehat u = \widehat f \in \mathcal B
$$
is said to be outgoing (or incoming), if $\mathcal U\widehat P_{ext}\widehat u$ satisfies  the condition (\ref{C4S4OutInRadCond}).

Recall the following lemma.


\begin{lemma}\label{Hardy}
For $f(x) \in L^1(0,\infty)$, put
$$
u(x) = \int_x^{\infty}f(t)dt.
$$
Then, for any $s > 1/2$, we have
$$
\int_0^{\infty}x^{2(s-1)}|u(x)|^2dx \leq \frac{4}{(2s-1)^2}\int_0^{\infty}
x^{2s}|f(x)|^2dx.
$$
\end{lemma}

This is well-known, and can be proven by using Hardy's inequality. For the proof see e.g. 
\cite{IsKu}, Chapter 3, Lemma 3.3.


\begin{lemma}\label{C5S3RadCondUnique}
Assume (A-2), (A-3), (A-5).  Then, for any compact interval $I \subset \sigma(\widehat H_0)\setminus\mathcal T_0$ and $s> 0$, there exists a constant $C_{s,I} > 0$ such that if $\widehat u \in \mathcal B^{\ast}$ satisfies $(\widehat H - \lambda)\widehat u =0$ on $\mathcal V$ and the outgoing (or incoming) radiation condition,
$$
\|\widehat u\|_{\ell^{2,s}} \leq C_{s,I}\|\widehat u\|_{\mathcal B^{\ast}}, \quad \forall \lambda 
\in  I.
$$
\end{lemma}

Proof. We put $\widehat u_e = \widehat P_{ext}\widehat u, \widehat u_i = \widehat P_{int}\widehat u$. First we show that
\begin{equation}
{\rm Im}((\widehat H - \lambda)\widehat u_e,\widehat u_e)) = 0.
\label{C5S3ImHuu=0}
\end{equation}
In fact, using the equation, we have
\begin{equation}
\begin{split}
& ((\widehat H - \lambda)\widehat u_e,\widehat u_e)) \\
= & - ((\widehat H - \lambda)\widehat P_{int}\widehat u,\widehat P_{int}\widehat u))  -((\widehat H - \lambda)\widehat P_{ext} \widehat u, \widehat P_{int} \widehat u)
 - ((\widehat H - \lambda)\widehat P_{int} \widehat u, \widehat P_{ext}\widehat u).
\end{split}
\nonumber
\end{equation}
By the computation
\begin{equation}
\begin{split}
&  ((\widehat H - \lambda)\widehat P_{int}\widehat u,\widehat P_{ext}\widehat u)) - 
 (\widehat P_{ext}\widehat u, (\widehat H - \lambda)\widehat P_{int}\widehat u) \\
& =  ((\widehat H - \lambda)\widehat P_{int}\widehat u,\widehat u - \widehat P_{int}\widehat u)) - 
 (\widehat u - \widehat P_{int}\widehat u, (\widehat H - \lambda)\widehat P_{int}\widehat u) \\
& =  ((\widehat H - \lambda)\widehat P_{int}\widehat u,\widehat P_{int}\widehat u)) - 
 (\widehat P_{int}\widehat u, (\widehat H - \lambda)\widehat P_{int}\widehat P_{int}\widehat u) \\
& =  (\widehat  P_{int}(\widehat H - \lambda)\widehat P_{int}\widehat u,\widehat P_{int}\widehat u)) - 
 (\widehat P_{int}\widehat u, \widehat P_{int}(\widehat H - \lambda)\widehat P_{int}\widehat P_{int}\widehat u)  \\
& = 0,
\end{split}
\nonumber
\end{equation}
we prove (\ref{C5S3ImHuu=0}).

We put $\widehat f_e = (\widehat H - \lambda)\widehat u_e = [\widehat H,\widehat P_{ext}]\widehat u$, which is compactly supported.  
Letting $\mathcal U\widehat u_e = u_e$ and $\mathcal U\widehat f_e = f_e$, we then have $(H_0 - \lambda)u_e = f_e$, and (\ref{C5S3ImHuu=0}) implies ${\rm Im}(u_e,f_e) = 0$. Lemma \ref{C4S4Im(upm,f)=|f|} then yields $f_e\big|_{M_{\lambda}}=0$.
By  Lemma \ref{C4S4RadcondinBast0}, it follows that $u_e \in {\mathcal B}^{\ast}_0$. 

Now, we argue as in the proof of Lemma \ref{C3S3u+lemma}. We take $y_1 = p(x,\lambda)$ as a new variable and pass to the Fourier transform. Then taking account of (\ref{Lemma4.3u+intf(xi)}), since $u_e$ satisfies the outgoing radiation condition, we have
$$
\|\widetilde u_e(\xi_1,\cdot)\|_t \leq C_t\int_{\xi_1}^{\infty}\|\widetilde f_e(\eta_1,\cdot)\|_td\eta_1, \quad \forall t > 0.
$$
Here, $\|\cdot\|_t$ denotes the $L^2({\bf R}^{d-1})$-norm with weight $(1 + |\xi'|^2)^{t/2}$. 
By Lemma \ref{Hardy}, we have
$$
\int_0^{\infty}\xi_1^{2(s-1)}\|\widetilde u_e(\xi_1,\cdot)\|^2_t d\xi_1\leq
C_s\int_0^{\infty}\xi_1^{2s}\|\widetilde f_e(\xi_1,\cdot)\|^2_td\xi_1, \quad 
\forall s > 1/2, \ \forall t > 0.
$$
Since $u_e \in \mathcal B_0^{\ast}$, 
$\widetilde u_e$ also satisfies the incoming radiation condition. Therefore we have, similarly
$$
\int_{-\infty}^0|\xi_1|^{2(s-1)}\|\widetilde u_e(\xi_1,\cdot)\|^2_t d\xi_1\leq
C_s\int_{-\infty}|\xi_1|^{2s}\|\widetilde f_e(\xi_1,\cdot)\|^2_td\xi_1, \quad 
\forall s > 1/2, \ \forall t > 0.
$$
Define the norm $\|\cdot\|_{s,t}$ by
$$
\|u\|_{s,t} = \|(1 + |\xi_1|^2)^{s/2}(1 + |\xi'|^2)^{t/2}\widetilde u(\xi_1,\xi')\|_{L^2}.
$$
Then, the above two inequalities imply
$$
\| u_e\|_{s-1,t} \leq C_{s,t}\| f_e\|_{s,t}, \quad 
\forall s > 1/2, \ \forall t > 0.
$$
Since $\widehat f_e$ is compactly supported, we have
$$
\|f_e\|_{s,t} \leq C_{s,t}\|u\|_{\mathcal B^{\ast}},
$$
and the lemma follows. \qed

\bigskip

 Let $\sigma_{rad}(\widehat H)$ 
be the set of $\lambda \in \sigma(\widehat H_0)\setminus{\mathcal T}_0$ for which there exists $0 \neq \widehat u \in \widehat{\mathcal B}^{\ast}$ satisfying $(\widehat H - \lambda)\widehat u =0$ on $\mathcal V$ and the outgoing radiation condition or incoming radiation condition. 


\begin{lemma}\label{mathcalEdiscrete}
Assume (A-2), (A-3), (A-5).  \\
\noindent
(1) $\ \sigma_{rad}(\widehat H) = \sigma_p(\widehat H)\cap \big(\sigma(\widehat H_0)\setminus {\mathcal T}_0\big)$. \\
(2) $\ \sigma_p(\widehat H)\cap \big(\sigma(\widehat H_0)\setminus {\mathcal T}_0\big)$ is discrete in $\sigma(\widehat H_0)\setminus \mathcal T_0$ with possible accumulation points in $\mathcal T_0$, 
and the multiplicity of each eigenvalue in $\sigma_p(\widehat H)\cap \big(\sigma(\widehat H_0)\setminus {\mathcal T}_0\big)$ is finite. Moreover, the associated eigenvector belongs to $\ell^{2,s}$, $\forall s > 0$.
\end{lemma}

Proof. The assertion (1) follows from Lemma \ref{C5S3RadCondUnique}.  If (2) is not true, 
there exists an infinite number of eigenvalues $\lambda_n \in \sigma(\widehat H_0)\setminus \mathcal T_0$, counting multiplicity,  converging to an interior point $\lambda \in \sigma(\widehat H_0)\setminus \mathcal T_0$. Let $\{\widehat u_n\}_{n=1}^{\infty}$ be the associated orthonormal system of eigenvectors. 
Lemma \ref{C5S3RadCondUnique} then implie that $\sup_n\|\widehat u_n\|_s < \infty$ for all $s > 0$.
One can then select a subsequence $\{\widehat u_n'\}$ converging strongly in $L^2(\mathcal V)$, which leads to a contradiction. The last assertion is proven in Lemma \ref{C5S3RadCondUnique}. \qed

\medskip
As a consequence, we have obtained the following uniqueness theorem.


\begin{lemma}\label{LemmaRadCondUnique}
Let $\lambda \in \sigma_e(\widehat H)\setminus\big(\mathcal T_0\cup \sigma_p(\widehat H)\big)$, and suppose $\widehat u \in \widehat{\mathcal B}^{\ast}$ satisfies $(\widehat H - \lambda)\widehat u=0$. Then $\widehat u=0$, if $\widehat u$ satisfies the outgoing or incoming radiation condition. 
In particular, if we assume (A-1) and $\lambda \not\in \mathcal T_0\cup\mathcal T_1$, the solution of $(\widehat H - \lambda)\widehat u = \widehat f$ satisfying the outgoing or incoming radiation condition is unique provided $\widehat H$ has the unique continuation property..
\end{lemma}

Proof. The 1st assertion follows from Lemma \ref{mathcalEdiscrete} (1). To prove the 2nd assertion, we assume that $\widehat f=0$. Then, $\lambda \in \sigma_{rad}(\widehat H)$, hence $\lambda \in \sigma_p(\widehat H)$, and $\widehat u \in \ell^2$. By Theorem \ref{C4S2RellichTh}, $\widehat u$ vanishes near infinity, hence in $\mathcal V$ by the unique continuation property. \qed

\bigskip
We put
\begin{equation}
\mathcal T = \mathcal T_0 \cup\sigma_p(\widehat H).
\end{equation}
\medskip
Let us introduce
\begin{equation}
\widehat Q_1(z) = (\widehat H_0 - z)\widehat P_{ext}\widehat R(z) = \widehat P_{ext} + \widehat K_1\widehat R(z),
\label{S7DefineQ1}
\end{equation}
\begin{equation}
\widehat K_1 = \widehat H_0\widehat P_{ext} - \widehat P_{ext}\widehat H,
\label{S7DefineK1}
\end{equation}
\begin{equation}
\widehat Q_2(z) = (\widehat H - z)\widehat P_{ext}\widehat R_0(z) = \widehat P_{ext} + \widehat K_2\widehat R_0(z),
\label{S7DefineK2}
\end{equation}
\begin{equation}
\widehat K_2 = \widehat H\widehat P_{ext} - \widehat P_{ext}\widehat H_0,
\label{DefineK2}
\end{equation}
Then, we have
\begin{equation}
\widehat P_{ext}\widehat R(z) = \widehat R_0(z)\widehat Q_1(z),
\label{S7PextRz=R0zQ1}
\end{equation}
\begin{equation}
\widehat P_{ext}\widehat R_0(z) = \widehat R(z)\widehat Q_2(z).
\end{equation}

We shall now derive the limiting absorption principle for $\widehat H$.


\begin{theorem}\label{C5S3LAPforwidehatH}
Assume (A-2), (A-3) and (A-5). \\
\noindent
(1) For any compact set $J \subset \sigma_{e}(\widehat H)\setminus\mathcal T$, there is a constant $C > 0$ such that
$$
\|\widehat R(z)\widehat f\|_{\widehat{\mathcal B}^{\ast}} \leq C\|\widehat f\|_{\widehat{\mathcal B}}, \quad 
{\rm Re}\,z \in J, \quad {\rm Im}\ z \neq 0.
$$
(2) For any $t > 1/2$, there exists a strong limit $\lim_{\epsilon\to 0}\widehat R(\lambda \pm   i\epsilon)\widehat f \in \ell^{2,-t}$ for $\widehat f \in \widehat{\mathcal B}$ and $\lambda \in J$. Moreover, $\widehat R(\lambda \pm i0)\widehat f \in \widehat{\mathcal B}^{\ast}$, 
$\lim_{\epsilon \to 0}(\widehat R(\lambda\pm i\epsilon)\widehat f,\widehat g) = (\widehat R(\lambda\pm 0)\widehat f, \widehat g)$ for any $\widehat f, \widehat g \in \widehat{\mathcal B}$, and the inequality
$$
\|\widehat R(\lambda\pm i0)\widehat f\|_{\widehat{\mathcal B}^{\ast}} \leq C\|\widehat f\|_{\widehat{\mathcal B}}, \quad \lambda \in J
$$
holds. \\
(3) For any $\widehat f, \widehat g \in \widehat{\mathcal B}$, 
$$
\sigma_{e}(\widehat H)\setminus \mathcal T \ni \lambda \to \widehat R(\lambda \pm i0)\widehat f \in \ell^{2,-t}, \quad t > 1/2,
$$
$$
\sigma_{e}(\widehat H)\setminus \mathcal T \ni \lambda \to (\widehat R(\lambda \pm i0)\widehat f,\widehat g)
$$
is continuous. \\
\noindent
(4) For $\widehat f \in \widehat{\mathcal B}$, $\widehat R(\lambda \pm i0)\widehat f$ satisfies the radiation condition. Moreover, letting $u_{\pm} = \mathcal U\widehat P_{ext}\widehat R(\lambda \pm i0)\widehat f$, and 
\begin{equation}
 Q_{1}(\lambda \pm i0) = \,\mathcal U\widehat Q_1(\lambda \pm i0),
\end{equation}
we have
\begin{equation}
P_ju_{\pm} \mp \frac{1}{\lambda_j(x) - \lambda \mp i0}\otimes
\big(P_jQ_{1}(\lambda \pm i0)\widehat f\big)\Big|_{M_{\lambda,j}}\in \mathcal B^{\ast}_0,
\end{equation}
where $P_j = P_j(x)$ is the eigenprojection associated with the eigenvalue $\lambda_j(x)$ of $H_0(x)$.
\end{theorem}

Proof. 
Letting $\widehat u(z) = \widehat R(z)\widehat f$, we have the following inequality
\begin{equation}
\|\widehat u(z)\|_{\widehat{\mathcal B}^{\ast}} \leq C\left(\|\widehat P_{ext}\widehat f\|_{\widehat{\mathcal B}}+ \|\widehat K_1\widehat u(z)\|_{\widehat{\mathcal B}} + \| \widehat P_{int}\widehat u(z)\|_{\widehat{\mathcal B}}\right).
\label{C5S3PerurbInequal}
\end{equation}
Note that $\widehat K_1$ and $\widehat P_{int}$ are finite rank operators. Moreover, we have
\begin{equation}
\widehat R(z) = \widehat P_{ext}\widehat R_0(z)\widehat P_{ext} + 
\widehat P_{ext}\widehat R_0(z)\widehat K_1 \widehat R(z) + \widehat P_{int}\widehat R(z).
\label{C5S3Rzresolvent}
\end{equation}

Suppose (1) does not hold. Then, there exist $\widehat f_j \in \widehat{\mathcal B}$, $z_j$ such that, letting $\widehat u_j = \widehat R(z_j)\widehat f_j$, we have, 
${\rm Re}\, z_j \in J$, 
$$
\|\widehat f_j\|_{\widehat{\mathcal B}} \to 0, \quad \|\widehat u_j\|_{\widehat{\mathcal B}^{\ast}} = 1.
$$
We can assume without loss of generality that $z_j \to \lambda + i0 \in J$. Take $t > 1/2$. Then the embedding $\widehat{\mathcal B}^{\ast} \subset \ell^{2,-t}$ is compact. Therefore, 
there exists a subsequence of $\{\widehat u_j\}$, which is again denoted by $\{\widehat u_j\}$, such that $\widehat u_j \to \widehat u$ in $\ell^{2,-t}$. In view of (\ref{C5S3PerurbInequal}), we have $\widehat u_j \to \widehat u$ in $\widehat{\mathcal B}^{\ast}$, hence $\|\widehat u\|_{\widehat{\mathcal B}^{\ast}} = 1$. Since $(\widehat H - z_j)\widehat u_j = \widehat f_j$, we have $(\widehat H -\lambda) \widehat u = 0$. Moreover, by (\ref{C5S3Rzresolvent}), 
$$
\widehat u = \widehat P_{ext}\widehat R_0(\lambda + i0)\widehat K_1 \widehat u + \widehat P_{int} \widehat u.
$$
Therefore, $\widehat u$ is outgoing. By Lemma \ref{LemmaRadCondUnique}, 
we have $\widehat u = 0$, which is a contradiction.

Next, let $\epsilon_j \to 0$ and $\widehat f \in \widehat{\mathcal B}$.
Then (1) and the compact embedding $\widehat {\mathcal B}^* \subset \ell^{2,-t}$ imply that there exists a subsequence $\epsilon_{j_p} \to 0$ and $\widehat v \in \ell^{2,-t}$ such that $\widehat v_p := \widehat R(\lambda + i\epsilon_{j_p}) \widehat f \to \widehat v$ in $\ell^{2,-t}$ as $p \to \infty$.
Therefore, $\widehat v$ satisfies $(\widehat H_0 + \widehat V - \lambda) \widehat v = \widehat f$.
In view of (\ref{C5S3PerurbInequal}), we have
\begin{equation}
  \widehat v = \widehat P_{ext} \widehat R_0(\lambda + i0) \widehat P_{ext} \widehat f + \widehat P_{ext} \widehat R_0(\lambda + i0) \widehat K_1 \widehat v + \widehat P_{int} \widehat v.
  \label{S7Limit_u}
\end{equation}
The third term on the right hand side is in $\widehat{\mathcal B}^*$, since $\widehat P_{int}$ is a finite rank operator.
$\widehat K_1$ is also a finite rank operator, which implies $\widehat K_1 \widehat v \in \widehat{\mathcal B}$; $\widehat f \in \widehat{\mathcal B}$ by the hypothesis.
Therefore, the first and the second terms in the right-hand side are also in $\widehat{\mathcal B}^*$, which means that the left hand side $\widehat v \in \widehat{\mathcal B}^*$.

Let us show the sequence $\{ \widehat R(\lambda \pm i\epsilon_j) \widehat f\}_{j=1}^{\infty}$ itself converges to $\widehat v$ in $\ell^{2,-t}$.
Assuming the contrary, we have another subsequence $\widehat v_{j_q} := \widehat R(\lambda + i\epsilon_{j_q}) \widehat f$ which satisfies
$$
  \| \widehat v - \widehat v_{j_q} \|_{\ell^{2,-t}} \ge \gamma, \quad (q = 1, 2, \cdots),
$$
for some $\gamma > 0$.
Then we can find a subsequence which is again denoted by $\{ \widehat v_{j_q}\}_{q=1}^{\infty}$ and $\widehat v^{\prime} \in \ell^{2,-t}$ such that $\widehat v_{j_{q}} \to \widehat v^{\prime}$ in $\ell^{2,-t}$ as $q \to \infty$,  $(\widehat H_0 + \widehat V - \lambda) \widehat v^{\prime} = \widehat f$, and
\begin{equation}
  \| \widehat v - \widehat v^{\prime} \|_{\ell^{2,-t}} \ge \gamma > 0.
  \label{S7Subsequence}
\end{equation}
In the same way as above, we have
\begin{equation}
  \widehat v^{\prime} = \widehat P_{ext} \widehat R_0(\lambda + i0) \widehat P_{ext} \widehat f + \widehat P_{ext} \widehat R_0(\lambda + i0) \widehat K_1 \widehat v^{\prime} + \widehat P_{int} \widehat v^{\prime},
  \label{S7Limit_u_prime}
\end{equation}
and $\widehat v^{\prime} \in \widehat {\mathcal B}^*$.

Subtracting (\ref{S7Limit_u_prime}) from (\ref{S7Limit_u}), we have
$$
\widehat v - \widehat v^{\prime} = \widehat P_{ext} \widehat R_0(\lambda + i0) \widehat K_1 (\widehat v - \widehat v^{\prime}) + \widehat P_{int} (\widehat v - \widehat v^{\prime}),
$$
which means $\widehat v - \widehat v^{\prime}$ is outgoing, and $(\widehat H - \lambda)(\widehat v - \widehat v^{\prime}) = 0$.
By Lemma \ref{LemmaRadCondUnique}, 
$\widehat v = \widehat v^{\prime}$, which contradicts (\ref{S7Subsequence}).

We define $\widehat R(\lambda + i0) \widehat f := \widehat v$. Let $\widehat f,\, \widehat g \in l^{2,t} \subset \widehat{\mathcal B}$.
Then $\widehat R(\lambda + i\epsilon) \widehat f \to \widehat R(\lambda + i0)\widehat f$ in $\ell^{2,-t}$ implies $(\widehat R(\lambda + i\epsilon) \widehat f, \widehat g) \to (\widehat R(\lambda + i0) \widehat f, \widehat g)$ as $\epsilon \to 0$, and
$$
|(\widehat R(\lambda + i0) \widehat f, \widehat g)| = \lim_{\epsilon \to 0}|(\widehat R(\lambda + i\epsilon) \widehat f, \widehat g)| \le C \|\widehat f\|_{\mathcal B}\,\|\widehat g\|_{\mathcal B}
$$
where $C>$ is a constant in (1).
Such a weak $\ast$ limit also exists for any $\widehat f,\,\widehat g \in \widehat{\mathcal B}$ and the inequality in (2) holds, since $\ell^{2,t}$ is dense in  $\widehat{\mathcal B}$.
We have thus proven (2).

As for (3), we can show the continuity of both mappings using the inequality in (2) instead of that in (1), the compact embedding $\widehat{\mathcal B}^* \subset \ell^{2,-t}$, and the denseness of $\ell^{2,t}$ in $\widehat{\mathcal B}$ in the same way as in (2).

By (\ref{S7PextRz=R0zQ1}), we have
for $\widehat f \in \widehat{\mathcal B}$
$$
\mathcal U \widehat P_{ext} \widehat R(\lambda \pm i0) \widehat f = \mathcal U \widehat R_0(\lambda \pm i0) \widehat Q_1(\lambda \pm i0) \widehat f.
$$
Applying Theorem \ref{C4S4TheoremR0lambda)} (4), we have the assertion (4). \qed


\subsection{Spectral representation}
We keep assuming (A-2), (A-3), (A-5). 
It is convenient to put for $\lambda \in \sigma_e(\widehat H)\setminus{\mathcal T}$
\begin{equation}
\widehat E'(\lambda) = \frac{1}{2\pi i}\left(\widehat R(\lambda + i0) - \widehat R(\lambda - i0)\right),
\label{S7DefineEprime0lambda}
\end{equation}
\begin{equation}
\widehat E_0'(\lambda) = \frac{1}{2\pi i}\left(\widehat R_0(\lambda + i0) - \widehat R_0(\lambda - i0)\right).
\label{DefineEprimelambda}
\end{equation}


\begin{lemma}\label{C5S3E'lambdaandE0'lambda}
We have  for $\lambda \in \sigma_e(\widehat H)\setminus{\mathcal T}$
\begin{equation}
(\widehat E_0'(\lambda)\widehat Q_1(\lambda\pm i0)\widehat f,\widehat Q_1(\lambda \pm i0)\widehat f) = (\widehat E'(\lambda)\widehat f,\widehat f),
\label{Lemma7.6Eprimelambda}
\end{equation}
\begin{equation}
(\widehat E'(\lambda)\widehat Q_2(\lambda\pm i0)\widehat f,\widehat Q_2(\lambda \pm i0)\widehat f) = (\widehat E_0'(\lambda)\widehat f,\widehat f).
\label{Lemma7.6Eprimelambda2}
\end{equation}
\end{lemma}

Proof. 
Recalling (\ref{S7DefineQ1}) and letting $z = \lambda + i\epsilon$, we have
\begin{equation}
\begin{split}
 \left((\widehat R_0(z)-\widehat R_0(\overline{z}))\widehat Q_1(z)\widehat f,
\widehat Q_1(z)\widehat f\right) 
= & 2i\epsilon \left(\widehat R_0(z)\widehat Q_1(z)\widehat f,\widehat R_0(z)\widehat Q_1(z)\widehat f\right) \\
= & 2i\epsilon \left(\widehat P_{ext}\widehat R(z)\widehat f,\widehat P_{ext}\widehat R(z)\widehat f\right) \\
= & 2i\epsilon \left(\widehat R(z)\widehat f,\widehat R(z)\widehat f\right) - 2i\epsilon \left(\widehat P_{int}\widehat R(z)\widehat f,\widehat R(z)\widehat f\right)  \\
= & \left((\widehat R(z)-\widehat R(\overline{z}))\widehat f,\widehat f\right) 
- 2i\epsilon \left(\widehat P_{int}\widehat R(z)\widehat f,\widehat R(z)\widehat f\right).
\end{split}
\nonumber
\end{equation}
Letting $\epsilon \to 0$, we obtain
$$
(\widehat E_0'(\lambda)\widehat Q_1(\lambda\pm i0)\widehat f,\widehat Q_1(\lambda \pm i0)\widehat f) = (\widehat E'(\lambda)\widehat f,\widehat f).
$$

Similarly, we obtain (\ref{Lemma7.6Eprimelambda2}). \qed

\bigskip
We define
\begin{equation}
\widehat{\mathcal F}_0(\lambda) = \mathcal F_0(\lambda)\mathcal U,
\end{equation}
\begin{equation}
\widehat{\mathcal F}_{\pm}(\lambda) = \widehat{\mathcal F}_0(\lambda)\widehat Q_{1}(\lambda \pm i0).
\label{C5S3DefineFpmlambdaFoQlambda}
\end{equation}
The formula (\ref{C4S3E07lambda=f0lambda}) implies $E_0'(\lambda) = \mathcal F_0(\lambda)^{\ast}\mathcal F_0(\lambda)$. Then in view of (\ref{Lemma7.6Eprimelambda}) and (\ref{C5S3DefineFpmlambdaFoQlambda}), we have the following lemma.


\begin{lemma}\label{C5S3Fpmlambdaparseval}
For $\lambda \in \sigma_{e}(\widehat  H)\setminus \mathcal T$, and $\widehat f, \widehat g \in \widehat{\mathcal B}$,
\begin{equation}
\big(\widehat E'(\lambda)\widehat f,\widehat g\big) = 
\big(\widehat{\mathcal F}_{\pm}(\lambda)\widehat f,\widehat{\mathcal F}_{\pm}(\lambda)\widehat g\big)_{{\bf h}_{\lambda}}.
\tag{\ref{Lemma7.6Eprimelambda}$'$}
\end{equation}
\end{lemma}

Theorem \ref{C5S3LAPforwidehatH} (2) then yields


\begin{lemma}
For any compact set $J \subset \sigma_{e}(\widehat H)\setminus\mathcal T$, there exists a constant $C >0$ such that
$$
\|\widehat{\mathcal F}_{\pm}(\lambda)\widehat f\|_{{\bf h}_{\lambda}} \leq C\|\widehat f\|_{\widehat{\mathcal B}}, \quad 
\lambda \in J.
$$
\end{lemma}

\bigskip
We are now in a position to apply the stationary scattering theory  for Schr{\"o}dinger operators (see e.g. \cite{Is}, and also the recent articles \cite{IsKu}, \cite{Yaf}). By following this abstract framework, one can derive the spectral representation (Theorem \ref{C5S3FpmSpecRepre}), the asymptotic completeness of wave operators (Theorem \ref{S7WaveOp}) and the unitarity of the S-matrix (Theorem \ref{S7Smatrix}). Since this is a well-known already established argument,  we only give the outline of the proof.  Let 
\begin{equation}
I = \sigma_e(H)\setminus{\mathcal T}.
\end{equation}
 Take any borel set $e \in I$, and integrate (\ref{Lemma7.6Eprimelambda}$'$) on $e$. Then letting $\widehat E(\lambda)$ be the spectral decomposition of $\widehat H$, we have
\begin{equation}
(\widehat E(e)\widehat f,\widehat g) = \int_e\big(\widehat{\mathcal F}_{\pm}(\lambda)\widehat f,\widehat{\mathcal F}_{\pm}(\lambda)\widehat g\big)_{{\bf h}_{\lambda}}d\lambda.
\label{S7Eprimeintegrate}
\end{equation}
We put
\begin{equation}
\begin{split}
I_j &= \lambda_j({\bf T}^d)\setminus \mathcal T, \\
I &= {\mathop\cup_{j=1}^s}I_j = \sigma(H_0)\setminus\mathcal T, \\
{\bf H}_j &= L^2(I_j,{\bf h}_{\lambda,j}a_j,d\lambda), \\
{\bf H} &= {\bf H}_1 \oplus \cdots \oplus {\bf H}_s = L^2(I,{\bf h}_{\lambda}, d\lambda).
\end{split}
\nonumber
\end{equation}
 Note that ${\bf H} = {\bf H}^{(0)}$, since $H$ is absolutely continuous on $I$ by virtue of Theorem \ref{C5S3LAPforwidehatH}.
We define $\widehat{\mathcal F}_{\pm}\widehat f \in {\bf H}$  by
$$
\big(\widehat{\mathcal F}_{\pm}\widehat f\big)(\lambda) = \widehat{\mathcal F}_{\pm}(\lambda)\widehat f.
$$
Then, by (\ref{S7Eprimeintegrate}), $\widehat{\mathcal F}_{\pm}$ is uniquely extended to an isometry from initial set $\widehat E(I){\ell}^2({\mathcal V})$ and the final set contained in ${\bf H}$.
Theorem \ref{C5S3LAPforwidehatH} implies $\widehat E(I){\ell}^2({\mathcal V}) = \mathcal H_{ac}(\widehat H)$.
By (\ref{Lemma7.6Eprimelambda2}), one derives
$$
(\widehat{\mathcal F}_{\pm}(\lambda)\widehat Q_2(\lambda\pm i0)\widehat f, \widehat{\mathcal F}_{\pm}(\lambda)\widehat Q_2(\lambda\pm i0)\widehat g) = 
(\widehat{ E}'_0(\lambda)\widehat f,\widehat g).
$$
Making use of this formula, one can construct a similar partial isometry $\widehat{\mathcal A}_{\pm}$ satisfying $\widehat{\mathcal F}_{\pm}\widehat{\mathcal A}_{\pm} = \widehat{\mathcal F}_0$. Therefore, the final set of $\widehat{\mathcal F}_{\pm}$ is equal to ${\bf H}$.
By (\ref{S7DefineQ1}) and (\ref{C5S3DefineFpmlambdaFoQlambda}), one can show
$$
\widehat{\mathcal F}_{\pm}(\lambda)(\widehat H - \lambda)\widehat f = 0, \quad 
\widehat f \in \widehat{\mathcal B},
$$
which implies $\widehat{\mathcal  F}_{\pm}\widehat H = \lambda\widehat{\mathcal F}_{\pm}$. In summary, one has


\begin{theorem}\label{C5S3FpmSpecRepre}
(1) $\widehat{\mathcal F}_{\pm}$ is a partial isometry with initial set $\widehat E(I){\ell}^2({\mathcal V})$ and final set $\bf H$. \\
\noindent
(2) $\big(\widehat{\mathcal F}_{\pm}\widehat H\widehat f\big)(\lambda) = \lambda \big(\widehat{\mathcal F}_{\pm}\widehat f\big)(\lambda), \quad \lambda \in \sigma_{e}(\widehat H)\setminus\mathcal T, \quad \widehat f \in \mathcal H $. \\
\noindent
(3) For $\lambda \in \sigma_{e}(\widehat H)\setminus{\mathcal T}$ and  $\phi \in {\bf h}_{\lambda}$,
$\big(\widehat H - \lambda\big)\widehat{\mathcal F}_{\pm}(\lambda)^{\ast}\phi = 0.$
\end{theorem}

Let us consider the time-dependent wave operators. The following theorem is well-known. See e.g. \cite{KatoKuroda} and \cite{Kuroda}.


\begin{theorem}\label{S7WaveOp}
(1) There exists a strong limit
\begin{equation}
{\mathop{\rm s-lim}_{t \to \pm \infty}}\,e^{it\widehat H}\widehat P_{ext}e^{-it\widehat H_0}\widehat P_{ac}(\widehat H_0) := \widehat W_{\pm},
\nonumber
\end{equation}
where $\widehat P_{ac}(\widehat H_0)$ is the projection onto the absolutely continuous subspace for $\widehat H_0$. \\
\noindent
(2) For any bounded Borel function $\psi$ on ${\bf R}$, we have
$$
\psi(\widehat H)\widehat W_{\pm} = \widehat W_{\pm}\psi(\widehat H_0).
$$
(3) $\widehat W_{\pm}$ is a partial isometry with initial set $\mathcal H_{ac}({\widehat H}_0)$ and final set
$\mathcal H_{ac}(\widehat H)$. 
\end{theorem}

The scattering operator $\widehat S$ is defined by
\begin{equation}
\widehat S = \big(\widehat W_+\big)^{\ast}\widehat W_-.
\label{S7ScatteringOpDefine}
\end{equation}
By Theorem \ref{S7WaveOp} (3), $\widehat S$ is unitary on ${\bf H}$.
We derive an expression of  $S = \widehat{\mathcal F}_0\widehat S(\widehat{\mathcal F}_0)^{\ast}$, where $\widehat{\mathcal F}_0 = \mathcal F_0\,\mathcal U$. 
 Since this is a well-known fact, we explain the procedure formally.
Assume that $\widehat f, \widehat g \in E_0(I)(L^2({\bf T}^d))^s$. Then, we have
\begin{equation}
\begin{split}
((\widehat S -1)\widehat f,\widehat g) &= 
((\widehat W_- - \widehat W_+)\widehat f, \widehat W_+\widehat g)\\
&= -i\int_{-\infty}^{\infty}(e^{it\widehat H}\widehat K_2e^{-it\widehat H_0}\widehat f,\widehat W_+\widehat g)dt \\
&= -i\int_{-\infty}^{\infty}(\widehat K_2e^{-it\widehat H_0}\widehat f,e^{-it\widehat H}\widehat W_+\widehat g)dt \\
&= -i\int_{-\infty}^{\infty}(\widehat K_2e^{-it\widehat H_0}\widehat f,\widehat P_{ext}e^{-it\widehat H_0}\widehat g)dt \\
&\ \ \ \ -\int_0^{\infty}ds\int_{-\infty}^{\infty}(\widehat K_2e^{-it\widehat H_0}\widehat f,e^{is\widehat H}\widehat K_2e^{-i(s+t)\widehat H_0}\widehat g)dt, 
\end{split}
\label{S7Formula0}
\end{equation}
where in the 3rd line, we have used $e^{-it\widehat H}\widehat W_+ = \widehat W_+e^{-it\widehat H_0}$ which follows from Theorem \ref{S7WaveOp} (2), and also
$$
\widehat W_+ = \widehat P_{ext}\widehat P_{ac}(\widehat H_0) + i
\int_0^{\infty}e^{is\widehat H_0}\widehat K_2e^{-is\widehat H_0}ds\widehat P_{ac}(\widehat H_0).
$$
Letting $I = \sigma(\widehat H_0)$, and passing to the spectral representation, we have
\begin{equation}
\begin{split}
& \int_{-\infty}^{\infty}(\widehat K_2 e^{-it\widehat H_0}\widehat f, e^{is\widehat H}\widehat K_2e^{-i(s+t)\widehat H_0}\widehat g)dt\\
&  = 
\int_{-\infty}^{\infty}\int_I(\widehat{\mathcal F}_0(\lambda)\widehat K_2^{\ast}e^{-is\widehat H}\widehat K_2e^{-it\widehat H_0}\widehat f,
e^{-i(s+t)\lambda}\widehat{\mathcal F}_0(\lambda)\widehat g)d\lambda dt .
\end{split}
\label{S7Formula1}
\end{equation}
using Theorem~\ref{S7WaveOp} (2),
Inserting $e^{-\epsilon|t|}$ and letting $\epsilon \to 0$, we see that this is equal to 
\begin{equation}
\begin{split}
& 2\pi\int_I(\widehat{\mathcal F}_0(\lambda)\widehat K_2^{\ast}e^{-i(\widehat H - \lambda)s}\widehat K_2\widehat E_0'(\lambda)\widehat f, \widehat{\mathcal F}_0(\lambda)\widehat g)d\lambda \\
& = 2\pi\int_I(\widehat{\mathcal F}_0(\lambda)\widehat K_2^{\ast}e^{-i(\widehat H - \lambda)s}\widehat K_2\widehat{\mathcal F}_0(\lambda)^{\ast}\widehat{\mathcal F}_0(\lambda)\widehat f, \widehat{\mathcal F}_0(\lambda)\widehat g)d\lambda. 
\end{split}
\label{S7Formula2}
\end{equation}
Therefore, the 2nd term of the most right-hand side of (\ref{S7Formula0}) is equal to 
$$
- 2\pi\int_0^{\infty}ds\int_I(\widehat{\mathcal F}_0(\lambda)\widehat K_2^{\ast}e^{-i(\widehat H - \lambda)s}\widehat K_2\widehat{\mathcal F}_0(\lambda)^{\ast}\widehat{\mathcal F}_0(\lambda)\widehat f, \widehat{\mathcal F}_0(\lambda)\widehat g)d\lambda.
$$
Inserting $e^{-\epsilon s}$, and letting $\epsilon \to 0$, this is equal to 
\begin{equation}
2\pi i \int_I(\widehat{\mathcal F}_0(\lambda)\widehat K_2^{\ast}\widehat R(\lambda + i0)\widehat K_2\widehat{\mathcal F}_0(\lambda)^{\ast}\widehat{\mathcal F}_0(\lambda)\widehat f, \widehat{\mathcal F}_0(\lambda)\widehat g)d\lambda.
\label{S7Formula3}
\end{equation}
Similarly, the 1st term of  the most right-hand side of (\ref{S7Formula0}) is equal to 
\begin{equation}
- 2\pi i \int_I(\widehat{\mathcal F}_0(\lambda)\widehat P_{ext}\widehat K_2\widehat{\mathcal F}_0(\lambda)^{\ast}\widehat{\mathcal F}_0(\lambda)\widehat f, \widehat{\mathcal F}_0(\lambda)\widehat g)d\lambda.
\label{S7Formula4}
\end{equation}
Summing up, we obtain
$$
((\widehat S -1)\widehat f,\widehat g) = - 2\pi i\int_I(A(\lambda)\widehat{\mathcal F}_0(\lambda)\widehat f,\widehat{\mathcal F}_0(\lambda)\widehat g)d\lambda,
$$
where
\begin{equation}
A(\lambda) = \widehat{\mathcal F}_0(\lambda)\widehat P_{ext}\widehat K_2\widehat{\mathcal F}_0(\lambda)^{\ast} - \widehat{\mathcal F}_0(\lambda)\widehat K_2^{\ast}\widehat R(\lambda + i0)\widehat K_2\widehat{\mathcal F}_0(\lambda)^{\ast}.
\nonumber
\end{equation}
Using (\ref{S7DefineQ1}), (\ref{DefineK2}) and (\ref{C5S3DefineFpmlambdaFoQlambda}), we arrive at
\begin{equation}
A(\lambda) = \widehat{\mathcal F}_0(\lambda)\widehat Q_1(\lambda+ i0)\widehat K_2\widehat{\mathcal F}_0(\lambda)^{\ast} = \widehat{\mathcal F}_+(\lambda)\widehat K_2\widehat{\mathcal F}_0(\lambda)^{\ast}.
\label{S7Alambda}
\end{equation}
The S-matrix is now defined by
\begin{equation}
S(\lambda) = 1 - 2\pi i A(\lambda), \quad \lambda \in \sigma_e(\widehat H)\setminus{\mathcal T},
\label{S7Slambda}
\end{equation}


\begin{theorem}\label{S7Smatrix}
(1) For $f \in {\bf H}$, we have
\begin{equation}
(Sf)(\lambda) = S(\lambda)f(\lambda),  \quad \lambda \in 
\sigma_e(\widehat H)\setminus{\mathcal T},
\label{S7DirectintofS}
\end{equation}
(2) $S(\lambda)$ is unitary on ${\bf h}_{\lambda}$, $\lambda \in \sigma(\widehat H_0)\setminus\mathcal T$.
\end{theorem}

Proof. The assertion (1) is proven above. Since $S$ is unitary on ${\bf H}$, $S(\lambda)$ is unitary for a.e. $\lambda \in  \sigma_e(\widehat H)\setminus{\mathcal T}$. However, $S(\lambda)$ is strongly continuous for $\lambda \in \sigma_e(\widehat H)\setminus{\mathcal T}$. Hence  $S(\lambda)$ unitary for $\lambda \in \sigma_e(\widehat H)\setminus{\mathcal T}$. \qed


\subsection{Helmholtz equation}

The last topic is the characterization of the solution space of the Helmholtz equation $ \{{\widehat u} \in \widehat{\mathcal B}^{\ast}\, ; \, 
(\widehat H - \lambda)\widehat u = 0\}$
in terms of spectral representation. Letting
\begin{equation}
\varphi_j= \big(P_j\phi\big)\big|_{M_{\lambda,j}},
\end{equation}
we have using (\ref{C5S3DefineFpmlambdaFoQlambda}) and (\ref{C4S5F0lambdaast}),
\begin{equation} 
 \widehat{\mathcal F}_{-}(\lambda)^{\ast}\phi =  \widehat P_{ext}{\mathcal U}^{\ast}
\sum_{j=1}^s\delta_{M_{\lambda,j}}\otimes\varphi_j 
 + \widehat R(\lambda + i0)\widehat K_1^{\ast}\widehat{\mathcal F}_0(\lambda)^{\ast}\phi. 
\end{equation}
Noting that
$$
\widehat R(\lambda + i0) \equiv \widehat P_{ext}\widehat R_0(\lambda + i0)\widehat Q_1(\lambda + i0) 
$$
modulo a regular term, and also
$$
\delta_{M_{\lambda,j}} = 
\frac{1}{2\pi i}\Big(\frac{1}{\lambda_j(x)-\lambda - i0} - \frac{1}{\lambda_j(x)-\lambda + i0} \Big),
$$
we then have
\begin{equation}
\begin{split}
\mathcal U\widehat P_{ext} \widehat{\mathcal F}_{-}(\lambda)^{\ast}\phi & \equiv \sum_{j=1}^s
\frac{1}{2\pi i}\Big(\frac{1}{\lambda_j(x)-\lambda - i0} - \frac{1}{\lambda_j(x)-\lambda + i0} \Big)\otimes\varphi_j \\
& + \sum_{j=1}^s\frac{1}{\lambda_j(x)-\lambda - i0}P_j \, \mathcal U\widehat Q_1(\lambda + i0) \, \widehat K_1^{\ast}\widehat{\mathcal F}_0(\lambda)^{\ast}\phi.
\end{split}
\end{equation}
This is rewritten as

\begin{equation}
\begin{split}
 & \mathcal U \widehat P_{ext} \widehat{\mathcal F}_{-}(\lambda)^{\ast}\phi \\
 \equiv&  - \frac{1}{2\pi i}\sum_{j=1}^s\left(\frac{1}{\lambda_j(x) - \lambda + i0}\otimes\varphi_j  -\frac{1}{\lambda_j(x) - \lambda- i0}\otimes\varphi^{out}_{j}\right),\\
\varphi^{out}_{j} &= P_j\phi + P_j\mathcal U\widehat Q_1(\lambda + i0)\widehat K_1^{\ast}\widehat{\mathcal F}_0(\lambda)^{\ast}\phi.
\end{split}
\end{equation}
We let
\begin{equation}
\varphi^{in} = \phi,\quad
\varphi^{out} = \phi + \mathcal U\widehat Q_1(\lambda + i0)\widehat K_1^{\ast}\widehat{\mathcal F}_0(\lambda)^{\ast}\phi.
\label{BeforeLemma712Definephiiniout}
\end{equation}
A direct computation using (\ref{C5S3DefineFpmlambdaFoQlambda}), (\ref{S7DefineK1}), (\ref{S7DefineK2}) and (\ref{BeforeLemma712Definephiiniout}) entails
\begin{equation}
\varphi^{out} = S(\lambda)\varphi^{in}.
\label{S7varphioutSvarphiin}
\end{equation}

\begin{lemma}
For $\lambda \in \sigma_{e}(\widehat H)\setminus{\mathcal T}$, 
$\widehat{\mathcal F}_{\pm}(\lambda){\widehat{\mathcal B}} = {\bf h}_{\lambda}$.
\end{lemma}

Proof. We prove this lemma for $\widehat{\mathcal F}_{-}(\lambda)$. Let $u = \mathcal U\widehat P_{ext}\widehat {\mathcal F}_{-}(\lambda)^{\ast}\phi$ and $u_j = P_ju$. Take $x^0 \in M_{\lambda,j}$.  Without loss of generality we can make a change of variables $x \to y = (y_1,y')$ around $x^0$ such that $y_1=\lambda_j(x)-\lambda$. Let $\chi \in C^{\infty}({\bf T}^d)$ of the form $\chi(y) = \chi_1(y_1)\chi_2(y')$ such that $\chi  = 1$ on a small neighborhood of $x^0$ and the support of $\chi$ is sufficiently small.  Then the Fourier transform of $\chi(y) u_j(y)$ becomes 
$$
(2\pi)^{-d/2}\widetilde{\chi u_j}(\xi) \equiv \int_{\xi_1}^{\infty}\widetilde\chi_1(\eta_1)d\eta_1a(\xi') - \int_{-\infty}^{\xi_1}\widetilde\chi_1(\eta_1)d\eta_1b(\xi')
$$
modulo a sufficiently regular term, where $a(\xi'), b(\xi')$ are Fourier transforms of $\chi_2(y')\varphi^{out}_{j}(y')$, $\chi_2(y')\varphi_{j}(y')$, respectively.  Then, we have
$$
\lim_{R\to\infty}\frac{1}{R}\int_{|\xi|<R}|\widetilde{\chi u_j}(\xi)|^2d\xi 
\geq C(\|\chi_2(y')\varphi^{out}_{j}(y')\|^2_{L^2} + \|\chi_2(y')\varphi_{j}(y')\|^2_{L^2}).
$$
We take a finite number of such $x^0$ and sum up the above inequality. Then, the right-hand side is estimated from below by $C(\|\varphi^{out}_{j}\|^2 + 
\|\varphi_{j}\|^2) = C\|\phi\|_{M_{\lambda,j}}^2$. On the other hand, the left-hand side is estimated from above by $C\|\widehat{\mathcal F}_-(\lambda)^{\ast}\phi\|^2_{\widehat{\mathcal B}^{\ast}}$. We have thus proven
$$
\|\phi\|_{M_{\lambda}} \leq C\|\widehat{\mathcal F}_-(\lambda)^{\ast}\phi\|_{\widehat{\mathcal B}^{\ast}},
$$
hence the range of $\widehat{\mathcal F}_-(\lambda)^{\ast}$ is closed. We can then argue as in Lemma \ref{C4S5OntoTrace} to complete the proof. \qed

\medskip
We have now arrived at the main theorem of this paper. Recall the definition of $A_{\pm}(\lambda)$ in (\ref{C4S4DefineApmlambda}).


\begin{theorem}\label{S7FinalTheorem}
Assume (A-2), (A-3) and (A-5). 
Let $\lambda \in \sigma_{e}(\widehat H)\setminus{\mathcal T}$. \\
\noindent
(1) $\{{\widehat u} \in \widehat{\mathcal B}^{\ast}\, ; \, 
(\widehat H - \lambda)\widehat u = 0\} = 
\widehat{\mathcal F}_-(\lambda)^{\ast}{\bf h}_{\lambda}.$ \\
\noindent
(2) For any $\phi^{in} \in {\bf h}_{\lambda}$, there exist unique $\phi^{out} \in {\bf h}_{\lambda}$ and $\widehat u \in \widehat{\mathcal B}^{\ast}$ satisfying\begin{equation}
(\widehat H - \lambda)\widehat u = 0, 
\label{Th713Equation}
\end{equation}
\begin{equation}
\mathcal U\widehat P_{ext}\widehat u + A_-(\lambda)\phi^{in} 
- A_+(\lambda)\phi^{out} \in \mathcal B^{\ast}_0.
\label{Th713Bast0}
\end{equation}
Moreover,
$S(\lambda)\phi^{in} = \phi^{out}.$
\end{theorem}

Proof. To prove (1), as in Lemma \ref{C4S4N(lambda)=F0lambdaasth}, we have only to show
$$
{\widehat u} \in \mathcal {\widehat B}^{\ast}, \ {\widehat f} \in \widehat{\mathcal B}, \ (\widehat{H}-\lambda)\widehat{u}=0, \ 
\widehat{\mathcal F}_-(\lambda)\widehat{f}=0 \Longrightarrow (\widehat{u},\widehat{f})=0.
$$
Let $\widehat{v} = \widehat{R}(\lambda - i0)\widehat{f}$. Then $(\widehat{H}-\lambda)\widehat{v}=\widehat{f}$, hence $(\widehat{u},\widehat{f}) = (\widehat{u}, (\widehat{H}-\lambda)\widehat{v})$. 

Take $\chi_{\infty} \in C^{\infty}({\bf R}^1)$ such that 
$\chi_{\infty}(t)=1$ for $t > R_0$, where $R_0>0$ is sufficiently large, and 
$\chi_{\infty}(t) = 0$ for $t < R_0-1$. We put for large $r > 0$
$$
 \widehat P_{\infty,r} = \mathcal U^{-1}\chi_{\infty}(|N|/r)\mathcal U\widehat P_{ext}, \quad
\widehat P_{0,r} = (1 - \mathcal U^{-1}\chi_{\infty}(|N|/r)\mathcal U)\widehat P_{ext} + \widehat P_{int},
$$
where $|N|$ is defined by (\ref{S4Define|N|}).

First note
\begin{equation}
(\widehat u,(\widehat H - \lambda) \widehat P_{0,r}\widehat v) = 
((\widehat H - \lambda)\widehat u,\widehat P_{0,r} \widehat v)=0.
\label{Th713Proof1}
\end{equation}
 Take $\chi \in C^{\infty}_0({\bf R})$ such that $\chi(t)=1$ for $t=1$. Then using $(\widehat H-\lambda)\widehat u=0$, we have
\begin{equation}
(\widehat u,{\mathcal U}^{-1}\chi(|N|/R){\mathcal U}(\widehat H-\lambda){\widehat P}_{\infty,r}\widehat v) = ([{\widehat H},({\mathcal U}^{-1}\chi(|N|/R){\mathcal U})^{\ast}]{\widehat u},{\widehat P}_{\infty,r}\widehat v).
\label{S7Th714ProofCommu}
\end{equation}
Now, let $u = \mathcal U \widehat u, \ v = \mathcal U \widehat v$. By virtue of Theorem \ref{C5S3LAPforwidehatH} (4) and the assumption that $\widehat{\mathcal F}_-(\lambda)\widehat f=0$, we have $v \in \mathcal B^{\ast}_0$. Moreover, as is well-known, one can apply the standard micro-local calculus, or semi-classical calculus using $R^{-1}$ as Plank's constant, for the commutator $[{\widehat H},({\mathcal U}^{-1}\chi(|N|/R){\mathcal U})^{\ast}]$ (see e.g. Chap. 14 of \cite{Zw}). Therefore, the right-hand side of (\ref{S7Th714ProofCommu}) is dominated from above by
$$
\frac{C}{R}\left(\int_{|\xi|<CR}|\widetilde u(\xi)|^2d\xi\right)^{1/2}
\left(\int_{|\xi|<CR}|\widetilde v(\xi)|^2d\xi\right)^{1/2},
$$
which tends to 0 as $R\to\infty$. Therefore, we have
\begin{equation}
(\widehat u,(\widehat H - \lambda) \widehat P_{\infty,r}\widehat v) = 0.
\label{Th713Proof2}
\end{equation}
The equalities (\ref{Th713Proof1}), (\ref{Th713Proof2}) prove (1). 

The uniqueness assertion of (2) follows from Theorem \ref{C5S3RadCondUnique}. Because if $\widehat u_1, \widehat u_2$ are such solutions, $\widehat u_1 - \widehat u_2$ satisfies the outgoing radiation condition.

To show the existence, let 
$$
\widehat u^{in} = \widehat P_{ext}\mathcal U^{-1}\mathcal F_0(\lambda)^{\ast}\phi^{in}, \quad
\widehat v = - \widehat R(\lambda + i0)(\widehat H - \lambda)\widehat u^{in}.
$$
Since $(H_0-\lambda)\mathcal F_0(\lambda)^{\ast}\phi^{in}=0$, we have
$$
\widehat v = - \widehat R(\lambda +i0)\widehat K_2\mathcal U^{-1}\mathcal F_0(\lambda)^{\ast}\phi^{in}.
$$
Then we have by Theorem~\ref{C5S3LAPforwidehatH} (4)
\begin{equation}
P_j\mathcal U\widehat v \simeq
- \frac{1}{\lambda_j(x) - \lambda - i0)}\otimes
 (P_jQ_1(\lambda + i0){\widehat K}_2\mathcal U^{-1}\mathcal F_0(\lambda)^{\ast}\phi^{in}\big|_{M_{\lambda,j}}).
\end{equation}
We now let $\widehat u = \widehat u^{in} + \widehat v$. It is easy to see that 
$\widehat u \in \widehat B^{\ast}$ and $(\widehat H - \lambda)\widehat u = 0$.
Moreover,
\begin{equation}
\begin{split}
& P_j\mathcal U \widehat u \\
 \simeq& \frac{1}{2\pi i}\Big(\frac{1}{\lambda_j(x) - \lambda - i0}\otimes(P_j\phi^{in}\big|_{M_{\lambda, j}}) -\frac{1}{\lambda_j(x) - \lambda + i0}\otimes(P_j\phi^{in}\big|_{M_{\lambda, j}})\Big) \\
& \ \ \ \  -
\frac{1}{\lambda_j(x)-\lambda+i0}\otimes
\Big(P_jQ_1(\lambda + i0)\widehat K_2\mathcal U^{-1}\mathcal F_0(\lambda)^{\ast}\phi^{in}\big|_{M_{\lambda,j}}\Big) \\
= &-\frac{1}{2\pi i}\frac{1}{\lambda_j(x) - \lambda + i0}\otimes(P_j \phi^{in}\big|_{M_{\lambda, j}})\\
&   + \frac{1}{2\pi i}\frac{1}{\lambda_j(x)-\lambda -i0}\otimes\Big(P_j\phi^{in}\big|_{M_{\lambda,j}} - 2\pi i
P_jQ_1(\lambda + i0)\widehat K_2\mathcal U^{-1}\mathcal F_0(\lambda)^{\ast}\phi^{in}\big|_{M_{\lambda,j}}\Big).
\end{split}
\nonumber
\end{equation}
The assertion (2) then follows from this formula and  (\ref{S7Alambda}), (\ref{S7DefineQ1}), (\ref{S7DefineK2}). \qed

\end{document}